\documentclass[reqno]{amsart}
\usepackage{amsfonts}
\usepackage{amsmath}
\usepackage{amsthm, amscd}
\usepackage{mathtools}
\usepackage{amssymb}
\usepackage{mathrsfs}
\usepackage{graphicx}
\usepackage{caption}
\usepackage{subcaption}
\usepackage{float}
\usepackage{xypic}
\usepackage[abs]{overpic}
\usepackage[alphabetic]{amsrefs}
\usepackage{etex}

\usepackage{hyperref}
\hypersetup{
	colorlinks,
	linkcolor=[rgb]{0,0,0.7},
	urlcolor=[rgb]{0,0,0.4},
	citecolor=[rgb]{0.4,0.1,0}
}

\allowdisplaybreaks

\theoremstyle{plain}\newtheorem{Theorem}{Theorem}[section]
\theoremstyle{plain}\newtheorem{Corollary}[Theorem]{Corollary}
\theoremstyle{plain}\newtheorem{Lemma}[Theorem]{Lemma}
\theoremstyle{plain}\newtheorem{Definition}[Theorem]{Definition}
\theoremstyle{plain}\newtheorem{Proposition}[Theorem]{Proposition}
\theoremstyle{plain}
\theoremstyle{plain}
\theoremstyle{plain}\newtheorem*{Theorem*}{Theorem}
\theoremstyle{plain}

\theoremstyle{remark}\newtheorem{remark}[Theorem]{Remark}
\theoremstyle{remark}
\theoremstyle{remark}\newtheorem*{Notation*}{Notation}

\theoremstyle{plain}
\makeatletter
\newtheorem*{rep@theorem}{\rep@title}
\newcommand{\newreptheorem}[2]{
\newenvironment{rep#1}[1]{
 \def\rep@title{#2 \ref{##1}}
 \begin{rep@theorem}}
 {\end{rep@theorem}}}
\makeatother
\newreptheorem{Theorem}{Theorem}
\newreptheorem{Proposition}{Proposition}
\newreptheorem{Corollary}{Corollary}

\numberwithin{equation}{section}

\usepackage{lipsum}

\DeclareMathOperator{\parr}{par}

\DeclareMathOperator{\rank}{rank}
\DeclareMathOperator{\II}{I}
\DeclareMathOperator{\ad}{ad}

\DeclareMathOperator{\muu}{\mu^{orb}}

\DeclareMathOperator{\Sp}{Sp}
\DeclareMathOperator{\SU}{SU}

\DeclareMathOperator{\spann}{span}
\DeclareMathOperator{\Hom}{Hom}
\DeclareMathOperator{\sHom}{\mathcal{H}om}
\DeclareMathOperator{\pt}{pt}
\DeclareMathOperator{\id}{id}

\DeclareMathOperator{\Ext}{Ext}

\newcommand{\bA}{\mathbb{A}}
\newcommand{\bC}{\mathbb{C}}
\newcommand{\bE}{\mathbb{E}}

\newcommand{\bU}{\mathbb{U}}
\newcommand{\bV}{\mathbb{V}}
\newcommand{\bW}{\mathbb{W}}
\newcommand{\bZ}{\mathbb{Z}}

\newcommand{\cM}{\mathcal{M}}

\author{Yi Xie}
\address{Beijing International Center for Mathematical Research, Peking University, Beijing 100871, China}
\email{yixie@pku.edu.cn}
\author{Boyu Zhang}
\address{Department of Mathematics, The University of Maryland at College Park, Maryland 20742, USA}
\email{bzh@umd.edu}

\title{Ring structures in singular instanton homology}

\begin{document}


\begin{abstract}
We calculate the ring structure of the singular instanton Floer homology of $(S^1\times \Sigma, S^1\times \{p_1,\dots,p_n\})$ with $\bC$--coefficients, where $\Sigma$ is a closed oriented surface. As an application, we prove an excision formula for singular instanton homology when $n=1$. This settles the last unknown case of excision formula for instanton Floer homology.
\end{abstract}

\maketitle
\section{Introduction}
\label{sec_intro}

Let $\Sigma$ be a closed oriented surface with genus $g$, and let $n$ be an odd positive integer.  Let $p_1,\dots,p_n$ be $n$ distinct points on $\Sigma$. This paper studies the singular instanton Floer homology of the pair $(S^1\times \Sigma, S^1\times\{p_1,\dots,p_n\})$ with $\bC$--coefficients. The Floer homology group is denoted by
\begin{equation}
	\label{eqn_II(S1_times_Sigma)}
\II(S^1\times \Sigma, S^1\times\{p_1,\dots,p_n\}, \emptyset;\bC),
\end{equation}
and it admits a ring structure (see Section \ref{subsec_ring_structure}). The main purpose of this paper is to give a  complete algebraic characterization of this ring structure.

The study of the ring structure of \eqref{eqn_II(S1_times_Sigma)} and its variations has a long history which can be dated back to Atiyah--Bott \cite{AB}. The ring structure of the non-singular instanton Floer homology of $S^1\times \Sigma$ was computed by Mu\~noz \cite{Munoz}. The ring structure of \eqref{eqn_II(S1_times_Sigma)} when $g=0$ was completely computed by Street \cite{Street} (with $\bC$--coefficients) and Kronheimer--Mrowka \cite{kronheimer2022relations} (with local coefficient systems). In a previous work \cite{XZ:excision}, the authors partially computed the ring structure of  \eqref{eqn_II(S1_times_Sigma)} for $g>0$ and used it to obtain several results on annular Khovanov homology. 

The ring structures of instanton homology on $S^1\times \Sigma$  characterize all the ``universal relations'' of instanton homology on general three-manifolds. Therefore, they reflect fundamental features of instanton Floer theory.
Earlier studies of these ring structures have yielded significant results in low-dimensional topology. For example, the sutured instanton homology introduced by Kronheimer--Mrowka \cite{KM:suture}, which is now a fundamental tool in low-dimensional topology, was defined based on Mu\~noz's computation \cite{Munoz}. The partial computation of the ring structure of \eqref{eqn_II(S1_times_Sigma)} in \cite{XZ:excision} also yielded numerous topological applications, including new detection results on Khovanov homology \cite{XZ:forest, xie2020links, li2021two, xie2020instantons, baldwin2021khovanov, li2022instanton}, and results on the $\SU(2)$--representation varieties of link complements \cite{xie2023meridian}. 

This paper gives a complete computation of the ring structure of \eqref{eqn_II(S1_times_Sigma)}. In Theorem \ref{thm_uniqueness_algebraic_structure_+} below, we will present a list of algebraic properties of the ring  \eqref{eqn_II(S1_times_Sigma)} and show that these algebraic properties uniquely determine the ring structure.  The proof of Theorem \ref{thm_uniqueness_algebraic_structure_+} also yields an algorithm that computes the ring structure of  \eqref{eqn_II(S1_times_Sigma)} for every given pair $(g,n)$ (see Remark \ref{rmk_algorithm}).  We will explicitly carry out this algorithm for $n=1$ and write down a recursive formula for the generating relations. As an application of this computation, we prove an excision formula for singular instanton Floer homology for $n=1$, which resolves the last open case of the excision formula in instanton homology theory.
 
The proof of Theorem \ref{thm_uniqueness_algebraic_structure_+} is based on earlier results from the literature as well as two major new ingredients. The first ingredient (Theorem \ref{thm_generators_of_Ign}) is an algebraic description of the \emph{singular} cohomology ring of the moduli space of critical points of the unperturbed Chern--Simons functional, which may be of independent interest. It is known that the singular cohomology ring of this moduli space gives the leading-order terms of the relations in the ring structure of \eqref{eqn_II(S1_times_Sigma)} (Propositions \ref{prop_I_realized_by_J}, \ref{prop_J_realized_by_I}). 
The second ingredient (Proposition \ref{prop_subleading}) gives a formula for the \emph{sub-leading} terms of the non-trivial relations at the lowest degree.  The proof of Proposition \ref{prop_subleading} uses a wall-crossing argument following a strategy of Kronheimer--Mrowka \cite{kronheimer2022relations}. 
 
 Many of the results in this paper can be extended to instanton Floer homology with local coefficient systems, which will be discussed in Section \ref{sec_local_coef}. We will give a complete calculation of the ring structure with local coefficients when $n=1$. 
 
 The paper is organized as follows. Section \ref{sec_preliminaries} reviews earlier results in the literature that will be needed  in this paper. Section \ref{sec_parabolic} gives the algebraic description of the singular cohomology ring of the moduli space of critical points of the unperturbed Chern--Simons functional. Section \ref{sec_subleading} proves the formula of sub-leading terms. Section \ref{sec_uniqueness} then states and proves the main theorem, which characterizes the ring structure of \eqref{eqn_II(S1_times_Sigma)}. Section \ref{sec_n=1} computes the ring structure for $n=1$ and uses it to prove the excision formula. Section \ref{sec_local_coef} discusses some generalizations to Floer homology with local coefficient systems. 
 
 \subsection*{Acknowledgment}
 We would like to express our most sincere gratitude to Peter Kronheimer and Tom Mrowka for explaining their work to us. We also want to thank Junliang Shen for helpful conversations about the Mumford relations. 

\section{Preliminaries}
\label{sec_preliminaries}
Recall that we use $\Sigma$ to denote a closed oriented surface with genus $g$, use $n$ to denote a positive odd integer, and $p_1,\dots,p_n$ denote $n$ distinct points on $\Sigma$. 
Let
$$\bV_{g,n} :=\II(S^1\times \Sigma, S^1\times \{p_1,\dots,p_n\},\emptyset;\bC)$$ 
be the singular instanton Floer homology of the pair $(S^1\times \Sigma, S^1\times \{p_1,\dots,p_n\})$ with trivial $w_2$ in $\bC$--coefficients. For the definition of singular instanton Floer homology, the reader may refer to \cite{KM:Kh-unknot}. The notation here follows \cite{XZ:excision}*{Section 2}.

In this section, we review the known properties of $\bV_{g,n}$
from the literature. 
We will only cite the results that are needed later in the paper; for a more complete historical review on this subject,  we refer the reader to \cite{XZ:excision}*{Sections 1, 2}.
All the notation and conventions in this paper will be consistent with \cite{XZ:excision}.

\subsection{The ring structure}
\label{subsec_ring_structure}
We first review the definition of the ring structure on $\bV_{g,n}$. Let $F$ be an oriented pair-of-pants cobordism from $S^1\sqcup S^1$ to $S^1$. Then the triple
$$
(F\times \Sigma, F\times \{p_1,\dots,p_n\}, \emptyset)
$$
defines a cobordism from $(S^1\times \Sigma, S^1\times \{p_1,\dots,p_n\},\emptyset)\sqcup (S^1\times \Sigma, S^1\times \{p_1,\dots,p_n\},\emptyset)$ to $(S^1\times \Sigma, S^1\times \{p_1,\dots,p_n\},\emptyset)$. As a result, it defines a homomorphism 
\begin{equation}
	\label{eqn_ring_structure_prod}
\bV_{g,n}\otimes_\bC \bV_{g,n}\to \bV_{g,n}.
\end{equation}
Similarly, let $D$ be an oriented disk, then the triple
$$
(D\times \Sigma, D\times \{p_1,\dots,p_n\}, \emptyset)
$$
defines a cobordism from $(\emptyset, \emptyset, \emptyset)$ to $(S^1\times \Sigma, S^1\times \{p_1,\dots,p_n\},\emptyset)$, which defines a homomorphism
\begin{equation}
	\label{eqn_ring_structure_unit}
	\bC\cong \II (\emptyset,\emptyset,\emptyset)\to \bV_{g,n}.
\end{equation}
In general, cobordism homomorphisms of singular instanton Floer homology are only well-defined up to signs. The maps in \eqref{eqn_ring_structure_prod} and \eqref{eqn_ring_structure_unit} have canonical choices of signs induced by the almost complex structures on $F$ and $D$. The ring structure on $\bV_{g,n}$ is defined by taking the multiplicative structure to be \eqref{eqn_ring_structure_prod} and taking the unit element to be the image of the canonical generator under \eqref{eqn_ring_structure_unit}. The fact that this is a ring structure on $\bV_{g,n}$ follows from the functoriality of cobordism maps.

\subsection{The homomorphism $\Phi$ and the $\mu$ maps}
\label{subsec_mu_maps}
Let 
$$
\bA_{g,n} = \bC[\alpha,\beta, \psi_1,\dots, \psi_{2g},\delta_1,\dots,\delta_n]
$$
be a free graded-commutative $\bC$--algebra, where the degrees of the generators are given by 
$$
\deg \alpha = 2, \quad \deg\beta = 4, \quad \deg\psi_i =3, \quad \deg \delta_i = 2.
$$
Let $\epsilon$ be a variable with degree $2$.
A surjective homomorphism
$$
\Phi: \bA_{g,n}[\epsilon]/(\epsilon^2-1)\twoheadrightarrow \bV_{g,n}
$$
was defined by \cite{XZ:excision}*{(18)}, where the surjectivity was proved in \cite{XZ:excision}*{Proposition 2.11}. The ring $\bV_{g,n}$ is relatively graded over $\bZ/4$, and the map $\Phi$ preserves the relative $\bZ/4$ gradings. 
We briefly review the definition of $\Phi$ to set up the notation.

The homomorphism $\Phi$ is defined via the $\mu$ maps on $\bV_{g,n}$. 
The $\mu$ maps are homomorphisms on gauge-theoretic Floer homology groups defined by taking intersections of the moduli spaces of solutions to the gauge-theoretic equations with certain homology classes. The idea was first introduced in \cite{donaldson1990polynomial, D-Floer}.
For the definition of the $\mu$ maps in the context of singular instanton Floer homology, see \cite{Street}*{Section 2.3.2} and \cite{XZ:excision}*{Section 2.1}. The definitions involve some choices of normalizing constants, and we will follow the convention in \cite{XZ:excision}, which agrees with the conventions in \cite{donaldson1990polynomial, KM:suture}.

Let $\{a_i\}_{1\le i\le 2g}$ be a set of oriented simple closed curves in $\Sigma$ such that the intersection number of $a_i$ and $a_j$ is $1$ if $j=i+g$, and is zero if $|i-j|\neq g$. 
The following $\mu$--maps were defined in \cite{XZ:excision}*{Section 2.1}: 
$$\mu(\Sigma),\quad \mu(\pt),\quad \mu(a_i),\quad  \sigma_p \text{ for }p\in S^1\times \{p_1,\dots,p_n\}.$$
Let
$E:\bV_{g,n}\to \bV_{g,n}$ be the cobordism map induced by the trivial cobordism pair $([0,1]\times S^1\times \Sigma, [0,1]\times S^1\times \{p_1,\dots,p_n\})$ with a Stiefel--Whitney class represented by the Poincar\'e dual of $\{1/2\}\times \{\pt\}\times \Sigma$ (see  \cite{XZ:excision}*{Section 2.4}). Then we have $E^2=\id$. The choice of sign for the map $E$ is discussed in Section \ref{subsec_sign_E} below.

Let $e\in \bV_{g,n}$ be the unit element of the ring structure.
The homomorphism $\Phi$ is defined by taking 
$$
\Phi(\alpha) = \mu(\Sigma)(e), \,\Phi(\beta) = \mu(\pt)(e), \,\Phi(\psi_i) = \mu(a_i)(e), \,\Phi(\delta_i) = \sigma_{(x,p_i)}(e), \, \Phi(\epsilon) = E(e),
$$
where $x\in S^1$ is an arbitrary point.
The fact that $\Phi$ is well-defined follows from the graded commutativity of the $\mu$ maps and the functoriality of cobordism maps.

Since $E^2 = 1$, 
the space $\bV_{g,n}$ decomposes as a direct sum
 $$\bV_{g,n} = \bV_{g,n}^+\oplus \bV_{g,n}^-,$$
 where $\bV_{g,n}^\pm$ is the eigenspace of $E$ with eigenvalue $\pm 1$. 
 Each $\bV_{g,n}^\pm$ is a ring itself with the unit element given by $(e\pm E(e))/2$. The map $\Phi$ induces two surjective ring homomorphisms
 $$
 \Phi^{\pm}:\bA_{g,n}\twoheadrightarrow \bV_{g,n}^\pm
 $$
 by
 $$
 \Phi^\pm(f) = \frac{\id \pm E}{2} (\Phi(f)).
 $$
 
 \subsection{The sign convention of $E$}
 \label{subsec_sign_E}
 Since cobordism maps between singular instanton Floer homology are in general only defined up to signs, we need to choose a sign for the map $E$. Our sign convention is defined as follows.
 
 Fix an orientation on the link $S^1\times \{p_1,\dots,p_n\}\subset S^1\times \Sigma$ using the standard orientation of $S^1$, and orient $[0,1]\times S^1\times \{p_1,\dots,p_n\}$ so that it is an oriented cobordism. By \cite{KM:YAFT}, these orientations (together 
 with the integral lift $[\{1/2\}\times \{\pt\}\times \Sigma]$ of the Stiefel--Whitney class) fix a choice of sign for the cobordism map, and we denote the (signed) map as $\hat E$. The map $E$ is defined by 
 $$E=(-1)^{(n-1)/2}\hat E.$$ 
 The reason for introducing the extra factor of $(-1)^{(n-1)/2}$ is that the cobordism map $\pi^{g,n+2}_{g,n}$ introduced in \cite{Street}*{Section 2.6} and the ``U-cobordism'' in \cite{XZ:excision}*{Lemma 4.5} are not compatible with the above link orientation, and hence they are anti-commutative with $\hat E$. The extra sign in the definition of $E$ makes these cobordism maps commutative with $E$. As a result, Propositions \ref{prop_incl_pi_n+2_n} and \ref{prop_sequence_of_eigenvalues} below hold with respect to this sign convention of the map $E$. 
 
  We will also denote $(-1)^{(n-1)/2} \epsilon \in \bA_{g,n}$ by $\hat{\epsilon}$, so that $\Phi(\hat\epsilon )=\hat E(e)$.
 
  \begin{remark}
 	\label{rmk_notation}
 	Our notation of $\epsilon$, $\alpha$ and $\delta_i$ is different from the notation in \cite{kronheimer2022relations}.  The variables $\epsilon, \alpha, \delta_i$ in \cite{kronheimer2022relations} correspond to $(-1)^{(n-1)/2}\epsilon=\hat \epsilon $, 
 	$\alpha+(\delta_1+\dots+\delta_n)/2$, $-\delta_i$ in our notation. We will later use $\omega$ to denote $\alpha+(\delta_1+\dots+\delta_n)/2$.
 \end{remark}
 
 \subsection{Decomposition by the $\Sp(2g;\bZ)$ action}
 Define 
 $$
 \gamma = \sum_{i=1}^g \psi_i\psi_{i+g} \in \bA_{g,n}.
 $$
 Let $[a_i] \in H_1(\Sigma;\bC)$ be the fundamental class of $a_i$ in $\Sigma$ (recall that the curves $a_i$ are defined in Section \ref{subsec_mu_maps}). 
 Identify $\spann\{\psi_1,\dots,\psi_{2g}\}$ with $H_1(\Sigma;\bC)$ by identifying $\psi_i$ with $[a_i]$. 
 Let $e^1,\dots,e^{2g}$ be the dual basis of $[a_1],\dots,[a_{2g}]$ in $H^1(\Sigma, \bC)$. 
 For $k=0,\dots, g$, let $\Lambda_0^k H_1(\Sigma;\bC)$ be the kernel of the map
 $$
 \Lambda^k H_1(\Sigma;\bC) \to \Lambda^{k-2} H_1(\Sigma;\bC)
 $$
defined by taking contraction with $\sum_{i=1}^g e^i e^{i+g}$, where $\Lambda^i H_1(\Sigma;\bC) = 0$ if $i<0$. Then we have (see \cite{XZ:excision}*{Section 3.1})
$$
\bA_{g,n} = \bigoplus_{k=0}^g \Lambda_0^k H_1(\Sigma;\bC)\otimes \bC[\alpha, \beta, \gamma, \delta_1,\dots, \delta_n]/(\gamma^{g-k+1}).
$$
We will abuse notation and use $\Phi^\pm$ to denote the composition of the projection
\begin{align}
	\label{eqn_projection_decom_Sp(2g)}
& \bigoplus_{k=0}^g \Lambda_0^k H_1(\Sigma;\bC)\otimes \bC[\alpha, \beta, \gamma, \delta_1,\dots, \delta_n]
\nonumber
\\
&\quad \twoheadrightarrow  \bigoplus_{k=0}^g \Lambda_0^k H_1(\Sigma;\bC)\otimes \bC[\alpha, \beta, \gamma, \delta_1,\dots, \delta_n]/(\gamma^{g-k+1}) 
=\bA_{g,n}
\end{align}
with
$\Phi^\pm$.
We have:
\begin{Proposition}[\cite{XZ:excision}*{Proposition 3.4}]
	\label{prop_kernel_Phi_decompse_Sp(2g)}
	The kernel of
	$$
	\Phi^\pm:  \bigoplus_{k=0}^g \Lambda_0^k H_1(\Sigma;\bC)\otimes \bC[\alpha, \beta, \gamma, \delta_1,\dots, \delta_n]\twoheadrightarrow \bV_{g,n}^\pm
	$$
	has the form
	$$
	\ker \Phi^\pm = \bigoplus_{k=0}^g \Lambda_0^k  H_1(\Sigma;\bC) \otimes J_{g,n,k}^{\pm},
	$$
	where $J_{g,n,k}^\pm$ are ideals in $ \bC[\alpha, \beta, \gamma, \delta_1,\dots, \delta_n]$.
\end{Proposition}
By definition, 
$\gamma^{g-k+1}\in J_{g,n,k}^\pm.$
We also have the following result:
\begin{Proposition}[\cite{XZ:excision}*{Corollary 4.2}]
 $J_{g,n,k}^\pm = J_{g-k,n,0}^\pm$.
\end{Proposition}

\begin{Definition}
	Define $J_{g,n}^\pm = J_{g,n,0}^\pm$.
\end{Definition}

Then we have $J_{g,n,k}^\pm = J_{g-k,n}^\pm$, and 
\begin{equation}
	\label{eqn_gamma^g+1_in_Jgn}
	\gamma^{g+1}\in J_{g,n}^\pm. 
\end{equation}

 \begin{Proposition}
 	\label{prop_J_gn_pm_symmetry}
 	Let $\varphi$ be the automorphism of $\bC[\alpha,\beta,\gamma,\delta_1,\dots,\delta_n]$ that takes $(\alpha,\beta,\gamma,\delta_1,\dots,\delta_n)$ to $(-\alpha,\beta,-\gamma,-\delta_1,\dots,-\delta_n)$.
 	Then $\varphi(J_{g,n}^+) = J_{g,n}^-.$
 \end{Proposition}

\begin{proof}
Since $\Phi$ preserves the relative $\bZ/4$ grading, 
$\ker \Phi$ is a  homogeneous ideal with respect to the relative $\bZ/4$ grading.  Hence $\ker \Phi \cap \bC[\alpha,\beta,\gamma,\delta_1,\dots,\delta_n,\epsilon]/(\epsilon^2-1)$ is invariant under the automorphism of $\bC[\alpha,\beta,\gamma,\delta_1,\dots,\delta_n,\epsilon]/(\epsilon^2-1)$ sending $(\alpha,\beta,\gamma,\delta_1,\dots,\delta_n,\epsilon)$ to $(-\alpha,\beta,-\gamma,-\delta_1,\dots,-\delta_n,-\epsilon)$. The desired result then follows from the definition of $J_{g,n}^\pm$. 
\end{proof}

\subsection{Flip symmetries}
\label{subsec_flip_symm}
There is a \emph{flip symmetry} on $\bC[\alpha,\beta,\gamma,\delta_1,\dots,\delta_n]$  defined as follows. Let
\begin{equation}
	\label{def_omega}
\omega = \alpha+ \frac{\delta_1+\dots+\delta_n}{2}\in  \bC[\alpha,\beta,\gamma,\delta_1,\dots,\delta_n],
\end{equation}
then $\omega$ is homogeneous with degree $2$, and $\bC[\alpha,\beta,\gamma,\delta_1,\dots,\delta_n]$ is a free $\bC$--algebra generated by $\omega,\beta,\gamma,\delta_i$. 

We will identify 
$\bC[\alpha,\beta,\gamma,\delta_1,\dots,\delta_n]$ with $\bC[\omega,\beta,\gamma,\delta_1,\dots,\delta_n]$ via \eqref{def_omega}.
Define 
$$
\tau_i:  \bC[\alpha,\beta,\gamma,\delta_1,\dots,\delta_n]\to  \bC[\alpha,\beta,\gamma,\delta_1,\dots,\delta_n]
$$
to be the automorphism that fixes $\omega, \beta, \gamma$ and all $\delta_j$'s with $j\neq i$, and takes $\delta_i$ to $-\delta_i$. Then $\{\tau_i\}_{1\le i\le n}$ generates an abelian group of automorphisms on  $\bC[\alpha,\beta,\gamma,\delta_1,\dots,\delta_n]$ that is isomorphic to $(\bZ/2)^n$. An element in this group is called a \emph{flip symmetry}.  We say that a flip symmetry is \emph{even} if it preserves $\delta_1\delta_2\dots\delta_n$, and we say that it is odd if it takes $\delta_1\delta_2\dots\delta_n$ to $-\delta_1\delta_2\dots\delta_n$.

For $I\in \{1,\dots,n\}$, we define
$\tau_I = \prod_{i\in I}\tau_i.$
Since the $\tau_i$'s are commutative to each other, the definition of $\tau_I$ does not depend on the ordering of the product. By definition, $\tau_I$ is even if and only if $I$ contains an even number of elements.

\begin{Proposition}[\cite{Street}*{Section 2.2.2}, see also \cite{XZ:excision}*{Section 2.3}]
	\label{prop_Jgn_flip}
The ideals $J_{g,n}^\pm$ are invariant under all even flip symmetries.
\end{Proposition}

We will also frequently consider rings of the form 
$\bC[\alpha,\beta,\gamma,\delta_1,\dots,\delta_n]/\mathfrak{a},$
where $\mathfrak{a}$ is an ideal invariant under flip symmetries.  In this case, the flip symmetries on $\bC[\alpha,\beta,\gamma,\delta_1,\dots,\delta_n]$ induce automorphisms on the quotient ring.  We will abuse notation and call the induced automorphisms as \emph{flip symmetries}, and we will use $\tau_I$ to denote the corresponding maps on the quotient ring. If 
\begin{equation}
	\label{eqn_prod_delta_i_notin_I}
	\delta_1\dots\delta_n\neq 0 \in \bC[\alpha,\beta,\gamma,\delta_1,\dots,\delta_n]/\mathfrak{a},
\end{equation}
we say that a flip symmetry is \emph{even} if it preserves $\delta_1\delta_2\dots\delta_n$, and we say that it is odd if it takes $\delta_1\delta_2\dots\delta_n$ to $-\delta_1\delta_2\dots\delta_n$.

\subsection{The moduli space of parabolic bundles and the map $\Psi$}
There is another surjective homomorphism
$$
\Psi:\bA_{g,n} \twoheadrightarrow H^*(N_{g,n}^d;\bC)
$$
introduced by \cite{XZ:excision}*{Section 2.5}, which we now review.  The space $N_{g,n}^d$ is the moduli space of rank-2 stable parabolic bundles on $(\Sigma,\{p_1,\dots,p_n\})$ with degree $d$, where the weights are specified by \cite{XZ:excision}*{Definition 2.5}. When $d$ is even, the moduli space of the critical points of of the Chern--Simons functional defining $\bV_{g,n}$ is diffeomorphic to two copies of $N_{g,n}^0$ by \cite{MSesh} (see also \cite{XZ:excision}*{Theorem 2.9}). 

Similar to the maps $\Phi^\pm$, we will abuse notation and use $\Psi$ to denote the composition of $\Psi$ with the projection \eqref{eqn_projection_decom_Sp(2g)}. We have the following proposition.
\begin{Proposition}[\cite{XZ:excision}*{Proposition 3.1}]
	\label{prop_decompose_ker_Psi_Sp(2g)}
The kernel of
\begin{equation}
	\label{eqn_Psi_decompose_Sp(2g)}
\Psi:  \bigoplus_{k=0}^g \Lambda_0^k H_1(\Sigma;\bC)\otimes \bC[\alpha, \beta, \gamma, \delta_1,\dots, \delta_n]\twoheadrightarrow H^*(N_{g,n}^d;\bC)
\end{equation}
has the form
$$
\ker \Psi = \bigoplus_{k=0}^g \Lambda_0^k  H_1(\Sigma;\bC) \otimes I_{g,n,k}^{d},
$$
where $I_{g,n,k}^{d}$ are ideals in $ \bC[\alpha, \beta, \gamma, \delta_1,\dots, \delta_n]$.
\end{Proposition}
The definition of $I_{g,n,k}^d$ implies that $\gamma^{g-k+1}\in I_{g,n,k}$.

\begin{Proposition}[\cite{XZ:excision}*{Theorem 3.3}]
	\label{prop_singular_coh_ideal_g-k}
	The ideal $I_{g,n,k}^d$ only depends on $g-k$, $n$, and $d$. 
\end{Proposition}

\begin{Definition}
	Define $I_{g,n}^d = I_{g,n,0}^d$.
\end{Definition}
Then we have $I_{g,n,k}^d= I_{g-k,n}^d$. Moreover,
\begin{equation}
	\label{eqn_gamma^(g+1)_in_Ign}
	\gamma^{g+1}\in I_{g,n}^d.
\end{equation}
The map $\Psi$ is homogeneous with degree zero, so the ideals $I_{g,n}^d$ are homogeneous. 

The following result is also straightforward.
\begin{Proposition}
	\label{prop_singular_coh_flip_symmetry}
The ideal $I_{g,n}^d$ is invariant under even flip symmetries. 
If $d$ and $d'$ have the same parity, then 
$I_{g,n}^d = I_{g,n}^{d'}$.
If $d$ and $d'$ have different parities, then $I_{g,n}^d$ and  $I_{g,n}^{d'}$
differ by an odd flip symmetry.
\end{Proposition}
\begin{proof}
	There is a canonical diffeomorphism from $N_{g,n}^d$ to $N_{g,n}^{d+2}$ by taking tensor product with a fixed holomorphic line bundle on $\Sigma$ with degree $1$. It is straightforward to verify that the induced isomorphism on cohomology commutes with $\Psi$, therefore $I_{g,n}^d$ only depends on the parity of $d$. The effect of flip symmetry follows from \cite{Street}*{Section 1.4.3} and the correspondence of $N_{g,n}^d$ with the moduli space of flat connections (\cite{MSesh}, see also \cite{XZ:excision}*{Theorem 2.9}). 
\end{proof}

\begin{Definition}
	\label{def_leading_term}
Suppose $\bA$ is a graded ring over $\bC$ and $f\in \bA$ is a non-zero element. 
We say that $f_0$ is the \emph{leading-order term} of $f$ if $f_0$ is the homogeneous element in $\bA$ such that $\deg(f-f_0)<\deg(f)$.
\end{Definition}

The ideals $J_{g,n}^\pm$ can be thought of as graded deformations of the ideal $I_{g,n}^0$, in the following sense. 
 
\begin{Proposition}[\cite{Street}, see also \cite{XZ:excision}*{Corollary 2.13}]
	\label{prop_I_realized_by_J}
	Suppose $f\in I_{g,n}^0$ is a non-zero homogeneous element. Then there exist $\tilde f^+ \in J_{g,n}^+$ and $\tilde f^- \in J_{g,n}^-$ such that the leading-order terms of $\tilde f^+$ and $\tilde f^-$ are both $f$. 
\end{Proposition}

\begin{Proposition}[\cite{Street}, see also \cite{kronheimer2022relations}*{Proposition 3.14}]
		\label{prop_J_realized_by_I}
	Assume $f\in J_{g,n}^+$ or $J_{g,n}^-$. Then the leading-order term of $f$ is contained in $I_{g,n}^0$.
\end{Proposition}

Now we discuss the Poincar\'e polynomial of $H^*(N_{g,n}^d)$.
\begin{Definition}
	Assume $M$ is a complex vector space graded by non-negative integers, and let $M^{(i)}$ be the component of $M$ with degree $i$.  If $\dim_\bC M^{(i)}$ is finite for all $i$, we define the \emph{Poincar\'e series} of $M$ to be the series
	$$
	P_t(M) = \sum_{i=0}^\infty \dim_\bC M^{(i)} t^i.
	$$
	When the Poincar\'e series is of finite degree, it is also called the \emph{Poincar\'e polynomial}.
\end{Definition}

The Poincar\'e polynomial of $H^*(N_{g,n}^d)$ is computed by Street \cite{Street}.
\begin{Proposition}[\cite{Street}]
	\label{prop_Poincare_total}
	\begin{equation}
		\label{eqn_Poincare_total}
		P_t(H^*(N_{g,n}^d))=\frac{(1+t^2)^{n-1}(1+t^3)^{2g}- (2t)^{n-1} t^{2g}(1+t)^{2g}}{(1-t^2)^2}.
	\end{equation}
\end{Proposition}

\begin{proof}
	By Proposition \ref{prop_singular_coh_flip_symmetry}, $	P_t(H^*(N_{g,n}^d)) = 	P_t(H^*(N_{g,n}^0))$.
	By \cite{Street}*{Equation (1.3.12)},
	$$
	P_t(N_{g,1}^0)=\frac{(1+t^3)^{2g}-t^{2g}(1+t)^{2g}}{(1-t^2)^2}.
	$$
	Moreover, by \cite{Street}*{Corollary 1.3.6},
	$$
	P_t(N_{g,n+2}^0)=(1+t^2)^2P_t(N_{g,n}^0) +2^{n-1}t^{2g+n-1}(1+t)^{2g}.
	$$
	The desired result then follows from a straightforward induction argument.
\end{proof}

\subsection{Properties of $I_{g,n}$ and $J_{g,n}$}
\label{subsec_properties_Ign_Jgn}
We collect several algebraic properties of the ideals $I_{g,n}$ and $J_{g,n}$.
The next result is an immediately consequence of the definition of $\Psi$. 
\begin{Proposition}
	\label{prop_Ign_symmetric_delta}
	The ideal $I_{g,n}^d$ is symmetric with respect to permutations of $\delta_1,\dots,\delta_n$. 
\end{Proposition}

\begin{Proposition}[\cite{Street}]
	\label{prop_multiply_by_gamma_g+1}
	$I_{g+1,n}^d \cap (\gamma) = \gamma \cdot I_{g,n}^d$. 
\end{Proposition}
\begin{proof}
By Proposition \ref{prop_singular_coh_flip_symmetry}, we only need to consider the case when $d$ is even. In this case, $N_{g,n}^d$ is identified with the space $\mathcal{R}_{g,n}$ in \cite{Street}. By \cite{Street}*{Proposition 1.5.8} (see also \cite{KingNewstead}*{Theorem 3.2} for the non-singular case), for each $f\in \bC[\alpha,\beta,\gamma,\delta_1,\dots,\delta_n]$ and for each $j\le g+1$, we have 
$$
\int_{N_{g+1,n}^d}\Psi(\psi_j\psi_{j+g+1}f) = \int_{N_{g,n}^d} \Psi(f),
$$
where we use $\Psi$ to denote both the map to $H^*(N_{g+1,n}^d)$ and the map to $H^*(N_{g,n}^d)$. Since $\Psi(g)=0\in H^*(N_{g,n}^d)$ if and only if $\int_{N_{g,n}^d}\Psi(gh)=0$ for all $h$, the result follows from the fact that $\gamma = \sum_{j=1}^{g+1} \psi_j\psi_{j+g+1} \in \bA_{g+1,n}$.
\end{proof}

\begin{Proposition}[\cite{Street}*{Lemma 1.5.10}]
	\label{prop_sing_coh_delta^2_beta}
	For each $i$, we have $\delta_i^2+\beta \in I_{g,n}^d$.
\end{Proposition}

Now we review a computation of Mumford relations from \cite{XZ:excision}*{Section 3}. 
Let $m=(n-1)/2$.
Define
$$
F(t) = (1+\beta t^2)^{-1/2+m/2}\Big(\frac{1-t\sqrt{-\beta}}{1+t\sqrt{-\beta}}\Big)^{\frac{2\alpha\beta + \gamma}{4(-\beta)^{3/2}}} \exp\Big(-\frac{t\gamma}{2\beta}\Big),
$$
and write $F(t)$ as a formal power series in $t$:
\begin{equation}
	\label{eqn_defn_xi_kn}
	F(t) = \sum_{k=0}^\infty \xi_{k,n}(\alpha, \beta, \gamma) t^k.
\end{equation}
Since $F(t)$ satisfies the following equation
$$
(1+\beta t^2)\frac{F'(t)}{F(t)} = -\frac{\gamma}{2}t^2 + (m-1)\beta t + \alpha,
$$
we have a recursive formula for $\xi_{k,n}$ when $k\ge 2$:
\begin{equation}
		\label{eqn_recursive_xi}
(k+1)\xi_{k+1,n}  = \alpha \xi_{k,n} + (m-k)\beta \xi_{k-1,n}  - \frac{\gamma}{2} \xi_{k-2,n}.
\end{equation}
It is also straightforward to verify that $\xi_{0,n}=1$, $\xi_{1,n} = \alpha$, $\xi_{2,n} = \alpha^2/2 + (m-1)\beta/2$,
so $\xi_{k,n}$ is a homogeneous polynomial  in $\alpha, \beta, \gamma$ with degree $2k$.  We have:
\begin{Proposition}[\cite{XZ:excision}*{Proof of Proposition 3.5}]
	\label{prop_xi_in_ideal}
	If $d+m$ is odd, then 
$$
\xi_{k,n} \in I_{g,n}^d
$$
for all $k\ge g+m$. 
\end{Proposition}

We will also need the following results on $J_{g,n}^\pm$. 

\begin{Proposition}[\cite{Kr-ob}*{Proposition 4.1}, see also \cite{XZ:excision}*{(2)}]
For each pair $(g,n)$ and each $i$, we have $\delta_i^2+\beta-2 \in J_{g,n}^+$,  $\delta_i^2+\beta-2 \in J_{g,n}^-$. 
\end{Proposition}

\begin{Proposition}[\cite{XZ:excision}*{Lemmas 4.3, 4.4}]
	$\gamma J_{g,n}^+\subset J_{g+1,n}^+\subset J_{g,n}^+$, and $\gamma J_{g,n}^-\subset J_{g+1,n}^-\subset J_{g,n}^-$.
\end{Proposition}


Let 
\begin{equation}
		\label{eqn_defn_pi}
\pi_{g,n}^{g,n+2}: \bC[\alpha,\beta,\gamma,\delta_1,\dots,\delta_n,\delta_{n+1},\delta_{n+2}]\to  \bC[\alpha,\beta,\gamma,\delta_1,\dots,\delta_n]
\end{equation}
be the homomorphism that takes $\alpha$ to $\alpha$, $\beta$ to $\beta$, $\gamma$ to $\gamma$, and $\delta_i$ to $\delta_i$ ($i=1,\dots,n$), and takes $\delta_{n+1}$ to $-\delta_n$, takes $\delta_{n+2}$ to $\delta_n$. We have
\begin{Proposition}[\cite{Street}*{Corollary 2.6.8}]
	\label{prop_incl_pi_n+2_n}
	$$
	\pi^{g,n+2}_{g,n}(J_{g,n+2}^+)\subset J_{g,n}^+,
	$$
	$$
	\pi^{g,n+2}_{g,n}(J_{g,n+2}^-)\subset J_{g,n}^-.
	$$
\end{Proposition}

\subsection{Eigenvalues}
Since the map $\Phi^+$ induces an isomorphism
\begin{equation}
	\label{eqn_iso_induced_by_Phi+}
\bV_{g,n}^+ \cong \bC[\alpha,\beta,\gamma,\delta_1,\dots,\delta_n]/J_{g,n}^+,
\end{equation}
and $\bV_{g,n}^+$ is a finite-dimensional linear space, we may view the multiplications by $\alpha,\beta,\gamma,\delta_1,\dots,\delta_n$ as linear maps on $\bV_{g,n}^+$ and consider their eigenvalues and (generalized) eigenspaces. When there is no risk of confusion, we will abuse notation and use $\alpha,\beta,\gamma,\delta_1,\dots,\delta_n$ to denote the linear operators on $\bV_{g,n}^+$ defined by the multiplications by these elements via the isomorphism \eqref{eqn_iso_induced_by_Phi+}.

We recall the following result from \cite{XZ:excision}.
\begin{Proposition}
	\label{prop_sequence_of_eigenvalues}
	\cite{XZ:excision}*{Corollary 5.4 and Lemma 5.6}
	There exists a sequence of integers $\lambda_1,\lambda_2,\dots,$ such that the following holds
	\begin{enumerate}
		\item $|\lambda_i| = 2i-1$ for all $i$.
		\item For each pair $(g,n)$, let $m=(n-1)/2$, then the set of simultaneous eigenvalues for
		$$
		(\alpha, \gamma, \delta_1,\dots,\delta_n)
		$$
		on $\bV_{g,n}^+$ in the generalized eigenspace of $\beta$ with eigenvalue $2$ is given by $(\lambda,0,\dots,0)$ for $\lambda = \lambda_1,\dots,\lambda_{g+m}$. 
	\end{enumerate}
\end{Proposition}

A similar result holds for $\bV_{g,n}^-$ with the signs of the eigenvalues reversed.
\begin{Proposition}
	\label{prop_sequence_of_eigenvalues_V-}
	Let $\lambda_1,\lambda_2,\dots$ be as in Proposition \ref{prop_sequence_of_eigenvalues}. For each pair $(g,n)$, let $m=(n-1)/2$, then the set of simultaneous eigenvalues of
	$$
	(\alpha, \gamma, \delta_1,\dots,\delta_n)
	$$
	on $\bV_{g,n}^-$ in the generalized eigenspace of $\beta$ with eigenvalue $2$ is given by $(-\lambda,0,\dots,0)$ for $\lambda = \lambda_1,\dots,\lambda_{g+m}$. 
\end{Proposition}

Proposition \ref{prop_sequence_of_eigenvalues_V-} is an immediate consequence of Propositions \ref{prop_sequence_of_eigenvalues} and \ref{prop_J_gn_pm_symmetry}.

Eigenvalues directly reflect the properties of $J_{g,n}^\pm$ because of the following observation. We only state the result for $J_{g,n}^+$; the analogous result for $J_{g,n}^-$ is verbatim.
\begin{Lemma}
	\label{lem_eigenvalue_and_ideal}
	Suppose $(\lambda_\alpha,\lambda_\beta,\lambda_\gamma,\lambda_{\delta_1},\dots,\lambda_{\delta_n})$ is a tuple of simultaneous eigenvalues for the operators $(\alpha,\beta,\gamma,\delta_1,\dots,\delta_n)$ on $\bV_{g,n}^+$. Let $\varphi$ be the homomorphism from $\bC[\alpha,\beta,\gamma,\delta_1,\dots,\delta_n]$ to $\bC$ taking $(\alpha,\beta,\gamma,\delta_1,\dots,\delta_n)$ to $(\lambda_\alpha,\lambda_\beta,\lambda_\gamma,\lambda_{\delta_1},\dots,\lambda_{\delta_n})$. Then $\varphi(J_{g,n}^+)= (0)$.
\end{Lemma}
\begin{proof}
Let $f(\alpha,\beta,\gamma,\delta_1,\dots,\delta_n)\in J_{g,n}^+$, let $v\in \bV_{g,n}^+$ be a simultaneous eigenvector for $(\alpha,\beta,\gamma,\delta_1,\dots,\delta_n)$ with eigenvalues $(\lambda_\alpha,\lambda_\beta,\lambda_\gamma,\lambda_{\delta_1},\dots,\lambda_{\delta_n})$. Then we have $f(\alpha,\beta,\gamma,\delta_1,\dots,\delta_n)(v) = 0\in \bV_{g,n}^+$, so $f(\lambda_\alpha,\lambda_\beta,\lambda_\gamma,\lambda_{\delta_1},\dots,\lambda_{\delta_n})=0\in\bC$.
\end{proof}

%
%

\section{Generators of the ideal $I_{g,n}^d$}
\label{sec_parabolic}
	Throughout this section, we will use $n$ to denote a positive odd integer and define $m=(n-1)/2$. Let $g\ge 0$ be the genus of $\Sigma$. 
	
The main result of this section is the following theorem.
\begin{Theorem}
	\label{thm_generators_of_Ign}
	Assume $d+m$ is odd. Then the ideal $I_{g,n}^d$ is generated by $\delta_i^2+\beta$ $(i=1,\dots,n)$,  $\gamma^{g+1}$, and the even flip symmetries of $\xi_{g+m,n}$, $\xi_{g+m+1,n}$, $\xi_{g+m+2,n}$.
\end{Theorem}

\begin{remark}
	\label{rem_earl_kirwan}
	The proof of Theorem \ref{thm_generators_of_Ign} does not depend on the work of Earl and Kirwan \cite{earl2004complete}.
	Since $\xi_{k,n}$ are Mumford relations, it is natural to ask whether Theorem \ref{thm_generators_of_Ign} can also be proved using Earl and Kirwan's results. It turns out that there are three difficulties:
	\begin{enumerate}
		\item The Mumford relations used in \cite{earl2004complete} contain not only polynomials of the form $\xi_{k,n}$, but also polynomials with other forms (which are obtained by taking slant products of the universal Mumford relations with homology classes on the Jacobian with degrees \emph{less than} $2g$).
		\item 	Earl--Kirwan \cite{earl2004complete} only discussed parabolic bundles with one marked point, while we need to consider parabolic bundles with $n$ marked points.
		\item The polynomials $\xi_{k,n}$ are obtained by considering parabolic maps with slope difference $1/4$ (in other words, $\hat{d}/\hat{n} - d/n =1/4$ in the notation of \cite{earl2004complete}*{Theorem 2.1}). Earl--Kirwan's theorem also requires considering parabolic maps with slope difference  $3/4$.  
	\end{enumerate}
    The first difficulty turns out to be resolvable by reducing the structure of the generating set using Proposition \ref{prop_singular_coh_ideal_g-k}. 
	The second difficulty is hypothetically resolvable by extending the arguments of \cite{earl2004complete} to multiple marked points. 
	However, because of the third difficulty, the best result one can get from the hypothetical extended version of Earl--Kirwan's theorem is the following: $I_{g,n}^d$ is generated by all polynomials listed in Theorem \ref{thm_generators_of_Ign} together with all the odd flip symmetries of $\xi_{g+m+1,n+2}$, $\xi_{g+m+2,n+2}$, $\xi_{g+m+3,n+2}$. It is not clear to the authors how to remove the additional generators without going through a substantial part of the current proof of Theorem \ref{thm_generators_of_Ign}.
\end{remark}

\begin{remark}
	\label{rem_Ign_g=0}
	When $g=0$, the ideal $I_{0,n}^d$ has been studied thoroughly by \cite{Street} and \cite{kronheimer2022relations}. The special case of Theorem \ref{thm_generators_of_Ign} for $g=0$ was proved in \cite{kronheimer2022relations} using results from \cite{Street}. More precisely,  \cite{kronheimer2022relations}*{Proposition 4.8} proved that if $d+m$ is odd, then $I_{0,n}^d$ is generated by $\delta_i^2+\beta$ ($i=1,\dots,n$), $\gamma$, and the even flip symmetries of $\xi_{m,n}$. 
	By the recursive relation \eqref{eqn_recursive_xi}, we know that $\xi_{m+1,n}, \xi_{m+2,n} \in (\xi_{m,n},\gamma)$, therefore this statement is equivalent to the $g=0$ case of Theorem \ref{thm_generators_of_Ign}.
\end{remark}

The general idea of the proof of Theorem \ref{thm_generators_of_Ign} is to show that the ideal generated by the listed polynomials has the same Poincar\'e series as $I_{g,n}^d$. 
Since all the polynomials listed in Theorem \ref{thm_generators_of_Ign} are elements of $I_{g,n}^d$ by Section \ref{subsec_properties_Ign_Jgn} and Equation \eqref{eqn_gamma^(g+1)_in_Ign}, this will imply the the desired theorem.

The proof of Theorem \ref{thm_generators_of_Ign} is organized as follows. 
In Section \ref{subsec_Poincare_series}, we compute the Poincar\'e series of $I_{g,n}^d$. In Section \ref{subsec_modulo_gamma}, we show that we can simplify the  problem by taking quotient with the ideal generated by $\gamma$ and $\delta_i^2+\beta$ ($i=1,\dots,n$). Section \ref{subsec_poincare_series_Kgn} then computes the Poincar\'e series of the quotient image of $I_{g,n}^d$. Sections \ref{subsec_isotypic_decomposition} to \ref{subsec_properties_of_rho} investigate the properties of the quotient images of $\xi_{k,n}$ under the isotypic decomposition by the flip symmetry actions. Sections \ref{subsec_ideals_modulo_beta} and \ref{subsec_completion_of_proof_Ign_generators} will finish the proof of Theorem \ref{thm_generators_of_Ign} by combining an induction argument with an explicit construction of basis for part of the ideal.

We also record the following corollary of Theorem \ref{thm_generators_of_Ign} for later reference.
\begin{Corollary}
	\label{cor_generating_set_Ign}
		Assume $d+m$ is odd and $g+m\ge 1$. Then the ideal $I_{g,n}^d$ is generated by $\delta_i^2+\beta$ $(i=1,\dots,n)$,  $\gamma^{g+1}$, and the even flip symmetries of $\xi_{g+m,n}$, $\xi_{g+m+1,n}$, $\gamma\cdot \xi_{g+m-1,n}$.
\end{Corollary}
\begin{proof}
This is an immediate consequence of Theorem \ref{thm_generators_of_Ign} and the recursive formula \eqref{eqn_recursive_xi}.
\end{proof}

\subsection{The Poincar\'e series of $I_{g,n}^d$}
\label{subsec_Poincare_series}
Recall that $N_{g,n}^d$ in \eqref{eqn_Psi_decompose_Sp(2g)} denotes the moduli space of parabolic bundles over $(\Sigma,(p_1,\dots,p_n))$ with degree $d$ and weights specified by \cite{XZ:excision}*{Definition 2.5}. By Proposition \ref{prop_decompose_ker_Psi_Sp(2g)},  $H^*(N_{g,n}^d)$ is decomposed as an $\Sp(2g)$ representation as
\begin{equation}
	\label{eqn_decompose_H^*(R_g,n)}
	H^*(N_{g,n}^d) \cong \bigoplus_{k=0}^g \Lambda_0^k H_1(\Sigma, \bC) \otimes H^*(N_{g,n}^d)_k,
\end{equation}
where 
\begin{equation}
	\label{eqn_defn_H*Ngn}
H^*(N_{g,n}^d)_k = \bC[\alpha,\beta,\gamma,\delta_1,\dots,\delta_n]/I_{g,n,k}^d.
\end{equation}
The map $\Psi$ is defined by slant products with a degree $4$ cohomology class on the universal bundle, so \eqref{eqn_decompose_H^*(R_g,n)} preserves the degree if we assign elements in $H_1(\Sigma,\mathbb{C})$ with degree $4-1=3$ on the right-hand side.

By Propositions	\ref{prop_singular_coh_ideal_g-k} and Proposition \ref{prop_singular_coh_flip_symmetry}, the Poincar\'e polynomial of $H^*(N_{g,n}^d)_k$ only depends on $g-k$ and $n$.

To simplify notation, we use 
$P_t(g,n)$ to denote the Poincar\'e polynomial of $H^*(N_{g,n}^d)_0$. Then by the above discussion, the Poincar\'e polynomial of $H^*(N_{g,n}^d)_k$ is equal to $P_t(g-k,n)$. 
By the definitions, we have
\begin{equation}
	\label{eqn_Pt(Ignd)}
P_t(I_{g,n}^d) + P_t(g,n) = P_t(\bC[\alpha,\beta,\gamma,\delta_1,\dots,\delta_n]) = \frac{1}{(1-t^2)^{n+1}(1-t^4)(1-t^6)},
\end{equation}
so the Poincar\'e series of $I_{g,n}^d$ is determined by $P_t(g,n)$. 
\begin{Proposition}
	\label{prop_Ptgn}
\begin{equation}
	\label{eqn_Ptgn}
P_t(g,n) = \frac{(1+t^2)^{n-1}(1-t^{6g+6}) - (2t)^{n-1}t^{2g}(1-t^{2g+2})(1+t^2+t^4)}{(1-t^2)^2(1-t^6)}.
\end{equation}
\end{Proposition}

\begin{proof}
Note that for $k\le g$, we have
$$\dim_\bC \Lambda_0^k H_1(\Sigma,\bC) = {2g \choose k}-{2g\choose k-2},$$
where by definition, ${2g \choose a}=0$ if $a<0$.

Since \eqref{eqn_decompose_H^*(R_g,n)} preserves the degree if we assign elements in $H_1(\Sigma,\mathbb{C})$ with degree $3$, we have
\begin{equation}
	\label{eqn_total_decompose_Ptgn}
	P_t(H^*(N_{g,n}^d)) =\sum_{k=0}^{g}\Big[{2g \choose k}-{2g\choose k-2} \Big]  P_t(g-k,n)\cdot t^{3k}.
\end{equation}

Let $Q_t(g,n)$ be the right-hand side of \eqref{eqn_Ptgn}. We  first show that \eqref{eqn_total_decompose_Ptgn} holds if we replace $P_t(g,n)$ with $Q_t(g,n)$. In other words, we show that 
\begin{equation}
	\label{eqn_total_decompose_Qtgn}
	P_t(H^*(N_{g,n}^d)) =\sum_{k=0}^{g}\Big[{2g \choose k}-{2g\choose k-2} \Big]  Q_t(g-k,n)\cdot t^{3k}.
\end{equation}
By the definition of $Q_t(g,n)$, it is straightforward to verify that for $g\ge 2$,
\begin{equation}
	\label{eqn_Qt_g_g-2}
Q_t(g,n)-t^6 Q_t(g-2,n) = \frac{(1+t^2)^{n-1}(1+t^{6g})-(2t)^{n-1}(t^{2g}+t^{4g})}{(1-t^2)^2}.
\end{equation}
To simplify notation, we define $Q_t(g,n)=0$ if $g<0$. It is then straightforward to verify that \eqref{eqn_Qt_g_g-2} also holds for $g=1$. When $g=0$, we have
$$
\Big(Q_t(g,n) - t^6 Q_t(g-2,n)\Big)\Big|_{g=0} = \frac12 \frac{(1+t^2)^{n-1}(1+t^{6g})-(2t)^{n-1}(t^{2g}+t^{4g})}{(1-t^2)^2}\Big|_{g=0}.
$$ 

Therefore, the right-hand side of \eqref{eqn_total_decompose_Qtgn} equals
\begin{align}
	& \sum_{k=0}^g {2g \choose k} Q_t(g-k,n)\cdot t^{3k} - \sum_{k=0}^{g} {2g\choose k-2}Q_t(g-k,n)\cdot t^{3k} 
	\nonumber \\
= & \sum_{k=0}^g {2g \choose k} Q_t(g-k,n)\cdot t^{3k} - \sum_{k=0}^{g-2} {2g\choose k}Q_t(g-k-2,n)\cdot t^{3k+6}
\nonumber \\
= &  \sum_{k=0}^g {2g \choose k} Q_t(g-k,n)\cdot t^{3k} - \sum_{k=0}^{g} {2g\choose k}Q_t(g-k-2,n)\cdot t^{3k+6}
\nonumber \\
= & \sum_{k=0}^g {2g \choose k} \Big(Q_t(g-k,n)-t^6 Q_t(g-k-2,n)\Big) \cdot t^{3k}
\nonumber \\
= & \sum_{k=0}^{g-1} {2g \choose k} \Big(\frac{(1+t^2)^{n-1}(1+t^{6(g-k)})-(2t)^{n-1}(t^{2(g-k)}+t^{4(g-k)})}{(1-t^2)^2}\Big) \cdot t^{3k}
\nonumber \\
&\quad +\Bigg( {2g\choose k}\cdot \frac12\cdot \frac{(1+t^2)^{n-1}(1+t^{6(g-k)})-(2t)^{n-1}(t^{2(g-k)}+t^{4(g-k)})}{(1-t^2)^2}\cdot t^{3k}\Bigg)\Bigg|_{k=g}
\label{eqn_RHS_combination_of_Qt(g,n)}
\end{align}
Note that by the binomial formula, we have
$$(1+x)^{2g} = \sum_{k=0}^{g-1} {2g\choose k} (x^k + x^{2g-k}) + \Big(\frac12 {2g\choose k} (x^{2g} + x^{2g-k})\Big)\Big|_{k=g}.$$
 Therefore, \eqref{eqn_RHS_combination_of_Qt(g,n)} is equal to 
$$
\frac{(1+t^2)^{n-1}(1+t^3)^{2g}-(2t)^{n-1}t^{2g}(1+t)^{2g}}{(1-t^2)^2},
$$
which equals $P_t(H^*(N_{g,n}^d))$ by  Proposition \ref{prop_Poincare_total}. Therefore \eqref{eqn_total_decompose_Qtgn} is proved.

We now prove the desired equation \eqref{eqn_Ptgn} by induction on $g$. If $g=0$, we have $P_t(0,n)= P_t(H^*(R_{0,n}^d))$, and the desired result follows from Proposition \ref{prop_Poincare_total}. Now assume $P_t(g',n)=Q_t(g',n)$  for all $g'<g$, we show that $P_t(g,n)=Q_t(g,n)$. By \eqref{eqn_total_decompose_Ptgn} and \eqref{eqn_total_decompose_Qtgn}, we have 
\begin{align*}
P_t(g,n) 
&= P_t(H^*(N_{g,n}^d)) - \sum_{k=1}^{g}\Big[{2g \choose k}-{2g\choose k-2} \Big]  P_t(g-k,n)\cdot t^{3k} \\
&= P_t(H^*(N_{g,n}^d)) - \sum_{k=1}^{g}\Big[{2g \choose k}-{2g\choose k-2} \Big]  Q_t(g-k,n)\cdot t^{3k} \\
&= Q_t(g,n).
\end{align*}
Therefore the proposition is proved. 
\end{proof}

\subsection{Ideals modulo $\gamma$ and $\delta_i^2+\beta$}
\label{subsec_modulo_gamma}
\begin{Definition}
	Define
\begin{align}
	\bar{R}_n & = \bC[\alpha,\beta,\gamma,\delta_1,\dots,\delta_n]/(\gamma,\delta_1^2+\beta,\delta_2^2+\beta, \dots, \delta_n^2+\beta)
	\nonumber
	\\
	& = \bC[\alpha,\beta,\delta_1,\dots,\delta_n]/(\delta_1^2+\beta,\delta_2^2+\beta, \dots, \delta_n^2+\beta).
	\label{eqn_def_barR_n}
\end{align}
Define $K_{g,n}$ to be the quotient image of $I_{g,n}^d\subset \bC[\alpha,\beta,\gamma,\delta_1,\dots,\delta_n]$ in $\bar{R}_n$, where the parity of $d$ is taken so that $d+m$ is odd. 
\end{Definition}

It is straightforward to verify that the ideals in \eqref{eqn_def_barR_n} are invariant under flip symmetries, and the ring $\bar{R}_n$ satisfies \eqref{eqn_prod_delta_i_notin_I}. By the discussions in Section \ref{subsec_flip_symm}, it makes sense to refer to flip symmetries on $\bar{R}_n$, and each flip symmetry on $\bar{R}_n$ has a well-defined parity.

\begin{Definition}
	Let $\bar{\xi}_{k,n}$ be the image of $\xi_{k,n}$ in $\bar{R}_{n}$. 
	Let $K_{g,n}'$ denote the ideal in $\bar{R}_n$ generated by all the even flip symmetries of $\bar{\xi}_{g+m,n}$, $\bar{\xi}_{g+m+1,n}$. 
\end{Definition}

\begin{Lemma}
	\label{lem_K'gn_relation_in_g}
	$$
	\beta\,K_{g,n}'\subset K_{g+1,n}'\subset K_{g,n}'.
	$$
\end{Lemma}
\begin{proof}
	By \eqref{eqn_recursive_xi}, we have $\bar{\xi}_{g+m,n}\cdot\beta \in (\bar{\xi}_{g+m+1,n}, \bar{\xi}_{g+m+2,n})$. 
	Therefore the lemma follows from the definition of $K_{g,n}'$. 
\end{proof}

By Proposition \ref{prop_xi_in_ideal}, we know that 
\begin{equation}
	\label{eqn_Kgn_subset_Kgn'}
	K_{g,n}'\subset K_{g,n}.
\end{equation}
The next proposition states that the reversed inclusion is also true. The proof of this proposition will be postponed to the rest of Section \ref{sec_parabolic}.
\begin{Proposition}
	\label{prop_generators_of_Kgn}
	$K_{g,n} = K_{g,n}'$.
\end{Proposition}

We now prove Theorem \ref{thm_generators_of_Ign} assuming Proposition \ref{prop_generators_of_Kgn}.
\begin{proof}[Proof of Theorem \ref{thm_generators_of_Ign} assuming Proposition \ref{prop_generators_of_Kgn}]
	
	Let 
	$$I_{g,n}'\subset \bC[\alpha,\beta,\gamma,\delta_1,\dots,\delta_n]$$ be the ideal generated by $\delta_i^2+\beta$ $(i=1,\dots,n)$,  $\gamma^{g+1}$, and the even flip symmetries of $\xi_{g+m,n}$, $\xi_{g+m+1,n}$, $\xi_{g+m+2,n}$. Since we assume $d+m$ is odd, by the results in Section \ref{subsec_properties_Ign_Jgn} and Equation \eqref{eqn_gamma^(g+1)_in_Ign}, we have $I_{g,n}'\subset I_{g,n}^d.$
     By \eqref{eqn_recursive_xi} and the definition of $I_{g,n}'$, we have 
     	\begin{equation}
     	\label{eqn_inclusion_multiply_gamma}
     	 \gamma I_{g,n}'\subset I_{g+1,n}'.
     \end{equation}
   It also follows from the definition that the quotient image of $I_{g,n}'$ in $\bar{R}_n$ is $K_{g,n}'$.

	We prove Theorem \ref{thm_generators_of_Ign} by induction on $g$. The case when $g=0$ was proved by \cite{kronheimer2022relations}*{Proposition 4.8} (see Remark \ref{rem_Ign_g=0}).
	Suppose Theorem \ref{thm_generators_of_Ign} holds for $g$, we show that it also holds for $g+1$. Take $f\in I_{g+1,n}^d$. By Proposition \ref{prop_generators_of_Kgn}, there exists $f'\in I_{g+1,n}'$, such that $f-f'\in (\gamma)$. Since $f-f'\in I_{g+1,n}^d$, by Proposition \ref{prop_multiply_by_gamma_g+1}, there exists $h\in I_{g,n}^d$ such that $f-f' = \gamma \cdot h$. By the induction hypothesis, we have $h\in I_{g,n}'$, thus by \eqref{eqn_inclusion_multiply_gamma}, we have $\gamma\cdot h\in I_{g+1,n}'$. Therefore, $f\in I_{g+1,n}'$, and the theorem is proved.
\end{proof}

Sections \ref{subsec_poincare_series_Kgn} to \ref{subsec_completion_of_proof_Ign_generators} will be devoted to the proof of Proposition \ref{prop_generators_of_Kgn}.
Throughout the proof, we will always assume that $d$ is chosen so that $d+m$ is odd as required in the definition of the ideal $K_{g,n}'$.

\subsection{The Poincar\'e series of $K_{g,n}$}
\label{subsec_poincare_series_Kgn}
\begin{Lemma}
	\label{lem_Poincare_Kgn}
	The Poincar\'e series of $K_{g,n}$ is given by 
	\begin{equation}
	\label{eqn_Poincare_Kgn}
	P_t(K_{g,n}) = 2^{n-1} t^{2(g+m)}\frac{1+t^2-t^{2+2g}}{(1-t^2)^2}.
	\end{equation}
\end{Lemma} 
\begin{proof}
	Let $R_n = \bC[\alpha,\beta,\gamma,\delta_1,\dots,\delta_n]/(\delta_1^2+\beta,\dots,\delta_n^2+\beta).$
	Let $\hat{I}_{g,n}^d$ be the quotient image of $I_{g,n}^d$ in $R_n$. 
	By \eqref{eqn_defn_H*Ngn}, we have a short exact sequence
	$$
	0 \to \hat{I}_{g,n}^d\to R_n \to (H^*(R_{g,n}^d))_0 \to 0,
	$$
	where all maps above have degree zero. 
	Hence 
\begin{equation}
	\label{eqn_Pt_barI_gn0}
	P_t(\hat{I}_{g,n}^d) = P_t(R_n) - P_t((H^*(R_{g,n}^d))_0).
\end{equation}

	By Proposition  \ref{prop_multiply_by_gamma_g+1}, we have a short exact sequence
	$$
	 0 \to	\hat{I}_{g,n}^d \stackrel{\times\gamma}{\to} \hat{I}_{g+1,n}^d \to K_{g+1,n}\to 0,
	$$
	where the quotient map from $ \hat{I}_{g+1,n}^d$ to $K_{g+1,n}$ has degree zero.  Therefore,
\begin{equation}
	\label{eqn_Pt_Kgn_in_barI}
	P_t(K_{g+1,n}) = P_t(\hat{I}_{g+1,n}^d) - P_t(\hat{I}_{g,n}^d)\cdot t^6. 
\end{equation}

As a result, when $g\ge 1$, \eqref{eqn_Poincare_Kgn} follows from \eqref{eqn_Pt(Ignd)},  \eqref{eqn_Ptgn}, \eqref{eqn_Pt_barI_gn0}, \eqref{eqn_Pt_Kgn_in_barI},   and a straightforward computation.

When $g=0$, we have a short exact sequence
$$
0 \to K_{0,n} \to \bar{R}_n \to H^*(R_{0,n}^d)\to 0,
$$
where all the above maps are have degree zero. Hence 
$$
P_t(K_{0,n}) = P_t( \bar{R}_n) - P_t(H^*(R_{0,n}^d)),
$$
and the result follows from Proposition \ref{prop_Poincare_total}. 
\end{proof}

By definition, $K_{g,n}$ is supported only on even degrees. 
We will use 
$K_{g,n}^{\langle i \rangle}$ to denote the component of $K_{g,n}$ with degree $2(g+m+i)$.
By Lemma \ref{lem_Poincare_Kgn}, we have $K_{g,n}^{\langle i\rangle }=0$ whenever $i<0$, and 
\begin{equation}
	\label{eqn_dim_K_gn^0}
	\dim_\bC K_{g,n}^{\langle 0 \rangle } = 2^{n-1}.
\end{equation}
We also have the following result.
\begin{Lemma}
	\label{lem_graded_dim_Kgn_diff}
	When $g\ge 1$, we have
$$
\dim_\bC K_{g,n}^{\langle i \rangle} - \dim_\bC K_{g-1,n}^{\langle i-1 \rangle} = 
\begin{cases}
	0  & \text{ if } i < 0, \\
	2^{n-1}  & \text{ if } i = 0, \\
	2^n  & \text{ if } i > 0.
\end{cases}
$$
\end{Lemma}

\begin{proof}
By \eqref{eqn_Poincare_Kgn},
\begin{equation*}
P_t(K_{g,n}) - t^4 P_t(K_{g-1,n}) = 2^{n-1}\,t^{2(g+m)}\,\frac{1+t^2}{1-t^2} 
 = 2^{n-1} t^{2(g+m)} + 2^n\sum_{k=g+m+1}^\infty t^{2k}.
\end{equation*}
Therefore the result is proved.
\end{proof}

\subsection{Isotypic decomposition by even flip symmetries}
\label{subsec_isotypic_decomposition}
The set of even flip symmetries forms a finite abelian group that acts on $\bar{R}_n$. We consider its isotypic decomposition, which is given as follows. 

For $I\subset \{1,\dots,n\}$, we use $\delta^I$ to denote $\prod_{i\in I} \delta_i$, use $|I|$ to denote the number of elements in $I$, and use $I^c$ to denote the complement of $I$ in $\{1,\dots,n\}$. 

Recall that $\omega = \alpha+(\delta_1+\dots+\delta_n)/2$, and we use this equation to identify
$$
\bar{R}_n = \bC[\alpha,\beta,\delta_1,\dots,\delta_n]/(\delta_1^2+\beta,\dots,\delta_n^2+\beta) 
$$
with
$$
 \bC[\omega,\beta,\delta_1,\dots,\delta_n]/(\delta_1^2+\beta,\dots,\delta_n^2+\beta).
$$
For each $I$ with $|I|\le m$, define 
$$
\bar{R}_n^{I} = \bC[\omega,\beta]\cdot \delta^I \oplus \bC[\omega,\beta]\cdot \delta^{I^c}.
$$
Note that $\bar{R}_n^{I}$ is only a linear space; it is not an ideal in $\bar{R}_n$. 
Then the isotypic decomposition of $\bar{R}_n$ by the even flip symmetry actions is given by
\begin{equation}
	\label{eqn_bar_R_isotypic_decomp}
	\bar{R}_n = \bigoplus_{|I|\le m} \bar{R}_n^I.
\end{equation}
For $J\subset \{1,\dots,n\}$ such that $|J|$ is even, the action of $\tau_J$ on  $\bar{R}_n^I$ is the scalar multiplication by $(-1)^{|I\cap J|}$.

Since both $K_{g,n}$ and $K'_{g,n}$ are invariant under even flip symmetries, we have 
\begin{equation}
	\label{eqn_Kgn_Kgn'_isotypic_decompose}
\begin{split}
K_{g,n} & = \bigoplus_{|I|\le m} \bar{R}_n^I \cap K_{g,n}, \\
K_{g,n}' & = \bigoplus_{|I|\le m} \bar{R}_n^I \cap K_{g,n}'. 
\end{split}
\end{equation}

When $g=0$, Lemma \ref{lem_Poincare_Kgn} shows that 
$$
\dim_\bC K_{0,n}^{\langle i \rangle}= \begin{cases}
	0  & \text{ if } i < 0, \\
	2^{n-1}(i+1) & \text{ if } i \ge 0.
\end{cases}
$$
We have the following result about the isotypic decomposition of $K_{0,n}^{\langle i \rangle}$.

\begin{Lemma}
	\label{lem_g=0_isotypic_decomposition_Kgn}
If $i\ge 0$, then for all $I\subset\{1,\dots,n\}$ such that $|I|\le m$, we have
$$
\dim_\bC K_{0,n}^{\langle i \rangle} \cap \bar{R}_n^I = i+1.
$$
\end{Lemma}

\begin{proof}
   By \eqref{eqn_gamma^(g+1)_in_Ign}, we have $\gamma\in I_{0,n}$, so 
	$$
	\bar{R}_n/K_{0,n} \cong H^*(R_{0,n}^d)
	$$
	as graded rings. 
	The ring $H^*(R_{0,n}^d)$ has been studied thoroughly in \cite{Street,kronheimer2022relations}.
	In particular, let
	 $$V_m = \underset{a+b+|J|<m}{\spann}\{\omega^a\beta^b\delta^J\}$$
	 be a finite dimensional linear subspace of $\bar{R}_n$, then Lemma 1.6.3 and Proposition 1.6.4 in \cite{Street} showed that the composition map
	 $$
	 V_m \hookrightarrow \bar{R}_n \to \bar{R}_n/K_{0,n} \cong H^*(R_{0,n}^d)
	 $$
	 is a linear space isomorphism. The above maps are equivariant with respect to even flip symmetries. Therefore,
	$ \dim_\bC K_{0,n}^{\langle i \rangle} \cap \bar{R}_n^I $ is equal to the dimension of the degree $2(m+i)$ component of the linear space $\bar{R}_n^I/( \bar{R}_n^I\cap V_m)$. This is equal to the number of non-negative integers $a,b$ such that either 
	$$
	a+b+|I| \ge m, \quad a+2b+|I| = m+i,
	$$ 
	or 
	$$
	a+b+|I^c| \ge m, \quad a+2b+|I^c| = m+i.
	$$
	It is straightforward to verify that the number of such pairs $(a,b)$ is equal to $i+1$. 
\end{proof}

\subsection{Decomposition of $\bar{\xi}_{k,n}$}
\label{subsec_decompose_bar_xi}
For $s=0,\dots,n$, define 
$$\delta^{\{s\}} =\sum_{I\subset\{1,\dots,n\}, |I|=s} \delta^I =  \sum_{1\le a_1<a_2<\dots<a_s\le n} \delta^{a_1}\delta^{a_2}\dots\delta^{a_s}.$$
Since $\bar{\xi}_{k,n}$ is a polynomial in $\alpha$ and $\beta$, it is symmetric with respect to $\delta_1,\dots,\delta_n$. So there are unique polynomials $\rho_{k,n,s} \in \bC[\omega,\beta]$ such that
\begin{equation}
	\label{eqn_decompose_xi_by_rho}
\bar{\xi}_{k,n}= 2^{-m-1}\,\sum_{s=0}^n   \rho_{k,n,s} \cdot \delta^{\{s\}}.
\end{equation}
The term $2^{-m-1}$ is a normalizing factor that will be convenient later.

Note that since $\bar{\xi}_{k,n}$ has degree $2k$, we have $\rho_{k,n,s}=0$ whenever $s>k$. 
The following statements give a formula for $\rho_{k,n,s}$ when $s\le k$. 

\begin{Definition}
	\label{defn_rho_kn}
Suppose $r$ is an (either positive or negative) odd integer.
Define $\rho_{k,r}\in\bC[\omega, \beta]$ by the formal series expansion:
\begin{multline*}
\sum_{k=0}^\infty \rho_{k,r}\cdot 
t^k 
=
(1+\beta t^2)^{-3/4}\Big(\frac{1-t\sqrt{-\beta}}{1+t\sqrt{-\beta}}\Big)^{\frac{\omega}{2\sqrt{-\beta}}}
\\
\times \Big(\sqrt{1-t\sqrt{-\beta}} + \sqrt{1+t\sqrt{-\beta}}\Big)\cdot 
\big(1+\sqrt{1+\beta t^2}\big)^{(r-1)/2}.
\end{multline*}
\end{Definition}

\begin{Lemma}
	\label{lem_formula_rho_kn}
When $k\ge s$, we have $\rho_{k,n,s}=\rho_{k-s,n-2s}$. 
\end{Lemma}

\begin{remark}
	It is straightforward to verify that $\rho_{k,n}$ given by Definition \ref{defn_rho_kn} is a polynomial in $\omega$ and $\beta$ with degree $2k$. We will see that the proof of Lemma \ref{lem_formula_rho_kn} does not need this property, so this also follows from Lemma \ref{lem_formula_rho_kn}.
\end{remark}

The rest of this sub-section is devoted to the proof of  Lemma \ref{lem_formula_rho_kn}. The argument of the proof will not be used elsewhere in the paper.

We start by proving the following lemma.
\begin{Lemma} 
	\label{lem_even_odd_cancellation_delta_series}
	In the formal series ring $\Big(\bC[\sqrt{-\beta}, \delta]/(\delta^2+\beta)\Big)[[t]]$, we have
	\begin{equation}
		\label{eqn_even_odd_cancellation_delta_series_1}
		\Big(\frac{1-t\sqrt{-\beta}}{1+t\sqrt{-\beta}}\Big)^{\frac{\delta}{4\sqrt{-\beta}} } + \Big(\frac{1-t\sqrt{-\beta}}{1+t\sqrt{-\beta}}\Big)^{-\frac{\delta}{4\sqrt{-\beta}} }
		= \frac{\sqrt{1-t\sqrt{-\beta}} + \sqrt{1+t\sqrt{-\beta}}}{(1+\beta t^2)^{1/4}}
	\end{equation}
and
	\begin{multline}
		\label{eqn_even_odd_cancellation_delta_series_2}
		-\Big(\frac{1-t\sqrt{-\beta}}{1+t\sqrt{-\beta}}\Big)^{\frac{\delta}{4\sqrt{-\beta}} } + \Big(\frac{1-t\sqrt{-\beta}}{1+t\sqrt{-\beta}}\Big)^{-\frac{\delta}{4\sqrt{-\beta}} }
		\\
		=  (t\delta) \cdot \frac{1}{(1+\beta t^2)^{1/4}} \cdot \frac{2}{\sqrt{1-t\sqrt{-\beta}} + \sqrt{1+t\sqrt{-\beta}}}
	\end{multline}
\end{Lemma}

\begin{proof}
	Suppose
	$$
	f(\sqrt{-\beta},\delta)\in \bC[\sqrt{-\beta}, \delta].
	$$
	We say that $f$ is \emph{even} in $\delta$, if $f(\sqrt{-\beta},\delta)=f(\sqrt{-\beta},-\delta)$; we say that $f$ is \emph{odd} in $\delta$, if $f(\sqrt{-\beta},\delta)=-f(\sqrt{-\beta},-\delta)$. Let 
	$$
	\psi: \bC[\sqrt{-\beta}, \delta]\to \bC[\sqrt{-\beta}, \delta]/(\delta^2+\beta)
	$$
	be the quotient map. Then if $f$ is even, we have 
	$$
	\psi(f)(\sqrt{-\beta},\delta) = \psi(f)(\sqrt{-\beta}, \sqrt{-\beta}).
	$$
	If $f$ is odd, we have
	$$
	\psi(f)(\sqrt{-\beta},\delta) = \delta\cdot \psi(f)(\sqrt{-\beta}, \sqrt{-\beta})/ \sqrt{-\beta} .
	$$
	Note that the left-hand side of \eqref{eqn_even_odd_cancellation_delta_series_1} is even in $\delta$, and the left-hand side of \eqref{eqn_even_odd_cancellation_delta_series_2} is odd in $\delta$. Therefore, in $\Big(\bC[\sqrt{-\beta}, \delta]/(\delta^2+\beta)\Big)[[t]]$, we have
	\begin{align*}
			\Big(\frac{1-t\sqrt{-\beta}}{1+t\sqrt{-\beta}}\Big)^{\frac{\delta}{4\sqrt{-\beta}} } + \Big(\frac{1-t\sqrt{-\beta}}{1+t\sqrt{-\beta}}\Big)^{-\frac{\delta}{4\sqrt{-\beta}} }
 =\,& \Big(\frac{1-t\sqrt{-\beta}}{1+t\sqrt{-\beta}}\Big)^{\frac14 } + \Big(\frac{1-t\sqrt{-\beta}}{1+t\sqrt{-\beta}}\Big)^{-\frac14 }
		\\
		= \, & \frac{\sqrt{1-t\sqrt{-\beta}} + \sqrt{1+t\sqrt{-\beta}}}{(1+\beta t^2)^{1/4}}
	\end{align*}
	and 
	\begin{align*}
		& -\Big(\frac{1-t\sqrt{-\beta}}{1+t\sqrt{-\beta}}\Big)^{\frac{\delta}{4\sqrt{-\beta}} } + \Big(\frac{1-t\sqrt{-\beta}}{1+t\sqrt{-\beta}}\Big)^{-\frac{\delta}{4\sqrt{-\beta}} }
		\\
		= \, & \delta \cdot \Bigg[-\Big(\frac{1-t\sqrt{-\beta}}{1+t\sqrt{-\beta}}\Big)^{\frac14 } + \Big(\frac{1-t\sqrt{-\beta}}{1+t\sqrt{-\beta}}\Big)^{-\frac14 }\Bigg]/\sqrt{-\beta}
		\\
		= \, & \delta \cdot \frac{-\sqrt{1-t\sqrt{-\beta}} + \sqrt{1+t\sqrt{-\beta}}}{(1+\beta t^2)^{1/4}}/\sqrt{-\beta}
		\\
		= \, & (t\delta) \cdot \frac{1}{(1+\beta t^2)^{1/4}} \cdot \frac{2}{\sqrt{1-t\sqrt{-\beta}} + \sqrt{1+t\sqrt{-\beta}}}.
	\end{align*}
\end{proof}

\begin{proof}[Proof of Lemma \ref{lem_formula_rho_kn}]
By \eqref{eqn_decompose_xi_by_rho}, we have
$$
2^{-m-1}\cdot \rho_{k,n,s}\cdot \delta_1\cdots\delta_s = \frac{1}{2^{n}}\sum_{I\subset \{1,\dots,n\}}(-1)^{I\cap\{1,\dots,s\}}\cdot \tau_I (\bar{\xi}_{k,n}).
$$
By the definition of $\bar{\xi}_{k,n}$, we have 
$$
\sum_{k=0}^\infty \bar{\xi}_{k,n} t^k = (1+\beta t^2)^{-1/2+m/2}\Big(\frac{1-t\sqrt{-\beta}}{1+t\sqrt{-\beta}}\Big)^{\frac{\omega}{2\sqrt{-\beta}}} \cdot 
\Big(\frac{1-t\sqrt{-\beta}}{1+t\sqrt{-\beta}}\Big)^{-\frac{\delta_1+\cdots+\delta_n}{4\sqrt{-\beta}}}
$$
Therefore, by the equations above and Lemma \ref{lem_even_odd_cancellation_delta_series}, 
\begin{align*}
& 2^{-m-1}\,\sum_{k=0}^\infty \rho_{k,n,s}\cdot \delta_1\cdots\delta_s\cdot t^{k} \\
= & 
 \frac{1}{2^n}(1+\beta t^2)^{-1/2+m/2}\Big(\frac{1-t\sqrt{-\beta}}{1+t\sqrt{-\beta}}\Big)^{\frac{\omega}{2\sqrt{-\beta}}}  
\\
&\quad 
\times \prod_{i=1}^s\Bigg[\Big(\frac{1-t\sqrt{-\beta}}{1+t\sqrt{-\beta}}\Big)^{-\frac{\delta_i}{4\sqrt{-\beta}}} - \Big(\frac{1-t\sqrt{-\beta}}{1+t\sqrt{-\beta}}\Big)^{\frac{\delta_i}{4\sqrt{-\beta}}}\Bigg]
\\
&\quad \times  \prod_{i=s+1}^n\Bigg[\Big(\frac{1-t\sqrt{-\beta}}{1+t\sqrt{-\beta}}\Big)^{-\frac{\delta_i}{4\sqrt{-\beta}}} + \Big(\frac{1-t\sqrt{-\beta}}{1+t\sqrt{-\beta}}\Big)^{\frac{\delta_i}{4\sqrt{-\beta}}}\Bigg]
\\
= & \frac{1}{2^n}(t\delta_1)\cdots(t\delta_s)(1+\beta t^2)^{-1/2+m/2}\Big(\frac{1-t\sqrt{-\beta}}{1+t\sqrt{-\beta}}\Big)^{\frac{\omega}{2\sqrt{-\beta}}} 
 \\
& \times 
(1+\beta t^2)^{-n/4}\cdot 2^s \cdot 
\Big(\sqrt{1-t\sqrt{-\beta}} + \sqrt{1+t\sqrt{-\beta}}\Big)^{n-2s}
\end{align*}
Recall that $n=2m+1$. Therefore,
\begin{align*}
& 2^{m}\,\sum_{k=s}^\infty \rho_{k,n,s}\cdot t^{k-s} \\
 = &
(1+\beta t^2)^{-1/2+m/2}\Big(\frac{1-t\sqrt{-\beta}}{1+t\sqrt{-\beta}}\Big)^{\frac{\omega}{2\sqrt{-\beta}}} \cdot 
(1+\beta t^2)^{-n/4}\cdot 2^s 
\\
&\quad \times  
\Big(\sqrt{1-t\sqrt{-\beta}} + \sqrt{1+t\sqrt{-\beta}}\Big)^{n-2s}
\\ 
= & 
2^s (1+\beta t^2)^{-3/4}\Big(\frac{1-t\sqrt{-\beta}}{1+t\sqrt{-\beta}}\Big)^{\frac{\omega}{2\sqrt{-\beta}}}
\big(2+2\sqrt{1+\beta t^2}\big)^{m-s}
\\
&\quad \times
 \Big(\sqrt{1-t\sqrt{-\beta}} + \sqrt{1+t\sqrt{-\beta}}\Big)
\\
= & 
2^m (1+\beta t^2)^{-3/4}\Big(\frac{1-t\sqrt{-\beta}}{1+t\sqrt{-\beta}}\Big)^{\frac{\omega}{2\sqrt{-\beta}}}
\cdot \Big(\sqrt{1-t\sqrt{-\beta}} + \sqrt{1+t\sqrt{-\beta}}\Big)
\\
&\quad \times 
\big(1+\sqrt{1+\beta t^2}\big)^{m-s}.
\end{align*}
Hence the desired result is proved.
\end{proof}

\subsection{Properties of $\rho_{k,r}$}
\label{subsec_properties_of_rho}
We establish several properties of the polynomials $\rho_{k,r}$ introduced in Definition \ref{defn_rho_kn}. 

Since $\rho_{k,r}$ is a polynomial in $\omega$ and $\beta$, it makes sense to consider the value of $\rho_{k,r}$ when plugging in $\beta=0$.
\begin{Lemma}
		\label{lem_rho_leading_term_nonzero}
	For all $k\ge 0$, we have
	$$\rho_{k,r}|_{\beta = 0} = (-1)^k \frac{ 2^{(r+1)/2} }{k!} \cdot \omega^k\neq 0.$$
\end{Lemma}

\begin{proof}
Note that 
$$
\Big(\frac{1-t\sqrt{-\beta}}{1+t\sqrt{-\beta}}\Big)^{\frac{\omega}{2\sqrt{-\beta}}} \Bigg|_{\beta=0} = \exp\Bigg(\frac{\ln(1-t\sqrt{-\beta}) - \ln(1+t\sqrt{-\beta})}{2t\sqrt{-\beta}}\cdot \omega t\Bigg)\Bigg|_{\beta=0} = e^{-\omega t}.
$$
Therefore the lemma follows from Definition \ref{defn_rho_kn}.
\end{proof}

In particular, Lemma \ref{lem_rho_leading_term_nonzero} implies that 
\begin{equation}
	\label{eqn_rho_kr_nonzero}
\rho_{k,r}\neq 0.
\end{equation}

\begin{Lemma}
	\label{lem_power_series_coef_aks_linear_independent}
	Suppose $k,s$ are non-negative integers.
	Define $a_{k,s}\in \mathbb{R}$ by the following series expansion near $x=0$:
	\begin{equation}
		\label{eqn_Taylor_expansion_sqrt(1+x)}
		\sum_{k=0}^\infty a_{k,s} x^k = (1+\sqrt{1+x})^s.
	\end{equation}
	For each non-negative integer $N$, consider the square matrix $(a_{k,s})_{0\le k,s\le N}$. 
	Then we have $$\det(a_{k,s})_{0\le k,s\le N}\neq 0.$$
\end{Lemma}

\begin{proof}
	Note that the radius of convergence of \eqref{eqn_Taylor_expansion_sqrt(1+x)} is positive.
Assume we have $\det(a_{k,s})_{0\le k,s\le N}= 0$ for some $N$, then there exists a sequence of real numbers $b_0,\dots,b_N$, such that 
\begin{equation}
	\label{eqn_limit_from_det(aks)=0}
\lim_{z\to 0} \frac{\sum_{s=0}^N b_s (1+\sqrt{1+z})^s}{z^N}= 0,
\end{equation}
 and that $b_s\neq 0$ for some $s$, where the limit in \eqref{eqn_limit_from_det(aks)=0} is taken in $\bC$. 
 
Assume the polynomial  $\sum_{s=0}^N b_s z^s$ factorizes in $\bC$ as
 $$\sum_{s=0}^N b_s z^s=c \prod_{j=1}^d (z-c_i)$$
 for $d\le N$, $c\neq 0$. Then we have
$$
\lim_{z\to 0} \frac{\prod_{j=1}^d (\sqrt{1+z}+1-c_i)}{z^N}=0,
$$
where the limit is taken in $\bC$.
Since $\frac{d}{dz}|_{z=0}(\sqrt{1+z})\neq 0$, each factor $(\sqrt{1+z}+1-c_i)$ is either non-zero or has a simple zero at $z=0$. Since $d\le N$, we have a contradiction.
\end{proof}

\begin{Lemma}
	\label{lem_linear_independency_rho}
	For each pair of integers $(k,r)$ such that $k\ge 0$ and $r$ is odd, the polynomials
	$\rho_{k,r}$, $\rho_{k,r+2}$, $\dots$, $\rho_{k,r+2[k/2]}$, are linearly independent over $\bC$.
\end{Lemma}
\begin{proof}
	Let $$F_r(t) = \sum_{k=0}^\infty \rho_{k,r}\cdot t^k.$$ 
	Then by Definition \ref{defn_rho_kn}, we have $F_{r+2s}(t) = (1+\sqrt{1+\beta t^2})^s F_{r}(t)$. 
	Let $a_{k,s}$ be as in Lemma \ref{lem_power_series_coef_aks_linear_independent}. Then for $j\le [k/2]$, we have 
	$$
	\rho_{k,r+2j} = \sum_{i=0}^{[k/2]} a_{i,j}\cdot \rho_{k-2i,r}\cdot \beta^i.
	$$
	Therefore, by Lemma \ref{lem_power_series_coef_aks_linear_independent}, we only need to show that the polynomials of the form $\beta^i\cdot \rho_{k-2i,r}$ are linearly independent for $i=0,\dots,[k/2]$. 
  This follows immediately from Lemma \ref{lem_rho_leading_term_nonzero}. 
\end{proof}

\subsection{The ideals modulo $\beta$}
\label{subsec_ideals_modulo_beta}
Now we study the algebraic properties of the ring
$$
\bar{R}_n/(\beta) = \bC[\alpha,\delta_1,\dots,\delta_n]/(\delta_1^2,\dots,\delta_n^2).
$$
It is straightforward to verify that the ideal in $(\beta)\subset\bar{R}_n$ is invariant under flip symmetries, and the ring $\bar{R}_n/(\beta)$ satisfies \eqref{eqn_prod_delta_i_notin_I}. Therefore, it makes sense to refer to flip symmetries on $\bar{R}_n/(\beta)$, and each flip symmetry has a well-defined parity.

\begin{Lemma}
	\label{lem_linear_independency_alpha^s+flip}
	Suppose $s\ge m$. Then the images of $\alpha^s$ under all even flip symmetries are linearly independent in $\bC[\alpha,\delta_1,\dots,\delta_n]/(\delta_1^2,\dots,\delta_n^2)$.
\end{Lemma}

\begin{proof}
	We have
	$$
	\alpha^s = (\omega- \frac12 \sum \delta_i)^s = \sum_{I\subset \{1,\dots,n\}, |I|\le s}\delta^I \cdot \omega^{s-|I|}\cdot (-1/2)^{|I|} \cdot \frac{s!}{(s-|I|)!}.
	$$
	Assume $I\subset \{1,\dots,n\}$ satisfies $|I|\le m$.  Define $\theta_I\in \bC[\alpha,\delta_1,\dots,\delta_n]/(\delta_1^2,\dots,\delta_n^2)$ by the following equation
	\begin{equation}
		\label{eqn_def_theta_I}
		\theta_I = 
		\begin{cases}
			\delta^I \cdot \omega^{s-|I|}\cdot (-\frac12)^{|I|} \cdot   \frac{s!}{(s-|I|)!}& \text{ if } |I^c|>s, \\
			\delta^I \cdot \omega^{s-|I|}\cdot (-\frac12)^{|I|} \cdot \frac{s!}{(s-|I|)!} + \delta^{I^c} \cdot \omega^{s-|I^c|}\cdot (-\frac12)^{|I^c|}   \cdot \frac{s!}{(s-|I^c|)!}&\text{ if } |I^c|\le s.
		\end{cases}
	\end{equation}
	Note that since $|I|\le m \le s$, Equation \eqref{eqn_def_theta_I} indeed defines a polynomial.
	
	The set $\{\theta_I\}_{|I|\le m}$ has $2^{n-1}$ elements and is linearly independent over $\bC$, and we have
	$$
	\alpha^s = \sum_{I\subset \{1,\dots,n\}, ~ |I| \le m} \theta_I,
	$$
	which implies 
	\begin{equation}
		\label{eqn_theta_I_in_terms_of_alpha^s}
		\theta_I = \frac{1}{2^{n-1}} \sum_{J\subset \{1,\dots,n\}, ~ |J| \text{ is even}}(-1)^{|I\cap J|}\,\tau_J (\alpha^s).
	\end{equation}
	Therefore, the linear space spanned by all $\tau_J(\alpha^s)$ with $|J|$ even contains all $\theta_I$ with $|I|\le m$, and hence it has dimension at least $2^{n-1}$. So the desired result is proved.
\end{proof}

\begin{Lemma}
	\label{lem_ideal_mod_beta_full_not_minimal}
	Suppose $s\ge m$, let $\mathfrak{a}$ be the ideal in $\bC[\alpha,\delta_1,\dots,\delta_n]/(\delta_1^2,\dots,\delta_n^2)$ generated by $\alpha^s$ and all its even flip symmetries. Suppose 
	$$ f\in \bC[\alpha,\delta_1,\dots,\delta_n]/(\delta_1^2,\dots,\delta_n^2)$$ is a homogeneous polynomial with degree strictly greater than $2s$. Then $f\in \mathfrak{a}$. 
\end{Lemma}
\begin{proof}
	If $I\subset \{1,\dots,n\}$ and $|I|\le m$, let $\theta_I$ be defined as in  \eqref{eqn_def_theta_I}.  By \eqref{eqn_theta_I_in_terms_of_alpha^s}, $\theta_I\in \mathfrak{a}$.
	
	We only need to consider the case when $f$ has the form 
	$$
	f = \omega^{u}\delta^J.
	$$
	with $J\subset\{1,\dots,n\}$ and $u+|J|>s$. We may also assume without loss of generality that $u+|J|=s+1$. We show that $f\in \mathfrak{a}$ by discussing the following cases:

	Case 1: $1\le |J|\le m+1$. Take $j\in J$, and let $J' = J\backslash \{j\}$, then 
	$$
	f = \theta_{J'} \cdot \delta_j\cdot (-2)^{|J'|}\cdot \frac{(s-|J'|)!}{s!}\in \mathfrak{a}.
	$$
	
	Case 2: $m+1<|J|\le n$. By the assumptions, $|J|-1\le s$. Take $j\in J$, let $J'=J^c\cup\{j\}$, then $|J'| = n - |J| + 1\le m$, and $|(J')^c| = |J|-1 \le s$. Hence we have 
	$$
	f = \theta_{J'}\cdot \delta_j\cdot (-2)^{|J|-1}\cdot \frac{(s-|J|+1)!}{s!}\in \mathfrak{a}.
	$$
	
	Case 3: $|J|=0$, $s<n$. In this case, we have $J=\emptyset$, and 
	$$
	f = \omega^{s+1} = \theta_{\emptyset} \cdot \omega \in \mathfrak{a}.
	$$
	
	Case 4: $|J|=0$, $s\ge n$. In this case, 
	$$
	\theta_\emptyset= \omega^s + \delta_1\cdots\delta_n\cdot \omega^{s-n} \cdot(-\frac12)^n \cdot \frac{s!}{(s-n)!}.
	$$
	By Case 2, we have $\delta_1\cdots\delta_n\cdot \omega^{s-n+1} \in\mathfrak{a}.$
	Therefore 
	$$
	f = \omega^{s+1} = \omega\cdot 	\theta_\emptyset - \delta_1\cdots\delta_n\cdot \omega^{s-n+1} \cdot(-\frac12)^n \cdot \frac{s!}{(s-n)!}\in\mathfrak{a}.
	\phantom\qedhere\makeatletter\displaymath@qed
	$$
\end{proof}

\begin{Definition}
	Define $\bar{K}_{g,n}$ to be the image of $K_{g,n}$ in $\bar{R}_n/(\beta)$. 
	Define $\bar{K'}_{g,n}$ to be the image of $K_{g,n}'$ in $\bar{R}_n/(\beta)$. 
\end{Definition}

By definition, $\bar{K'}_{g,n}$ is the ideal in $\bar{R}_n/(\beta)$ generated by $\alpha^{2(g+m)}$ and all of its even flip symmetries.

We also use $\bar{K}_{g,n}^{\langle i \rangle}$ and $\bar{K'}_{g,n}^{\langle i \rangle}$ to denote the respective components of $\bar{K}_{g,n}$ and $\bar{K'}_{g,n}$ with degree $2(g+m+i)$. 

\begin{Lemma}
	\label{lem_Kgn=Kgn'_modulo_beta}
	We have $\bar{K}_{g,n} = \bar{K'}_{g,n}$, and 
	$$\dim_\bC \bar{K}_{g,n}^{\langle i \rangle} = \dim_\bC \bar{K'}_{g,n}^{\langle i \rangle} = 
	\begin{cases}
		0  &\text{ if } i< 0 \\
		2^{n-1}&\text{ if } i= 0 \\
		2^n - {\displaystyle \sum_{j= 0}^{m-g-i}{n\choose j}}&\text{ if } i>0.
	\end{cases}
	$$
\end{Lemma}
\begin{remark}
	By definition, $\sum_{j= 0}^{m-g-i}{n\choose j}=0$ if $m-g-i<0$. We also have ${n\choose k}=0$ whenever $k<0$ or $k>n$. 
\end{remark}
\begin{proof}
	If $i< 0$, by Lemma \ref{lem_Poincare_Kgn}, we have $\dim_\bC \bar{K}_{g,n}^{\langle i \rangle}\le \dim_\bC K_{g,n}^{\langle i \rangle}=0$.  By \eqref{eqn_Kgn_subset_Kgn'}, we have $\bar{K}_{g,n}'\subset \bar{K}_{g,n}$, which implies $\dim_\bC \bar{K}_{g,n}^{\prime \, \langle i \rangle}=0$.
	
	If $i=0$, by \eqref{eqn_dim_K_gn^0}, we have
	$\dim_\bC \bar{K}_{g,n}^{\langle 0\rangle}\le \dim_\bC K_{g,n}^{\langle 0 \rangle}=2^{n-1}.$ 
	By Lemma \ref{lem_linear_independency_alpha^s+flip}, we have $\dim_\bC \bar{K}_{g,n}^{\prime\,\langle 0 \rangle} \ge 2^{n-1}$. By \eqref{eqn_Kgn_subset_Kgn'}, we have $\bar{K'}_{g,n}^{\langle 0\rangle}\subset \bar{K}_{g,n}^{\langle 0 \rangle}$. Therefore, all the above inequalities must achieve equality, and we have  $\bar{K}_{g,n}^{\langle 0 \rangle}=\bar{K'}_{g,n}^{\langle 0 \rangle}$ with their dimension being equal to $2^{n-1}$. 
	
	If $i>0$, by Lemma \ref{lem_ideal_mod_beta_full_not_minimal} and the definition of $\bar{K'}_{g,n}$, we have $\bar{K'}_{g,n}^{\langle i \rangle} $ is equal to the component of $\bC[\alpha,\delta_1,\dots,\delta_n]/(\delta_1^2,\dots,\delta_n^2)$ with degree $2g+2m+2i$.
	Since $\bar{K'}^{\langle i \rangle}_{g,n}\subset \bar{K}^{\langle i \rangle}_{g,n}$, we must have $\bar{K}_{g,n}^{\langle i \rangle}  = \bar{K'}_{g,n}^{\langle i \rangle}$, and their dimension is equal to
	$$
	 \sum_{j=0}^{g+m+i}{n\choose j} = 2^n - \sum_{j=0}^{m-g-i}{n\choose j}.
	\phantom\qedhere\makeatletter\displaymath@qed
	$$
\end{proof}

\subsection{Completion of the proof of Theorem \ref{thm_generators_of_Ign}}
\label{subsec_completion_of_proof_Ign_generators}
We now prove Proposition  \ref{prop_generators_of_Kgn}, which will finish the proof of Theorem \ref{thm_generators_of_Ign}.
We first set up some notation.

\begin{Notation*}
In Section \ref{subsec_completion_of_proof_Ign_generators}, $i,j$ will always denote non-negative integers, and $I$ will denote subsets of $\{1,\dots,n\}$.

Recall that for $s\in \mathbb{Z}$, we use $K_{g,n}^{\langle s \rangle}$ to denote the component of $K_{g,n}$ with degree $2(g+m+s)$. We will also use $K_{g,n}^{\prime \, \langle s \rangle}$ to denote the component of $K_{g,n}'$ with degree $2(g+m+s)$. 
\end{Notation*}

By Lemma \ref{lem_Poincare_Kgn}, we have $K_{g,n}^{\langle s \rangle}=0$ when $s<0$. Since $K_{g,n}'\subset K_{g,n}$, we also have $K_{g,n}^{\prime \, \langle s \rangle}=0$ when $s<0$. The next few lemmas explicitly construct several families of elements in $K_{g,n}'$ using the polynomials $\rho_{k,r}$ introduced in Section \ref{subsec_decompose_bar_xi}.

\begin{Lemma}
	\label{lem_isolated_products_in_lowest_degree_K'gn}
	Assume $g+|I|\le m$. Then for all $I\subset\{1,\dots,n\}$ with $|I|=s$, we have
	$$
	\rho_{g+m-|I|,n-2|I|}\,\delta^I \in K_{g,n}^{\prime \,\langle 0\rangle}.
	$$
\end{Lemma}

\begin{proof}
Since  $g+|I|\le m$, we have $|I^c|>g+m$, so
$\rho_{g+m,n,n-s}=0$. By \eqref{eqn_decompose_xi_by_rho},
$$
2^{-m-1}\cdot \rho_{g+m,n,|I|} \cdot \delta^I = \frac{1}{2^{n-1}}\sum_{J\subset \{1,\dots,n\}, |J| \text{ even }}(-1)^{|I\cap J|}\cdot \tau_J(\bar{\xi}_{g+m,n}) \in K_{g,n}^{\prime \,\langle 0\rangle}.
$$
By Lemma \ref{lem_formula_rho_kn}, we have 	$\rho_{g+m-|I|,n-2|I|} = \rho_{g+m,n,|I|}$.  
\end{proof}

\begin{Lemma}
	\label{lem_delta^I*beta^i_from_lowest_degree}
	Assume  $g+i+|I|\le m$. We have
	$$
	\rho_{g+m-i-|I|,n-2|I|-2i}\,\delta^I \beta^{i} \in K_{g,n}^{\prime \, \langle i \rangle}.
	$$
\end{Lemma}
\begin{proof}
	By the assumption, $i+|I|\le m <n$. Let $a_1,\dots,a_i$ be distinct elements in $I^c$, and let $J=I\cup\{a_1,\dots,a_i\}$. Then by Lemma \ref{lem_isolated_products_in_lowest_degree_K'gn}, 
	$$
	\rho_{g+m-|J|,n-2|J|}\,\delta^J\in K_{g,n}^{\prime \, \langle 0 \rangle}.
	$$
	Multiplying the left-hand side by $\delta_{a_1}\cdots\delta_{a_i}$ yields 
	$$
	\rho_{g+m-|J|,n-2|J|}\delta^I\beta^i \in K_{g,n}^{\prime \, \langle i \rangle}.
	\phantom\qedhere\makeatletter\displaymath@qed
	$$
\end{proof}

\begin{Lemma}
		\label{lem_delta^I*beta^i_from_lower_g}
	Assume $g+i+|I|\le m$ and $j\le \min\{i,g\}$. We have
$$
\rho_{g+m-i-|I|,n-2i-2|I|+2j}\,\delta^I \beta^{i} \in K_{g,n}^{\prime \, \langle i \rangle}.
$$
\end{Lemma}
\begin{proof}
	By Lemma \ref{lem_K'gn_relation_in_g}, $\beta^j\cdot K_{g-j,n}^{\prime \, \langle i-j \rangle}\subset K_{g,n}^{\prime \, \langle i \rangle}$.
	Applying Lemma \ref{lem_delta^I*beta^i_from_lowest_degree} with $g$ replaced with $g-j$ and $i$ replaced with $i-j$, we have 
	$$
	\rho_{g+m-i-|I|,n-2i-2|I|+2j}\,\delta^I \beta^{i-j} \in K_{g-j,n}^{\prime \, \langle i-j \rangle}.
	$$
	Therefore the result is proved.
\end{proof}

\begin{Lemma}
	\label{lem_rho_linearly_independent_for_j<=g}
	Assume $g+i+|I|\le m$, then the polynomials in $\bC[\omega,\beta]$ given by 
	$$\rho_{g+m-i-|I|,n-2i-2|I|+2j}, \qquad  j \le g$$ 
	are linearly independent. 
\end{Lemma}

\begin{proof}
	Note that since $g+i+|I|\le m$, we have $2g \le g+m-i-|I|$. Therefore the desired result follows from Lemma \ref{lem_linear_independency_rho}.
\end{proof}

Recall that in Section \ref{subsec_isotypic_decomposition}, we introduced a linear space decomposition 
$$
 \bar{R}_n = \bigoplus_{|I|\le m} \bar{R}_n^I.
$$
This is the isotypic decomposition of $\bar{R}_n$ with respect to the even flip symmetry actions. The isotypic decompositions of $K_{g,n}$ and $K_{g,n}'$ by even flip symmetry actions are given by \eqref{eqn_Kgn_Kgn'_isotypic_decompose}.

Proposition \ref{prop_generators_of_Kgn} is implied by the following result.

\begin{Proposition}
	\label{prop_proof_Kgn=Kgn'}
	For all $i \ge 0$, we have
	\begin{enumerate}
		\item $K_{g,n}^{\prime \,\langle i \rangle } = K_{g,n}^{\langle i \rangle }$.
		\item If $g+i+|I|\le m$, then $K_{g,n}^{\prime \,\langle i \rangle } \cap \beta^i \bar{R}_n^I$ is the linear space spanned by the elements of the  form:
        \begin{equation}
        	\label{eqn_basis_of_Ki_cap_beta^iRI}
        	     \rho_{g+m-i-|I|,n-2i-2|I|+2j}\,\delta^I \beta^{i} \quad \text{ with } j\le \min \{i,g\}.
        \end{equation}
	\end{enumerate}
\end{Proposition}

\begin{remark}
	The notation $\beta^i \bar{R}_n^I$ denotes the linear space consisting of all elements of the form $\beta^i\cdot x$ with $x\in \bar{R}_n^I$. It is not an ideal.
\end{remark}

\begin{remark}
	\label{rem_linear_independency_of_constructed_elements}
	By Lemma \ref{lem_delta^I*beta^i_from_lower_g}, the polynomials listed in \eqref{eqn_basis_of_Ki_cap_beta^iRI} are elements of $K_{g,n}^{\prime \,\langle i \rangle } \cap \beta^i \bar{R}_n^I$. By Lemma \ref{lem_rho_linearly_independent_for_j<=g}, these elements are  linearly independent.
\end{remark}

\begin{proof}
We prove the proposition by induction on the pair $(i,g)$. The proof is divided into three steps: we first show that the proposition holds when $i=0$ or $g=0$; then, we show that for a fixed value of $n$ and given $i>0, g>0$, if the proposition holds for $(i-1,g-1)$, then it also holds for $(i,g)$.\\

\textbf{Step 1: $i=0$.}
In this case, by \eqref{eqn_dim_K_gn^0} and Lemma \ref{lem_Kgn=Kgn'_modulo_beta}, we have
$$
\dim_\bC K_{g,n}^{\langle 0 \rangle} = 2^n = \dim_\bC \bar{K'}_{g,n}^{\langle 0 \rangle} \le  \dim_\bC K_{g,n}^{\prime \, \langle 0 \rangle} \le \dim_\bC K_{g,n}^{\langle 0 \rangle}.
$$
Therefore, all the above inequalities must achieve equalities, and we have $K_{g,n}^{\prime \,\langle 0 \rangle } = K_{g,n}^{\langle 0 \rangle }$.

To verify the second part of the proposition, note that $K_{g,n}^{\prime \,\langle 0 \rangle }$ is the linear space spanned by the images of $\bar{\xi}_{g,n}$ under even flip symmetries. By \eqref{eqn_Kgn_Kgn'_isotypic_decompose} and \eqref{eqn_decompose_xi_by_rho}, we see that $K_{g,n}^{\prime \,\langle 0 \rangle }$ is spanned by the polynomials of the form
$$
\rho_{g+m,n, |I|}\cdot \delta^I + \rho_{g+m,n,|I^c|} \cdot \delta^{I^c},\quad \text{ with } |I|\le m.
$$

If $g+|I|\le m$, we have $|I^c|>g+m$, thus $\rho_{g+m,n,|I^c|} = 0$. By Lemma \ref{lem_formula_rho_kn}, $\rho_{g+m,n, |I|} = \rho_{g+m-|I|, n-2|I|}$. Hence $K_{g,n}^{\prime \,\langle 0 \rangle }\cap \bar{R}_n^I$ is a $1$--dimensional linear space spanned by $\rho_{g+m-|I|, n-2|I|}\delta^I$. \\

\textbf{Step 2: $g=0$.} In this case, the first part of the proposition was proved in \cite{kronheimer2022relations}, see Remark \ref{rem_Ign_g=0}. We prove the second part of the proposition.  Assume $i, I$ are given such that $i+|I|\le m$. 

By \eqref{eqn_recursive_xi}, we have $\bar{\xi}_{m+1,n} =(\alpha\,\bar{\xi}_{m,n})/(m+1)$. Therefore the ideal $K'_{0,n}$ is generated by $\bar{\xi}_{m,n}$ and all its even flip symmetries.  
By \eqref{eqn_Kgn_Kgn'_isotypic_decompose} and \eqref{eqn_decompose_xi_by_rho}, the ideal $K_{0,n}'$ is generated by polynomials of the form
$$
\rho_{m,n, |J|}\cdot \delta^J + \rho_{m,n,|J^c|} \cdot \delta^{J^c}
$$
for all $J$ such that $|J|\le m$. Since $|J^c| = n-|J|>m$, we have $\rho_{m,n,|J^c|} =0$.  Since $\rho_{m,n, |J|} = \rho_{m-|J|,n-2|J|}$, we conclude that the ideal $K_{0,n}'$ is generated by polynomials of the form
$$
\rho_{m-|J|,n-2|J|}\cdot \delta^J
$$
for all $J$ such that $|J|\le m$. 

For each $s=0,\dots,i$, take $J$ to be a set with $|I|+i$ elements such that $I\subset J$. By the assumptions, we have $|J|\le m$, therefore $\rho_{m-|J|,n-2|J|}\cdot \delta^J\in K_{0,n}'$. Multiplying this element by $\omega^{i-s}\delta^{J\backslash I}$, we have 
\begin{equation}
	\label{eqn_delta^Iw^(i-s)beta^s_factor_in_K0n^iI}
\rho_{m-|I|-s,n-2|I|-2s}\cdot \delta^I\omega^{i-s} \beta^s \in K_{0,n}'.
\end{equation}
Computing the gradings of the left-hand side of \eqref{eqn_delta^Iw^(i-s)beta^s_factor_in_K0n^iI}, we see that it is an element in $K_{0,n}^{\prime \, \langle i \rangle}$. Therefore, we have constructed $i+1$ elements in $K_{0,n}^{\prime \, \langle i \rangle}\cap \bar{R}_n^I:$
\begin{equation}
	\label{eqn_basis_K0n'_intersects_RI}
\rho_{m-|I|-s,n-2|I|-2s}\cdot \delta^I\omega^{i-s} \beta^s,\quad s = 0,1,\dots,i.
\end{equation}
By Lemma \ref{lem_rho_leading_term_nonzero}, these elements are linearly independent. By Lemma \ref{lem_g=0_isotypic_decomposition_Kgn}, $\dim_\bC K_{0,n}^{\langle i \rangle}\cap \bar{R}_n^I=i+1$, therefore \eqref{eqn_basis_K0n'_intersects_RI} is a basis of $K_{0,n}^{\prime \, \langle i \rangle}\cap \bar{R}_n^I.$ By Lemma \ref{lem_rho_leading_term_nonzero} again, we see that
$$
K_{0,n}^{\prime \, \langle i \rangle}\cap \beta^i \bar{R}_n^I
$$
is the $1$--dimensional linear space spanned by $\rho_{m-|I|-i,n-2|I|-2i}\cdot \delta^I\beta^{i}$.\\

\textbf{Step 3: $i>0$, $g>0$, and assume the proposition holds for $(i-1,g-1)$.}  
We show that the proposition holds for $(i,g)$. 
Assume $I$ is a subset of $\{1,\dots,n\}$ such that $g+i+|I|\le m$. If $g+i>m$, then no such $I$ exists, and the discussions below about $I$ will be vacuously true.

By the induction hypothesis, $K_{g-1,n}^{\prime\, \langle i-1\rangle}\cap \beta^{i-1} \bar{R}_n^I$ is the linear space spanned by 
$$
 \rho_{g+m-i-|I|,n-2i-2|I|+2+2j}\,\delta^I \beta^{i-1}, \qquad 0\le j \le \min\{g-1,i-1\}.
$$
Therefore, $\beta K_{g-1,n}^{\prime\, \langle i-1\rangle}\cap \beta^i \bar{R}_n^I$ is the linear space spanned by 
\begin{equation}
	\label{eqn_basis_in_beta^i_RI_from_g-1}
 \rho_{g+m-i-|I|,n-2i-2|I|+2j}\,\delta^I \beta^{i}, \qquad 1\le j \le \min\{g,i\}.
\end{equation}
Hence by Lemma \ref{lem_rho_linearly_independent_for_j<=g} (see also Remark \ref{rem_linear_independency_of_constructed_elements}), we have 
\begin{equation}
	\label{eqn_base_element_in_Kgn^i_cap_RI_not_from_g-1}
\rho_{g+m-i-|I|,n-2i-2|I|}\,\delta^I \beta^{i}  \notin \beta K_{g-1,n}^{\prime\, \langle i-1\rangle}\cap \beta^i \bar{R}_n^I.
\end{equation}
By Lemma \ref{lem_delta^I*beta^i_from_lowest_degree}, 
\begin{equation}
	\label{eqn_base_element_in_Kgn^i_cap_RI}
	\rho_{g+m-i-|I|,n-2i-2|I|}\,\delta^I \beta^{i} \in K_{g,n}^{\prime\, \langle i \rangle}\cap \beta^i \bar{R}_n^I.
\end{equation}

By Lemma \ref{lem_K'gn_relation_in_g}, we have $\beta K_{g-1,n}^{\prime\, \langle i-1\rangle}\subset K_{g,n}^{\prime \,\langle i \rangle }$, thus 
$$\beta K_{g-1,n}^{\prime\, \langle i-1\rangle}\cap \beta^i \bar{R}_n^I \subset K_{g,n}^{\prime \,\langle i \rangle }\cap \beta^i \bar{R}_n^I.$$
Therefore, \eqref{eqn_base_element_in_Kgn^i_cap_RI_not_from_g-1} and \eqref{eqn_base_element_in_Kgn^i_cap_RI} imply
\begin{equation}
	\label{eqn_intersect_beta^i_RI_strict_inclusion}
\beta K_{g-1,n}^{\prime\, \langle i-1\rangle}\cap \beta^i \bar{R}_n^I \subsetneqq K_{g,n}^{\prime \,\langle i \rangle }\cap \beta^i \bar{R}_n^I.
\end{equation}

Since $\beta^i \bar{R}_n^I=(\beta^i)\cap  \bar{R}_n^I$,  
 \eqref{eqn_intersect_beta^i_RI_strict_inclusion} can be rewritten as 
 $$
 (\beta^i)\cap \beta K_{g-1,n}^{\prime\, \langle i-1\rangle}\cap \bar{R}_n^I \subsetneqq (\beta^i)\cap K_{g,n}^{\prime \,\langle i \rangle }\cap \bar{R}_n^I,
 $$
 which implies
$$
\beta K_{g-1,n}^{\prime\, \langle i-1\rangle}\cap \bar{R}_n^I =(\beta)\cap  \beta K_{g-1,n}^{\prime\, \langle i-1\rangle}\cap \bar{R}_n^I\subsetneqq (\beta)\cap K_{g,n}^{\prime \,\langle i \rangle }\cap  \bar{R}_n^I,
$$
therefore
$$
\dim_\bC \frac{(\beta)\cap K_{g,n}^{\prime \,\langle i \rangle }\cap  \bar{R}_n^I}{\beta K_{g-1,n}^{\prime\, \langle i-1\rangle}\cap \bar{R}_n^I}\ge 1.
$$

Since the above argument works for every $I$ with $|I|\le m-g-i$, by \eqref{eqn_Kgn_Kgn'_isotypic_decompose}, we have
\begin{equation}
	\label{eqn_dim_diff_Kgni_minus_beta_K(g-1)_lower_bound}
\dim_\bC \frac{(\beta)\cap K_{g,n}^{\prime \,\langle i \rangle }}{\beta K_{g-1,n}^{\prime\, \langle i-1\rangle} }\ge \sum_{j=0}^{m-g-i}{n\choose j}.
\end{equation}
Note that the right-hand side of \eqref{eqn_dim_diff_Kgni_minus_beta_K(g-1)_lower_bound} is zero if $g+i>m$. 

By the induction hypothesis, Lemma \ref{lem_Kgn=Kgn'_modulo_beta}, Lemma \ref{lem_graded_dim_Kgn_diff}, and the short exact sequence
\begin{equation}
	\label{eqn_short_exact_sequence_Kgn'_mod_beta}
0 \to (\beta)\cap K_{g,n}^{\prime\, \langle i \rangle} \to K_{g,n}^{\prime\, \langle i \rangle} \to \bar{K'}_{g,n}^{\langle i \rangle} \to 0,
\end{equation}
we have
\begin{align*}
\dim_\bC K_{g,n}^{\prime\, \langle i \rangle} & = \dim_\bC \, (\beta)\cap K_{g,n}^{\prime\, \langle i \rangle}  + \dim_\bC \bar{K'}_{g,n}^{\langle i \rangle}
\\
& \ge \sum_{j=0}^{m-g-i}{n\choose j} + \dim_\bC \beta K_{g-1,n}^{\prime\, \langle i-1\rangle} + \dim_\bC \bar{K'}_{g,n}^{\langle i \rangle}
\\
& =  \sum_{j=0}^{m-g-i}{n\choose j} + \dim_\bC  K_{g-1,n}^{\prime\, \langle i-1\rangle} + \dim_\bC \bar{K'}_{g,n}^{\langle i \rangle}
\\
& =  \sum_{j=0}^{m-g-i}{n\choose j} + \dim_\bC  K_{g-1,n}^{\langle i-1\rangle} + \dim_\bC \bar{K'}_{g,n}^{\langle i \rangle}
\\
& = \dim_\bC  K_{g-1,n}^{\langle i-1\rangle} + 2^n 
\\
& = \dim_\bC K_{g,n}^{\langle i \rangle},
\end{align*}
where the first line follows from the short exact sequence \eqref{eqn_short_exact_sequence_Kgn'_mod_beta}, the second line follows from \eqref{eqn_dim_diff_Kgni_minus_beta_K(g-1)_lower_bound}, the third line follows the fact that $\beta$ is a non-zero divisor in $\bar{R}_n$, the fourth line follows from the induction hypothesis, the fifth line follows from Lemma \ref{lem_Kgn=Kgn'_modulo_beta}, and the last line follows from Lemma \ref{lem_graded_dim_Kgn_diff}.

Since $K_{g,n}^{\prime\, \langle i \rangle}  \subset K_{g,n}^{\langle i \rangle}$, we must have $K_{g,n}^{\prime\, \langle i \rangle}  = K_{g,n}^{\langle i \rangle}$, and all the inequalities above achieve equalities. 
This proves the first part of the proposition. For the second part, since all inequalities above achieve equalities, the dimension difference between the two sides of \eqref{eqn_intersect_beta^i_RI_strict_inclusion} must be $1$. Therefore, $K_{g,n}^{\prime \,\langle i \rangle }\cap \beta^i \bar{R}_n^I$ is spanned by $\beta K_{g-1,n}^{\prime\, \langle i-1\rangle}\cap \beta^i \bar{R}_n^I $ and $\rho_{g+m-i-|I|,n-2i-2|I|}\,\delta^I \beta^{i}$. Since a basis of $\beta K_{g-1,n}^{\prime\, \langle i-1\rangle}\cap \beta^i \bar{R}_n^I $ is given by \eqref{eqn_basis_in_beta^i_RI_from_g-1}, the second part of the proposition is proved.
\end{proof}

\section{The sub-leading term}
\label{sec_subleading}
The ideal $I^0_{g,n}$ provides the leading terms for polynomials in $J_{g,n}$ according to Proposition \ref{prop_I_realized_by_J}
and Proposition \ref{prop_J_realized_by_I}. In this section, we will derive a formula for the \emph{sub-leading} terms of the minimal-degree relations in $J_{g,n}$. This generalizes the analogous result for the genus $0$ case given by \cite{kronheimer2022relations}*{Proposition 5.3}. We will use a wall-crossing argument that follows the strategy of
 \cite{kronheimer2022relations}*{Section 5}. In fact, most of the arguments in \cite{kronheimer2022relations} extend to the positive-genus case without essential change. The new difficulty is the appearance of a non-trivial Jacobian, which will be dealt with  in Section \ref{subsec_M0M1_subleading}.

In this section, we will write $N_{g,n}^0$ as $N_{g,n}$, and we will use the notation $d$ to denote the degrees of line bundles. We use $\Sigma_{g,n}$ to denote a closed oriented surface with genus $g$ and $n$ labeled points $\{p_1,\dots,p_n\}$, and use $\Sigma_g$ to denote a closed oriented surface with genus $g$. We will use $\eta$ to denote subsets of $\{p_1,\dots,p_n\}$.

\subsection{The Mumford relations}\label{subsec_construct_w0w1}
We briefly review the computation of Mumford relations from \cite{XZ:excision}*{Section 3}. The only difference between the computation below and the one in \cite{XZ:excision} is that here, we compute the Mumford relations in a slightly more general setting by considering all possible choices of $\eta$. In fact, it only differs from the results in \cite{XZ:excision}*{Section 3} by flip symmetries. 

Nevertheless, we sketch a direct argument without resorting to flip symmetries because flip symmetries no longer exist for Floer homology with local coefficients, so this will be necessary for the discussions in Section \ref{sec_local_coef}. We  also deduce a new formula \eqref{eqn_chern_class_from_dual_bundle} along the way, which will be important later.

Let $(\widetilde{\mathbb{E}},\widetilde{\mathbb{L}} )\to \widetilde{N}_{g,n}\times \Sigma_{g,n}$ be the universal family of stable parabolic bundles of degree 
0 over $\Sigma_{g,n}$, where the weights at labeled points are given by \cite{XZ:excision}*{Definition 2.5}.
Given $\eta$,
suppose $L_0$ is a parabolic line bundle of degree $d$ over $\Sigma_{g,n}$ with weight $\frac14$ at all $p\in \eta$ and 
$-\frac14$ at all $p\in \eta'$. Then we have 
$$
\parr\deg L_0=d+\frac14 |\eta| -\frac14 |\eta'|=d+ \frac12 |\eta|-\frac14 n=d+ \frac12 |\eta|-\frac12 m -\frac14.
$$
Unless otherwise specified, we will always assume that $d$ and $\eta$ are chosen so that $\parr\deg L_0=\frac14$. 
Notice that this means $|\eta|$ must be odd when $m$ is even, and $|\eta|$ must be even when $m$ is odd.
When we need to emphasize the role of $\eta$,
we will write $L_0$ as $L_0^\eta$.

Let $a_i$, $\psi_i$, $e_i$ be defined as in Section \ref{sec_preliminaries}. Let $\sigma\in H^2(\Sigma;\bC)$ be the Poincar\'e dual of $[\pt]$. We use $J_g$ to denote the Jacobian of $\Sigma_g$ and let $\mathscr{L}\to J_g\times \Sigma_g$ be a universal family of line bundles of degree $0$ over $\Sigma_g$ so that $\det \widetilde E$ is the pullback of $\mathscr{L}$ from $J_g$.
let $d_j=H^1(J_g;\bC)$ be defined by $d_j:=c_1(\mathscr{L})/[a_j]$. Then we have
$$
c_1(\widetilde{\bE})=\sum_{j=1}^{2g} d_j\otimes e_j + x\otimes 1
$$
for some $x\in H^2(J_g;\bC)$, where the right-hand side is viewed as the pullback from $H^\ast(J_g\times \Sigma_g)$ to $H^\ast(\widetilde{N}_{g,n}\times \Sigma)$.

We define a bundle $\widetilde {\mathbb{E}}_\eta$ by the exact sequence
$$
0\to \widetilde{\mathbb{E}}_\eta \to \widetilde{\mathbb{E}}\to  i_\ast(\widetilde{\mathbb{E}}|_{\widetilde{N}_{g,n}\times \eta}/\widetilde{\mathbb{L}})\to 0,
$$
where $i:\widetilde{N}_{g,n}\times \eta \hookrightarrow \widetilde{N}_{g,n}\times \Sigma_g$ is the inclusion map, and $i_*$ denotes the push-forward of sheaves.
Let $\pi_2:\widetilde{N}_{g,n}\times \Sigma_g \to \Sigma_g$ be the projection.
Then we have (see also \cite[Section 3.2]{XZ:excision})
\begin{align}
\parr\sHom(\pi_2^\ast L_0, \widetilde{\mathbb{E}} )&=\sHom(\pi_2^\ast L_0, \widetilde{\mathbb{E}}_\eta ), \nonumber\\
c_1(\sHom(\pi_2^\ast L_0,\widetilde{\bE}_\eta)) &= 
-(m+1)\otimes \sigma + \sum_{j=1}^{2g}d_j \otimes e_j+ x\otimes 1 \label{eq_c1_Hom_L,E}, \\
-\frac14 p_1(\ad \widetilde \bE) &= \alpha \otimes \sigma + \sum_{j=1}^{2g} \psi_j\otimes e_j+\beta\otimes 1, \nonumber\\
-\frac14 p_1(\ad \widetilde \bE_\eta)&=-\frac14 \tau_\eta^\ast p_1(\ad \widetilde \bE)=\alpha_\eta \otimes \sigma + \sum_{j=1}^{2g} \psi_j\otimes e_j+\beta\otimes 1 \label{eq_p1_E_eta},
\end{align}
where $\tau_\eta$ is the flip symmetry and $\alpha_\eta=\tau_\eta ^\ast \alpha$. We also have
\begin{equation}\label{eq_c2_Hom_L_E)}
c_2(\sHom(\pi_2^\ast L_0,\widetilde{\bE}_\eta))  = \frac14 c_1(\sHom(\pi_2^\ast L_0,\widetilde{\bE}_\eta))^2 -\frac14 p_1(\ad \widetilde \bE_\eta) .
\end{equation}

The stability condition implies 
$$
\parr\Hom( L_0, E)=0 
$$
for every stable parabolic bundle $E\in \widetilde{N}_{g,n}$. By the Riemann--Roch theorem, we have
$$
\chi(\parr\sHom(\pi_2^\ast L_0, E ))= 2g-2+2d+|\eta|=2g+m-1.
$$ 
Therefore by cohomology and base change we know that the derived pushforward 
${R}^0\pi_\ast\parr\sHom(\pi_2^\ast L_0, \widetilde{\mathbb{E}} )=0$ and 
${R}^1\pi_\ast\parr\sHom(\pi_2^\ast L_0, \widetilde{\mathbb{E}} )$ is  a vector bundle of rank
$2g+m-1$.
Hence $c_{i}({R}^1\pi_\ast\parr\sHom(\pi_2^\ast L_0, \widetilde{\mathbb{E}} ) )=0$ when $i\ge 2g+m$. On the other hand,
by the Grothendieck--Riemann--Roch theorem, these Chern classes can be written as a polynomial
in terms of the generators of $H^\ast(\widetilde{N}_{g,n})=H^\ast({N}_{g,n})\otimes H^\ast(J_g)$. 
According to \cite[Proposition 3.9]{XZ:excision} and its proof, when $\eta=\emptyset$ we have
$$
c_{g+k}({R}^1\pi_\ast\parr\sHom(\pi_2^\ast L_0, \widetilde{\mathbb{E}} ) )/[J_g]=\frac{1}{2^g}\xi_{k,n}
$$
where  $\xi_{k,n}\in \mathbb{C}[\alpha,\beta,\gamma]$
satisfies the following recursive relation:
\begin{equation}\label{eq_relation_xi_kn_0}
(k+1)\xi_{k+1,n}=\alpha \xi_{k,n}+(m-k)\beta \xi_{k-1,n} -\frac{\gamma}{2}\xi_{k-2,n}
\end{equation}
with
\begin{align*}
\xi_{0,n}&=1 \\
\xi_{1,n}&=\alpha \\
\xi_{2,n}&=\frac{\alpha^2}{2}+\frac{(m-1)\beta}{2}
\end{align*}
where $n=2m+1$. 

Since \eqref{eq_c1_Hom_L,E} does not depend on $\eta$ and \eqref{eq_c2_Hom_L_E)} can be obtained by substituting $\alpha$ in the formula of  $c_2(\sHom(\pi_2^\ast L_0,\widetilde{\bE}))$
by $\alpha_\eta$,
 the calculation of the Grothendieck--Riemann--Roch formula in  \cite[Section 3]{XZ:excision} also implies  
$$
c_{g+k}({R}^1\pi_\ast\parr\sHom(\pi_2^\ast L_0, \widetilde{\mathbb{E}} ) )/[J_g]=\frac{1}{2^g}\tau_\eta^\ast \xi_{k,n}
=\xi_{k,n}(\alpha_\eta,\beta,\gamma)
$$
for a general $\eta$.

We define 
$$
\widetilde{w}_0^\eta:=2^g c_{2g+m} ({R}^1\pi_\ast\parr\sHom(\pi_2^\ast L_0, \widetilde{\mathbb{E}} ) )
$$
 and 
$$
w_0^\eta:=\widetilde{w}_0^\eta/[J_g]=\tau_\eta^\ast\xi_{g+m,n}=\xi_{g+m,n}(\alpha_\eta,\beta,\gamma).
$$

Let $L_1$ be the dual of $L_0$. Then the weight of $L_1$ is
 $\frac14$ at $p\in \eta'$ and $-\frac14$ at $p\in \eta$. we have
 $$
\parr \deg L_1=-\parr \deg L_0=-\frac14.
 $$
 Similar as before, we have
 $$
\parr\sHom(\pi_2^\ast L_1, \widetilde{\mathbb{E}} )=\sHom(\pi_2^\ast L_1, \widetilde{\mathbb{E}}_{\eta'} ).
 $$
By the Riemann--Roch theorem, we have
$$
\chi(\parr\sHom(\pi_2^\ast L_1, E ))= 2g-2+2d+|\eta|=2g+m-2.
$$ 
for any $E\in \widetilde{N}_{g,n}$. We have
\begin{equation}
	\label{eqn_c1_L1_E}
c_1(\sHom(\pi_2^\ast L_1,\widetilde{\bE}_{\eta'})) = 
-m\otimes \sigma + \sum_{j=1}^{2g}d_j \otimes e_j+ x\otimes 1.
\end{equation}
Note that the right-hand side of \eqref{eqn_c1_L1_E} only differs from the right-hand side of
 \eqref{eq_c1_Hom_L,E} by substituting $m$ with $m-1$. Therefore,
 the calculation in  \cite[Section 3]{XZ:excision} implies that
\begin{equation}
	\label{eqn_chern_class_from_dual_bundle}
c_{g+k}(-\mathbf{R}\pi_\ast\parr\sHom(\pi_2^\ast L_1, \widetilde{\mathbb{E}} ) )/[J_g]=\frac{1}{2^g}\tau_{\eta'}^\ast \xi_{k,n-2}
=\xi_{k,n-2}(\alpha_{\eta'},\beta,\gamma),
\end{equation}
where $\xi_{k,n}$ is also well-defined when $n<0$.
We define 
$$
\widetilde{w}_1^\eta:=2^g c_{2g+m-1} (-\mathbf{R}\pi_\ast\parr\sHom(\pi_2^\ast L_1, \widetilde{\mathbb{E}} ) )
$$
 and 
$$
w_1^\eta:=\widetilde{w}_1^\eta/[J_g]=\tau_{\eta'}^\ast\xi_{g+m-1,n-2}=\xi_{g+m-1,n-2}(\alpha_{\eta'},\beta,\gamma).
$$

\subsection{The instanton moduli spaces}
\label{subsec_M0M1_subleading}

In this section, we study the moduli space of stable parabolic bundles over $\mathbb{CP}^1\times \Sigma_{g,n}$.
To define the stable parabolic bundle, we also need to fix a K\"ahler metric over $\mathbb{CP}^1\times \Sigma_{g}$.
Let $\omega_B$ and $\omega_C$ be K\"ahler forms of volume 1 over $B:=\mathbb{CP}^1$ and $C:=\Sigma_g$ respectively. 
We take $\omega_t=t\omega_B + \omega_C$ ($t>0$) as the K\"ahler metric used to define stable parabolic bundles over 
$B\times C$. When $t>2$, we say this metric lies in the $B$--chamber. When $0<t<2$ is small, we say
this metric is in the $C$--chamber.

We use $\widetilde \cM_0$ to denote the moduli space of stable parabolic bundles over 
$\mathbb{CP}^1\times \Sigma_{g,n}$
of energy $0$ \emph{without} fixing the determinant line bundle. We use 
$\cM_0$ to denote the moduli space of stable parabolic bundles over 
$\mathbb{CP}^1\times \Sigma_{g,n}$
of energy $0$ \emph{with} a fixed determinant line bundle.
Notice that by the energy condition, $\widetilde \cM_0$ and $ \cM_0$ are the same as  $\widetilde N_{g,n}$ and
$N_{g,n}$ respectively. Similarly, when the metric is $\omega_t=t\omega_B+\omega_C$,
we use 
we use $\widetilde \cM_1^t$ to denote the moduli space of stable parabolic bundles $E$ over 
$\mathbb{CP}^1\times \Sigma_{g,n}$
of energy $\frac14$ and $c_1(E)=P.D.[\{\pt\}\times \Sigma_g]$ \emph{without} fixing the determinant line bundle and we use 
$\cM_1^t$ to denote the corresponding moduli space of stable parabolic bundles with a fixed determinant.
To simplify the notation, we use $E\boxtimes F$ to denote the bundle $\pi_1^\ast E\otimes \pi_2^\ast F$ where $E$ and $F$
are bundles over $\mathbb{CP}^1$ and $\Sigma_g$ respectively, and $\pi_1,\pi_2$ are projections from 
$\mathbb{CP}^1\times \Sigma_{g}$ to $\mathbb{CP}^1$ and $\Sigma_g$ respectively.
\begin{Proposition}\label{prop_extension_tilde_M1}
When $t>2$ and $t$ is close to $2$, the moduli space $\widetilde \cM_1^t$ consists of parabolic bundle $E$ which comes from 
 a non-split extension of parabolic bundles
\begin{equation}\label{eq_extension_tilde_M1}
0\to \mathcal{O}(1)\boxtimes (L_1^\eta\otimes \zeta^\ast\otimes  \xi) \to E \to L_0^\eta\otimes \zeta \to 0
\end{equation}
where $\xi,\zeta\in J_g$. Therefore we have a fiber bundle structure
\begin{align*}
\coprod_{\eta~\text{even}}\mathbb{P}( \Ext^1_{\mathbb{CP}^1\times \Sigma_{g,n}}(L_0^\eta\otimes \zeta, \mathcal{O}(1)\boxtimes (L_1^\eta\otimes \zeta^\ast \otimes \xi))) 
\to 
\widetilde \cM_1 &\to J_g \times J_g \\
E &\mapsto  (\xi, \zeta)
\end{align*}
Moreover, the space $\widetilde \cM_1^t$ ($\cM_1^t$) is regular when $t$ is close to $2$.
\end{Proposition}
Before proving the above proposition, we prove some lemmas first.

\begin{Lemma}\label{lemma_stability_extension}
The parabolic bundle $E$ in the non-split extension \eqref{eq_extension_tilde_M1} is stable whenever $t>2$.
\end{Lemma}
\begin{proof}
Notice that
\begin{align*}
 \Ext^1_{B\times C}(L_0^\eta\otimes \zeta, \mathcal{O}(1)\boxtimes (L_1^\eta\otimes \zeta^\ast \otimes \xi))
 &=H^1(B\times C; \mathcal{O}(1)\boxtimes ( (L_1^\eta)^{\otimes 2}\otimes (\zeta^\ast)^{\otimes 2} \otimes \xi )\\
 &=H^0(B;\mathcal{O}(1))\otimes H^1(C;(L_1^\eta)^{\otimes 2}\otimes (\zeta^\ast)^{\otimes 2} \otimes \xi )
\end{align*}
and similarly 
$$
\Ext^1_{C}(L_0^\eta\otimes \zeta, (L_1^\eta\otimes \zeta^\ast \otimes \xi))=
H^1(C;(L_1^\eta)^{\otimes 2}\otimes (\zeta^\ast)^{\otimes 2} \otimes \xi).
$$
Using the above isomorphisms, the restriction map
$$
\Ext^1_{B\times C}(L_0^\eta\otimes \zeta, \mathcal{O}(1)\boxtimes (L_1^\eta\otimes \zeta^\ast \otimes \xi))
\to \Ext^1_{\{[(z_0,z_1)]\}\times C}(L_0^\eta\otimes \zeta, (L_1^\eta\otimes \zeta^\ast \otimes \xi))
$$
can be described as
\begin{align*}
H^0(B;\mathcal{O}(1))\otimes H^1(C;(L_1^\eta)^{\otimes 2}\otimes (\zeta^\ast)^{\otimes 2} \otimes \xi ) &\to
H^1(C;(L_1^\eta)^{\otimes 2}\otimes (\zeta^\ast)^{\otimes 2} \otimes \xi) \\
 r\otimes s &\mapsto r(z_0,z_1)s  
\end{align*}
where we view $r\in H^0(B;\mathcal{O}(1))$ as a linear homogeneous polynomial in two variables. From this description 
we see that a non-zero extension class is still non-zero after restricting to a generic slice $\{\pt\}\times C$. 

Suppose $E$ is the parabolic bundle in the non-split extension \eqref{eq_extension_tilde_M1}. 
We have
$$
\parr\deg_{B\times C} \mathcal{O}(1)\boxtimes (L_1^\eta\otimes \zeta^\ast\otimes  \xi)=
t\parr\deg_C L_1^\eta + 1=-\frac14 t+ 1
$$
and
$$
\parr\deg_{B\times C}  L_0^\eta\otimes \zeta=
t\parr\deg_C L_0^\eta=\frac14 t, \parr\deg E=1.
$$

Suppose there is a non-zero sheaf homomorphism $f:L\to E$ where $L$ is a line bundle over $B\times C$. Then $L$
has the form $\mathcal{O}(k)\boxtimes M$ where $M$ is a holomorphic line bundle over $C$. 
Now we equip $L$ with a parabolic structure so that $f$ is a homomorphism of parabolic bundles. 
If $f$ factors through $ \mathcal{O}(1)\boxtimes (L_1^\eta\otimes \zeta^\ast\otimes  \xi)$, then the parabolic 
degree of $L$ must be no greater than the parabolic degree of  $ \mathcal{O}(1)\boxtimes (L_1^\eta\otimes \zeta^\ast\otimes  \xi)$, hence the parabolic slope of $L$ is less than the parabolic slope of $E$. 
If $f$ does not factor through $ \mathcal{O}(1)\boxtimes (L_1^\eta\otimes \zeta^\ast\otimes  \xi)$, then we obtain a 
non-zero parabolic homomorphism $\bar{f}:L\to L_0^\eta\otimes \zeta$. Since the restriction to a generic slice $B\times \{\pt\}$
gives us a non-zero map $\mathcal{O}(k)\to \mathcal{O}$,  we have $k\le 0$. 
Since the restriction to a generic slice $\{\pt\}\times C$ gives us a non-zero map $M \to L_0^\eta\otimes \zeta$, we have
$$
\parr \deg_C M \le \parr\deg L_0^\eta\otimes\zeta=\frac14
$$
and the equality holds if and only if
the map $M\to L_0^\eta\otimes \zeta$ is an isomorphism of parabolic bundles. 
If the equality holds, then the extension \eqref{eq_extension_tilde_M1} splits
after restricting to a generic slice $\{\pt\}\times C$, which contradicts the previous conclusion.
 If the equality does not hold,
then we have $\parr \deg_C M \le -\frac14$ since the parabolic degrees of line bundles over $C$ can only differ by half integers. In this case we have
$$
\parr\deg_{B\times C} L=t\parr\deg M + k=-\frac14 t +k < \parr\text{slope}(E)=\frac12.
$$
In conclusion,  $E$ is a stable parabolic bundle. 
\end{proof} 

\begin{Lemma}\label{lemma_Mt_t_le_2}
The space $\widetilde \cM_1^t$ ($\cM_1^t$) is empty when $0<t\le 2$.
\end{Lemma}

\begin{proof}

A parabolic line bundle over $B\times C$ must have the form $\mathcal{O}(k)\boxtimes M$ where $M$ is a parabolic line bundle
over $C$. We have
$$
\parr\deg_{B\boxtimes C} \mathcal{O}(k)\otimes M=t\parr\deg_C M +k.
$$
We use $\cM_1^t$ to denote the moduli space of parabolic bundles with fixed determinant under the metric 
$\omega_t=t\omega_B+\omega_C$. Suppose $E\in \cM_1^t$ for some $0<t< 2$. Then given any non-zero parabolic homomorphism
$f:\mathcal{O}(k)\otimes M\to E$, we have
\begin{equation}\label{eq_stability_omegat}
\parr\deg_{B\times C} \mathcal{O}(k)\boxtimes M= t\parr\deg_C M +k < \parr\text{slope} E=\frac12.
\end{equation}
The restriction of $f$ to a generic slice $B\times \{\pt\}$ gives us a non-zero map 
$\mathcal{O}(k)\to E|_B$. According to the Birkhoff–Grothendieck theorem, 
$E|_B$ has the form $\mathcal{O}(a_1)\oplus \mathcal{O}(a_2)$. Therefore we have $k\le \max\{a_1,a_2\}$.
The restriction of $f$ to a generic slice $\{\pt\}\times C$ gives us a non-zero map $M\to E|_C$. Now we forget
the parabolic structure of $M, E|_C$ and view them as an ordinary holomorphic vector  bundles temporarily.  We have
$$
0\neq \Hom_C(M, E|_C)\cong H^0(C; M^{-1}\otimes E|_C^{-1})\cong H^1(C; M\otimes E|_C^{-1}\otimes K_C)^\ast
$$
Fix a degree 1 line bundle $L$ over $C$ and write $M$ as $L^k\otimes \ell$ where $\ell\in J_g$.
We want to show that
the the (ordianry) degree of $M$ is bounded above. Suppose not, then $k$ can be arbitrarily large. We have
$$
0\neq H^1(C; L^k\otimes \ell \otimes E|_C^{-1}\otimes K_C).
$$
On the other hand, for each fixed $\ell\in J_g$ we could find a large enough $k$ so that the above  cohomology
vanishes by the ampleness of $L$ and the Serre's vanishing theorem. Moreover, since the vanishing of cohomology is an open 
condition and $J_g$ is compact, we could find a large enough $k$ so that
$H^1(C; L^k\otimes \ell \otimes E|_C^{-1}\otimes K_C)=0$ for all $\ell\in J_g$, which is a contradiction.
Since
$$
\parr \deg_C M=\deg M+ \sum \text{weights}\le \deg M+ \frac14 n,
$$
we have $\parr\deg_C M$ is bounded above. 

Now we want to prove that there is a neighborhood of $(t-h,t+h)$ of $t$ so that  
$E$ is still stable under the metric $\omega_{t'}$ for all $t'\in (t-h,t+h)$.
If not then we could find a sequence $(M_n,k_n,t_n)$ such that 
\begin{itemize}
\item $M_n$ is a parabolic line bundle over $C$, $k_n\in\mathbb{Z}$, $t_n\in \mathbb{R}^+$;
\item $\lim_{n\to \infty} t_n= t$;
\item there exists a non-zero parabolic map $\mathcal{O}(k_n)\boxtimes M_n\to E$
\end{itemize} 
and 
\begin{equation}\label{eq_destabilize_omegatn}
\parr\deg_{B\times C}^{t_n} \mathcal{O}(k_n)\boxtimes M_n = t_n\parr\deg_C M_n + k_n \ge 
\parr\deg_{B\times C}^{t_n} E/\rank E=\frac12,
\end{equation}
where $\parr\deg_{B\times C}^{t_n}$ denotes the parabolic degree under the metric $\omega_{t_n}$.
We have proved that both $k_n$ and $\parr\deg M_n$ are bounded above. Inequality \eqref{eq_destabilize_omegatn}
implies that both  $k_n$ and $\parr\deg_C M_n$ are bounded below. Thus we conclude that $k_n$ and $\parr\deg_C M_n$
are bounded. Since $\parr\deg_C M_n$ can only be a multiple of $\frac14$, there are only finitely many
possibilities of  $k_n$ and $\parr\deg_C M_n$. Then \eqref{eq_destabilize_omegatn} contradicts \eqref{eq_stability_omegat}
when $n\to \infty$.

Now choose an arbitrary $t\in (0,2)$ and choose an arbitrary $E\in \cM_1^t$. We have proved that
$$
S=\{s\in (0,2)|E\in \cM_1^t\}
$$
is an open subset of $(0,2)$. If $S$ is not closed, then we could find $s_n\to r$ where $s_n\in S$ and $r\in (0,2)-S$.
Since $E\notin \cM_1^r$, we could find a non-zero parabolic map $\mathcal{O}(k)\boxtimes M\to E$ such that
$$
\parr\deg_{B\times C}^r \mathcal{O}(k)\boxtimes M=r\parr\deg_C M+ k \ge \parr\deg_{B\times C}^{r} E/\rank E=\frac12.
$$
On the other hand, we have
$$
\parr\deg_{B\times C}^{t_n} \mathcal{O}(k)\boxtimes M=t_n\parr\deg_C M+ k <\parr\deg_{B\times C}^{t_n} E/\rank E=\frac12.
$$
Therefore we must have
$$
\parr\deg_{B\times C}^r \mathcal{O}(k)\boxtimes M=r\parr\deg_C M+ k = \parr\deg_{B\times C}^{t_n} E/\rank E=\frac12
$$
which implies $E$ is strictly semistable. On the other hand, $\cM_1^r$ can be identified with the moduli space $M_1^{\ast r}$
of irreducible anti-self-dual $SO(3)$ singular connections of energy $\frac14$. For energy reason no bubble would appear in this moduli space. And it is proved in \cite{kronheimer2022relations}*{Lemma 5.5} that whenever $r\neq 2$, there is no reducible 
anti-self-dual connections. Therefore $M_1^{\ast r}$ (hence $\cM_1^r$) is compact. A semistable parabolic bundle is the limit of parabolic stable bundles, but the space $\cM_1^{ r}$ of stable parabolic bundles  is already compact, which is a contradiction.
So we prove that $S=(0,2)$ hence the space $\cM_1^{t}$ ($t\in (0,2)$) does not depend on $t$. 

It is proved in \cite{kronheimer2022relations}*{Lemma 5.4}\footnote{Indeed Lemma 5.4 and Lemma 5.5 in \cite{kronheimer2022relations} are only stated for genus $0$ case. but their proofs do not depend on the genus.}
 that when $t$ is small enough, $\cM^t_1=M_1^{\ast t}$ is empty. Therefore
we conclude that $\cM^t_1$ is empty for all $t\in (0,2)$. This also implies $\cM_1^2=M_1^{\ast 2}$ is empty: otherwise a stable parabolic 
bundle in $\cM_1^2$ would also lie in $\cM_1^t$ for any $t$ closed to $2$ by the previous discussion. 
\end{proof}

\begin{proof}[Proof of Proposition \ref{prop_extension_tilde_M1}]
By Lemma \ref{lemma_stability_extension}, the extension construction gives us a subspace of $\widetilde{\cM}_1$
with $2^{n-1}$ connected components. We want to prove that this subspace comprises the entire moduli space $\widetilde{\cM}_1$.
Since both this subspace and $\widetilde \cM_1$ can be acted by $J_g$ by tensor product,
it suffices to prove that bundles with fixed determinant from this extension construction  comprise the entire moduli space $\cM_1$.

According to Lemma \ref{lemma_Mt_t_le_2}, under the metric $\omega_2$,
the space $M_1^2$ of anti-seilf-dual $SO(3)$ connections under consists of reducible anti-self-dual connections. In the genus 0 case, these reducible connections
have been classified in \cite{kronheimer2022relations}*{Lemma 5.6}. The general case is very similar. We could find 
a singular $S^1$ connection $A$ on a complex line bundle $K$ with limiting holonomy $-1$ around $B\times \{p_1,\cdots,p_n\}$ such that 
$$
\frac{1}{2\pi i} {F_A}=-\omega_B+\frac12 \omega_C
$$
This curvature form is an anti-self-dual harmonic 2-form under the metric $\omega_2=2\omega_B+\omega_C$.
Then the connection $A\oplus \theta$ on $K\oplus \mathbb{R}$ is a reducible anti-self-dual  singular $SO(3)$ connection.
Twisting $A$ with flat line bundles over $C$ and using flip symmetries would give us all the connections in $M_1^2$. So we could  conclude that $M_1^2$ is the disjoint union $2^{n-1}$ copies of $J_g$.

The connection $A$ can be identified with a parabolic line bundle $L$ over $B\times C$ with weights $\pm \frac12$.  
We also have an identification of the homology of the deformation complex:
$$
H^0_A\cong H^0(B\times C;L), H^2_A\cong H^1(B\times C;L), H^2_A=H^2_+(B\times C; A)\cong H^2(B\times C;L)
$$

The bundle $L$ has the form $\mathcal{O}(k)\boxtimes M$ as before. The curvature of $A$ implies that $k=-1$.
By K\"unneth formula, We have
$$
H^2_A\cong H^2(B\times C;\mathcal{O}(-1)\boxtimes M )=\bigoplus_{i+j=2}
H^i(B;\mathcal{O}(-1))\otimes H^j(C;M)=0
$$
Hence the obstruction in the deformation complex of $A\oplus \theta$ is 
$$
H^2_{A\oplus \theta}=H^2_{A}\oplus H^2_{\theta}=H^2_+(B\times C;\theta)=H^2_+(B\times C;\mathbb{R}).
$$
The obstruction $H^2_+(B\times C;\mathbb{R})$ is a 1-dimensional space generated by the K\"ahler form $\omega_2$.
Now we consider the parameterized moduli space
$$
\mathscr{M}_1:=\bigcup_{t>0} M_1^t
$$
over the 1-parameter families of metrics $\{\omega_t\}$. The obstruction space of $\mathscr{M}_1$ at $A$ can be identified 
with the cokernel of 
\begin{align*}
D:\mathbb{R} &\to H^2_+(B\times C;\mathbb{R})\\
s &\mapsto  P^+_2 \frac{d}{dh}\Bigg|_{h=0}P^+_{2+h s} \frac{F_A}{2\pi i}
\end{align*}
where $P^+_t:\Omega^2(B\times C)\to H^2_{t,+}(B\times C;\mathbb{R})$ denotes the projection to the anti-self-dual 
harmonic 2-forms of $B\times C$ under the metric $\omega_t$ (cf. \cite{Froyshov_inequality_h}*{Section 4.1}). 
it is straightforward to verify that the map $D$ is non-zero. Hence $\mathscr{M}_1$ is regular at 
any $A\oplus \theta \in M_1^2$.
Since regularity is an open condition, $\mathscr{M}_1$ is also regular near $M_1^2$. 
Given  any $A\oplus \theta \in M_1^2$, we have
$$
\dim H^0_{A\oplus \theta}=1, \dim H^2_{A\oplus \theta}=1, 
$$
$$
\dim H^1_{A\oplus \theta}=\dim H^1_A+\dim H^1_\theta =(4g+2n-2)+2g=6g+2n-2
$$
The description of moduli space near reducibles (cf. \cite{DK}*{Section 4.3.6})
shows that when $t$ is close to $2$, $M_2^t$ is a $\mathbb{CP}^{4g+2n-4}$--bundle over $2^{n-1}$ copies of $J_g$.

On the other hand,  an direct index calculation shows that the subspace of $M_1^t$ constructed in Lemma \ref{lemma_stability_extension} is also a $\mathbb{CP}^{4g+2n-4}$--bundle over $2^{n-1}$ copies of $J_g$. So this subspace must
be the entire moduli space. We finishe the proof of the first part of Proposition \ref{prop_extension_tilde_M1}.

We pick $t$ close to $2$ so that $\mathscr{M}_1$ is regular at all $A\in M_1^t$. We want to prove that $M_1^t$ is also
regular. The moduli space $\mathscr{M}'=\bigcup_{2<s<t+h} M_1^s$ can be defined (at least near $A$) as $F^{-1}(0)$
where 
$$
F:\mathcal{B}\times I\to \mathcal{C}
$$ 
is a smooth Fredholm map between Banach manifolds, $0\in \mathcal{C}$ is a specific point
and  $I=(2,t+h)$ parameterizes the metrics we use. According to the first part of Proposition \ref{prop_extension_tilde_M1}, we have
$$
 G: M_1^t\times I\cong \mathscr{M}' \subset \mathcal{B}\times I.
$$
If $M_1^t$ is not regular at $A$, then $TF_{A,t}$ is not surjective after restricted to 
$T_A\mathcal{B}\times \{0\}\subset T_{A,t}(\mathcal{B}\times I)$.

Since $F\circ G(\{A\}\times I)=\{0\}$, we have 
$$
TF \circ TG(T_{A,t}(\{A\}\times I  )))=\{0\}
$$
which implies that there exists $v\in T_A\mathcal{B}$ such that
$$
TF(v,u)=0
$$
where $u$ is a generator of $T_tI$. Therefore for any $w\in T_A\mathcal{B}$, we have
$$
TF(w,u)=TF(w,0)+TF(0,u)=TF(w-v,0)\in TF_{A,t}(T_A\mathcal{B}\times \{0\} ).
$$
Therefore $TF_{A,t}: T_{A,t}\mathcal{B}\times I \to T_0\mathcal{C}$ is not surjective, which contradicts the regularity of
$\mathscr{M}_1$ at $(A,t)$.
%
 %
\end{proof}

From now on we fix $t>2$ which is close to $2$ and simply write $\cM_1$ for $\cM_1^t$. 
By Proposition \ref{prop_extension_tilde_M1} we see that if we fix the bundle $\xi$  in the extension \eqref{eq_extension_tilde_M1} to be the trivial 
bundle, then 
we could obtain all the bundles in $\cM_1$ from such extensions. We have a fiber bundle
\begin{align}
\coprod_{\eta~\text{even}}\mathbb{P}( \Ext^1_{\mathbb{CP}^1\times \Sigma_{g,n}}(L_0^\eta\otimes \zeta, \mathcal{O}(1)\boxtimes (L_1^\eta\otimes \zeta^\ast ))) 
\to 
 \cM_1 &\to  J_g \label{fiber_bundle_of_M1}\\
E &\mapsto  \zeta. \nonumber
\end{align}
More precisely, $\cM_1$ has $2^{n-1}$ connected components and each component is the projectivization
of a vector bundle over $J_g$. The fiber of this vector bundle over $\zeta\in J_g$ is 
$$
\Ext^1_{\mathbb{CP}^1\times \Sigma_{g,n}}(L_0^\eta\otimes \zeta, \mathcal{O}(1)\boxtimes (L_1^\eta\otimes \zeta^\ast )). 
$$

We have a fiber bundle
$$
\mathcal{M}_1 \to \widetilde{\mathcal{M}}_1 \xrightarrow{\det} J_g
$$
and $2^{2g}$-to-1 covering maps
\begin{align*}
F: \mathcal{M}_1\times J_g  &\to \widetilde{\mathcal{M}}_1\\
(E, \ell, p) &\mapsto (E\otimes \ell ,p).
\end{align*}
We define
$$
\Gamma_g:=\{\ell \in J_g|  \ell^2=\mathcal{O}  \}.
$$
Then we have an action of $\Gamma_g$ on $\mathcal{M}_1\times J_g$ defined by
\begin{align*}
\Gamma_g\times \mathcal{M}_1\times J_g &\to  \mathcal{M}_1\times J_g\\
(\ell, E,\ell') &\mapsto (E\otimes \ell, \ell'\otimes \ell^\ast) 
\end{align*}
Then it is clear that $\Gamma_g$ is the deck transformation group of the covering map $F$.
\begin{Proposition}\label{prop_homology_product_bundle}
The group $\Gamma_g$ acts trivially on $H^\ast(\mathcal{M}_1\times J_g;\mathbb{C})$. Therefore we have 
$$
H^\ast(\widetilde{\mathcal{M}}_1;\mathbb{C})\cong H^\ast( \mathcal{M}_1;\mathbb{C})\otimes H^\ast(J_g;\mathbb{C}).
$$ 
\end{Proposition}
\begin{proof}
Given an element $\ell\in \Gamma_g$, its tensor product action on $J_g$ is homotopic to the identity map. Therefore
it induces the identity map on $H^\ast( J_g;\mathbb{C})$. 

Suppose $E\in \cM_1$ comes from the extension
\begin{equation*}
0\to \mathcal{O}(1)\boxtimes (L_1^\eta\otimes \zeta^\ast) \to E \to L_0^\eta\otimes \zeta \to 0
\end{equation*}
then tensoring with $\ell$ we obtain an extension
\begin{equation*}
0\to \mathcal{O}(1)\boxtimes (L_1^\eta\otimes (\zeta\otimes \ell)^\ast) \to E\otimes \ell \to L_0^\eta\otimes \zeta\otimes \ell \to 0
\end{equation*}
Notice that we have a canonical isomorphism
$$
\Ext^1_{\mathbb{CP}^1\times \Sigma_{g,n}}(L_0^\eta\otimes \zeta, \mathcal{O}(1)\boxtimes (L_1^\eta\otimes \zeta^\ast ))
\cong
\Ext^1_{\mathbb{CP}^1\times \Sigma_{g,n}}(L_0^\eta\otimes \zeta\otimes \ell, \mathcal{O}(1)\boxtimes 
(L_1^\eta\otimes (\zeta\otimes \ell)^\ast ))
$$
The fiber bundle $\eqref{fiber_bundle_of_M1}$ is obtained by taking the $\mathbb{P}(H)$ where $H$
is a vector bundle over $J_g$ with fiber  
$\Ext^1_{\mathbb{CP}^1\times \Sigma_{g,n}}(L_0^\eta\otimes \zeta, \mathcal{O}(1)\boxtimes (L_1^\eta\otimes \zeta^\ast ))$.
Therefore tensoring with $\ell$ gives a bundle automorphism of $H$ which cover the map on $J_g$ also 
defined by tensoring with $\ell$. Now by the Leray-Hirsch theorem it is clear that tensoring with $\ell$ 
induces the identity map on $H^\ast(\mathcal{M}_1)$. This finishes the proof of the first part of the theorem.
The second part follows from the standard fact that the rational homology of the base is the invariant part
of the rational homology of a finite covering space. The same result also holds for complex coefficients.
\end{proof}

We can use the $\mu$-maps in Donaldson theory to define cohomology classes $\mu(w_i^\eta)$ ($i=0,1$) in the configuration space 
$\mathscr{B}^\ast(\Sigma_{g,n})$
of orbifold connections. Now pick a point $b\in \mathbb{CP}^1$. Then 
we could pullback these cohomology classes to the space of orbifold connections 
$\mathscr{B}^\ast(\mathbb{CP}^1\times \Sigma_{g,n})$
using the restriction map 
$\mathscr{B}^\ast(\mathbb{CP}^1\times \Sigma_{g,n}) \to 
\mathscr{B}^\ast(\{b\}\times \Sigma_{g,n})$. The duals of these cohomology classes can be represented by specific loci
according to the discussion in \cite{kronheimer2022relations}*{Section 4.1}. We denote these loci by $\widetilde{W}^\eta_i$ ($i=0,1$).
The space of stable parabolic bundles can be identified with the space of anti-self-dual orbifold connections. Therefore
we could view $\widetilde{\cM}_i$ ($i=0,1$) as subspaces of   $\mathscr{B}^\ast(\mathbb{CP}^1\times \Sigma_{g,n})$. 

The following proposition is an analogue of \cite{kronheimer2022relations}*{Lemma 5.9}. 
After we have Proposition \ref{prop_extension_tilde_M1}, its proof can be carried over almost verbatim from the proof 
of \cite{kronheimer2022relations}*{Lemma 5.9}, so we omit its proof.
\begin{Proposition}
We have
\begin{itemize}
\item[(a)]
A parabolic bundle $E\in \widetilde \cM_0 \cap \widetilde W_1^\eta$ if and only if $E\in \widetilde \cM_0$
$\parr\Hom(L_1, E)\neq 0$. Moreover in this case we have
a non-split extension of parabolic bundles 
\begin{equation}\label{eq_M0_extension}
0\to L_1\to E\to \xi \otimes L_0\to 0
\end{equation}
where $\xi=\det E\in J_g$.
\item[(b)]
A parabolic bundle $E\in  \widetilde \cM_1 \cap \widetilde W_0^\eta$ if and only if there is a non-split extension of parabolic bundles
$$
0\to \mathcal{O}(1)\boxtimes (L_1\otimes \xi) \to E \to L_0\to 0
$$
where $\xi\in J_g$ such that the extension splits restricting to $\{b\}\times \Sigma_{g,n}$.
\end{itemize}
\qed
\end{Proposition}
Now we have 
\begin{align*}
\Ext^1_{\mathbb{CP}^1\times \Sigma_{g,n}}(\xi\otimes L_0, L_1)&=H^1(\mathbb{CP}^1\times \Sigma_{g,n}, \xi^\ast \otimes  L_1^{\otimes 2}) 
=H^1(\Sigma_{g,n}, L_1^{\otimes 2} \otimes  \xi^\ast ) 
\end{align*}
and
\begin{align*}
\Ext^1_{\mathbb{CP}^1\times \Sigma_{g,n}}(L_0, \mathcal{O}(1)\boxtimes (L_1\otimes \xi))&= 
H^1(\mathbb{CP}^1\times \Sigma_{g,n}, \mathcal{O}(1)\boxtimes (L_1^{\otimes 2} \otimes \xi))\\
&= H^0(\mathbb{CP}^1, \mathcal{O}(1))\otimes H^1(\Sigma_{g,n},L_1^{\otimes 2} \otimes \xi   ).
\end{align*}
We have two fiber bundles 
\begin{align*}
\mathbb{P}(H^1(\Sigma_{g,n},   L_1^{\otimes 2} \otimes \xi^\ast ) )\to \widetilde \cM_0 \cap \widetilde W_1^\eta
\xrightarrow{\det} J_g \\
\mathbb{P}(Z_b\otimes H^1(\Sigma_{g,n},   L_1^{\otimes 2}  \otimes \xi) )\to \widetilde \cM_1 \cap \widetilde W_0^\eta
\xrightarrow{\det} J_g 
\end{align*}
To be more precise, let $\mathscr{L}\to \Sigma_{g,n} \times J_g$ be the universal line bundle and 
$p:\Sigma_{g,n} \times J_g\to J_g$ be the projection map. Then we have
$$
\widetilde \cM_0 \cap \widetilde W_1^\eta\cong \mathbb{P} (R^1p_\ast (L_1\otimes \mathscr{L}^{-1}) ),~
\cM_1 \cap \widetilde W_0^\eta \cong  \mathbb{P} (R^1p_\ast (L_1\otimes \mathscr{L}) )
$$

We have 
$H^\ast(N_{g,n})$ is generated by $\alpha,\beta,\psi_i,\delta_j$ and $H^\ast(\widetilde N_{g,n})$
is generated by $\alpha,\beta,\psi_i,\delta_j,d_k$ where $d_k$ is pulled back from $H^1(J_g)$ by
the determinant map. We can view  $\widetilde{w}_i^\eta$ as a polynomial in 
$\mathbb{C}[\alpha,\beta,\delta_i,\psi_j,d_k]$ and 
${w}_i^\eta$ as a polynomial in 
$\mathbb{C}[\alpha,\beta,\gamma,\delta_i]$.

Consider the map
\begin{align*}
\widetilde \cM_0 &\to \widetilde \cM_0 \\
E&\mapsto E\otimes (\det E)^\ast
\end{align*}
This map is an involution and under this map the extension \eqref{eq_M0_extension} becomes the extension
$$
0\to L_1\otimes \xi'\to E'\to L_0 \to 0 
$$
where $E'=E\otimes (\det E)^\ast$ and $\xi'=\det E'=(\det E)^\ast$. Therefore $ \widetilde \cM_0 \cap \widetilde W_1^\eta$
can also be described as a fiber bundle
$$
\mathbb{P}(H^1(\Sigma_{g,n},   L_1^{\otimes 2} \otimes \xi ) )\to \widetilde \cM_0 \cap \widetilde W_1^\eta
\xrightarrow{\det} J_g
$$
Now it is clear that we have the following.
\begin{Proposition}\label{prop_tilde_w0=w1}
We have a diffeomorphism 
$$\widetilde \cM_0 \cap \widetilde W_1^\eta \cong  \widetilde \cM_1 \cap \widetilde W_0^\eta $$
and 
$$
\langle \mu(\widetilde{w}_1^\eta)\smile \mu(g), [\widetilde \cM_0 ] \rangle =
\langle \mu(\widetilde{w}_0^\eta) \smile \mu(g), [\widetilde \cM_1 ] \rangle 
$$
for every $g\in \mathbb{C}[\alpha,\beta,\delta_i,\psi_j,d_k]$.
\end{Proposition}

The folowing is an analogue of \cite{kronheimer2022relations}*{Proposition 5.3}.
\begin{Proposition}
	\label{prop_subleading}
We have
$$
\langle \mu({w}_1^\eta)\smile \mu(g), [ \cM_0 ] \rangle =
\langle \mu({w}_0^\eta) \smile \mu(g), [ \cM_1 ] \rangle 
$$
for any $g\in \mathbb{C}[\alpha,\beta,\delta_i,\psi_j]$.

\end{Proposition}\label{prop_subleading_term}
\begin{proof}
Recall that  $w_i^\eta$ is the coefficient of $\prod_i d_i\smile d_{i+g}$ in $\widetilde{w}_i^\eta$
using the decomposition 
$H^\ast(\widetilde{N}_{g,n})\cong H^\ast(N_{g,n})\otimes H^\ast(J_g)$. Now the proposition follows
from Proposition \ref{prop_tilde_w0=w1} and the second part of Proposition \ref{prop_homology_product_bundle}.
\end{proof}

\begin{Corollary}
    \label{cor_subleading_term_formula}
There exists $W^\eta=w^\eta_0+(-1)^{g}\hat{\epsilon} w_1^\eta +w_2^\eta + \cdots \in J_{g,n}$ where 
$\hat \epsilon$ is defined in Section \ref{subsec_sign_E}, 
$w_0^\eta$, $w_1^\eta$ are defined as before and $ w_i^\eta\in \mathbb{C}[\alpha,\beta, \gamma,\delta_1,\cdots,\delta_n,\epsilon]/(\epsilon^2-1)$ is a homogeneous polynomial of degree $2g+2m-2i$. 
\end{Corollary}
\begin{proof}
According to \cite{kronheimer2022relations}*{Corollary 3.16}, there is a unique relation 
$$
U=U(0)+\hat{\epsilon}U(1)+U(2)+\hat{\epsilon}U(3)+\cdots \in J_{g,n}
$$
where $\deg U(i)=\deg U(0)-2i=2g+2m-2$ and 
 whose leading term  $U(0)=w^\eta_0$. Suppose $g\in \mathbb{A}_{g,n}$ is a
homogeneous element of degree $4g+2m$. Then we have
 $$
0=\langle \Phi(U), \Phi(g) \rangle=D(Ug)
=D_1(U(0)g)+D_0(U(1)g) 
 $$ 
where $D_i$ denotes the singular Donaldson invariants defined using
 moduli space of  singular (projectively) anti-self-dual $U(2)$ connections over 
$\mathbb{CP}^1\times \Sigma_{g.n}$ of energy $\frac{i}{4}$.  The first Chern class of the $U(2)$
bundle is $0$ when $i$ is even and is $P.D.[\{\pt\}\times \Sigma_g]$ when $i$ is odd.
Therefore we obtain
\begin{equation}\label{eq_M0M1_U0U1}
\langle \mu(U(1))\mu(g), [M_0]\rangle = - \langle \mu(U(0))\mu(g)), [M_1]  \rangle
\end{equation}
The space $M_0$ can be identified with $N_{g,n}^0$ and $\mu(U(1))\mu(g)=\Psi(U(1))\Psi(g)$. Since
$\Psi: \mathbb{A}_{g,n}\to H^\ast(N_{g,n}^0)$ is an isomorphism when $\deg\le 2g+2m-2$, \eqref{eq_M0M1_U0U1} determines
$U(1)$ uniquely. Comparing \eqref{eq_M0M1_U0U1} with Proposition \ref{prop_subleading} and using the identifications
$\cM_0=M_0, \cM_1=M_1$, we obtain $U(1)=-w_1^\eta$. However, there is a subtlety in the identification $\cM_1=M_1$:
the complex orientation and the orientation used to define the Donaldson invariants may be different under this identification.
The orientation to define the Donaldson invariants depends on the first Chern class $c=c_1(E)$
of the $U(2)$ bundle. Suppose $X$ is a K\"ahler manifold, the difference between the complex orientation and 
the  orientation to define the Donaldson invariants is 
$$
(-1)^{({c^2+K_X\cdot c})/{2}}
$$
according to \cite{Donaldson-orientation}*{Proposition 3.25} and its proof.
If we take $X=\mathbb{CP}^1\times \Sigma_g$ and $c=P.D.[\{\pt\}\times \Sigma_g]$, then we obtain $(-1)^{g-1}$. 
Even though \cite{Donaldson-orientation}*{Proposition 3.25} is only stated for the regular Donaldson invariants,
the singular case can be reduced to the regular case using the excision property of indices of elliptic operators. In summary
we obtain $U(1)=(-1)^g w_1^\eta$.
\end{proof}

\section{Characterization of the ring structure}
\label{sec_uniqueness}

This section gives a complete algebraic characterization for the ideals $J_{g,n}^\pm$. The main result is stated in Theorems \ref{thm_uniqueness_algebraic_structure_+} below. The proof of the theorem also yields an algorithm that computes the generators of $J_{g,n}^\pm$, which is described in Remark \ref{rmk_algorithm}.

%

Recall that $\bC[\alpha,\beta,\gamma,\delta_1,\dots,\delta_n]$ is a graded algebra with 
$$
\deg \alpha = 2, \quad \deg\beta = 4, \quad \deg\gamma = 6, \quad \deg\delta_i =2,
$$
and the ideals $I_{g,n}^d$ are homogeneous ideals. Suppose $0\neq f\in \bC[\alpha,\beta,\gamma,\delta_1,\dots,\delta_n]$,  the \emph{leading-order term} of $f$ is the homogeneous polynomial $f_0$ such that $\deg f_0 = \deg f$ and $\deg (f-f_0)<\deg f$. Also recall that the polynomials $\xi_{k,n}$ are defined in \eqref{eqn_defn_xi_kn}, and this definition remains valid even if $n$ is negative. For $I\subset \{1,\dots,n\}$, we use $\tau_I$ to denote the flip symmetry with respect to $I$. We also use $m$ to denote $(n-1)/2$. The homomorphism $\pi^{g,n+2}_{g,n}$ was defined in \eqref{eqn_defn_pi}. 

The formula of sub-leading term given by Corollary \ref{cor_subleading_term_formula} allows us to determine the signs of the eigenvalues in Proposition \ref{prop_sequence_of_eigenvalues}.

\begin{Lemma}
	\label{lem_fix_signs_eigenvalue}
	Let $\lambda_1,\lambda_2,\dots$ be given by Proposition  \ref{prop_sequence_of_eigenvalues}. Then 
	$$\lambda_i = (-1)^{i-1}(2i-1).$$
\end{Lemma}

\begin{proof}
	Consider the homomorphism from $\bC[\alpha,\beta,\gamma,\delta_1,\dots,\delta_n]$ to $\bC[\alpha]$ that takes $\alpha$ to $\alpha$, takes $\beta$ to $2$, and takes $\gamma,\delta_1,\dots,\delta_n$ to $0$. Then when $\epsilon =1$ (which holds on $\bV_{g,n}^+$), the image of $W^\eta$ in Corollary  \ref{cor_subleading_term_formula} is a polynomial of degree $g+m$ with the form 
	$$
	\frac{1}{(g+m)!} \alpha^{g+m} + (-1)^{g+m} \frac{1}{(g+m-1)!}\alpha^{g+m-1} + \text{lower order terms}.
	$$
	Hence it has at most $g+m$ distinct roots. By Proposition \ref{prop_sequence_of_eigenvalues} and Lemma \ref{lem_eigenvalue_and_ideal}, these roots are given by $\lambda_1,\dots,\lambda_{g+m}$. By Vieta's formula, this implies
	$$
	\lambda_1+\dots+\lambda_{g+m} = (-1)^{g+m+1} (g+m).
	$$
	Hence the result follows.
\end{proof}

The following theorem gives a complete characterization of the ideals $J_{g,n}^+$.

\begin{Theorem}
	\label{thm_uniqueness_algebraic_structure_+}
	Suppose for each pair $(g,n)$ such that $g$ is a non-negative integer and $n$ is a positive odd integer, there is an ideal $\hat J_{g,n}\subset \bC[\alpha,\beta,\gamma,\delta_1,\dots,\delta_n]$. Then $\hat J_{g,n}=J_{g,n}^+$ for all $(g,n)$ if and only if all of the following conditions hold:
	\begin{enumerate}
		\item $\hat J_{0,1}=(1)$.
		\label{cond_J01}
		\item 
		\label{cond_leading_terms_ideal}
		For every $f\in \hat J_{g,n}$, the leading-order term of $f$ is contained in $I_{g,n}^0$; for every homogeneous element $f_0\in I_{g,n}^0$, there exists $f\in \hat J_{g,n}$ such that the leading-order term of $f$ is equal to $f_0$.
		\item $\hat J_{g,n} $ is invariant under even flip symmetries.
		\label{cond_flip}
		\item 
		\label{cond_Jgn_symmetric_delta}
		$\hat J_{g,n}$ is symmetric with respect to $\delta_1,\dots,\delta_n$.
		\item 
		\label{cond_inclusion_hatJ_in_g}
		$\gamma \hat J_{g,n} \subset \hat J_{g+1,n}\subset \hat J_{g,n} $.
		\item  
		\label{cond_gamma_delta^2+beta-2_in_J}
		$\gamma^{g+1}\in \hat J_{g,n}$, $\delta_i^2+\beta-2 \in \hat J_{g,n}$. 
		\item 
		\label{cond_n+2_to_n}
		$\pi^{g,n+2}_{g,n}(\hat J_{g,n+2})\subset \hat J_{g,n}$.
		\item 
\label{cond_lowest_deg_polynomials}
If $(g,n)\neq (0,1)$, there exists a polynomial of the form
\begin{equation}
	\label{eqn_def_f_gn_in_Condition_for_hatJ}
\hat f_{g,n} = \xi_{g+m,n}+(-1)^{g+m}\tau_{\{1,\dots,n\}} (\xi_{g+m-1,n-2}) + h_{g,n},
\end{equation}
 where $\deg h_{g,n}\le 2(g+m-2),$
such that the following holds. If $m$ is odd, then $\tau_I(\hat f_{g,n})\in \hat J_{g,n}$ for all $I\subset\{1,\dots,n\}$ with $|I|$ even. If $m$ is even, then $\tau_I(\hat f_{g,n})\in \hat J_{g,n}$ for all $I\subset\{1,\dots,n\}$ with $|I|$ odd. 
\item 
\label{cond_eigenvalues}
Define $\hat{\bV}_{g,n}= \bC[\alpha,\beta,\gamma,\delta_1,\dots,\delta_n]/\hat J_{g,n}$. We use $x\in \hat \bV_{g,n}$ to denote the linear operator on $\hat \bV_{g,n}$ defined by the multiplication by $x$. Then the set of simultaneous eigenvalues of 
$$
(\alpha, \gamma,\delta_1,\dots,\delta_n)
$$
in $\hat \bV_{g,n}$ in the generalized eigenspace of $\beta$ with eigenvalue $2$ is given by 
$$
(\lambda,0,\dots,0)
$$
for $\lambda = (-1)^{i+1} (2i-1)$ with $i= 1, \dots, g+m$.
	\end{enumerate}
\end{Theorem}

\begin{remark}
	The ideal $J_{g,n}^-$ is characterized by Theorem \ref{thm_uniqueness_algebraic_structure_+} and Proposition \ref{prop_J_gn_pm_symmetry}. One can also state an analogous result of Theorem \ref{thm_uniqueness_algebraic_structure_+} for $J_{g,n}^-$, in which Conditions \ref{cond_J01} to \ref{cond_n+2_to_n} remain the same, and the sign of the sub-leading term in \eqref{eqn_def_f_gn_in_Condition_for_hatJ} and the signs of the eigenvalues in Condition \ref{cond_eigenvalues} are reversed.
\end{remark}

\begin{remark}
If $\hat J_{g,n} = J_{g,n}^+$, then by the earlier results, we know that $\hat J_{g,n}$ satisfies all the conditions listed in Theorem \ref{thm_uniqueness_algebraic_structure_+}: Condition \ref{cond_J01} follows immediately from the fact that $\bV_{0,1}=0$, Condition \ref{cond_leading_terms_ideal} follows from Propositions \ref{prop_I_realized_by_J} and \ref{prop_J_realized_by_I}, Condition \ref{cond_flip} follows from Proposition \ref{prop_Jgn_flip}, Condition \ref{cond_Jgn_symmetric_delta} follows from the definitions of $J_{g,n}^+$, Conditions \ref{cond_inclusion_hatJ_in_g} to \ref{cond_n+2_to_n} are given by Section \ref{subsec_properties_Ign_Jgn} and \eqref{eqn_gamma^g+1_in_Jgn}, Condition \ref{cond_lowest_deg_polynomials} is given by Corollary \ref{cor_subleading_term_formula}, and Condition  \ref{cond_eigenvalues} is given by Lemma \ref{lem_fix_signs_eigenvalue}. Therefore, we only need to show that if $\hat J_{g,n}$ satisfies all the conditions in Theorem \ref{thm_uniqueness_algebraic_structure_+}, then $\hat J_{g,n} = J_{g,n}^+$.
\end{remark}

For the rest of this section, we will assume that $\hat J_{g,n}$ satisfies all the conditions in Theorem \ref{thm_uniqueness_algebraic_structure_+}. The goal is to prove that $\hat J_{g,n} = J_{g,n}^+$.

\subsection{The ring $R_{g,n}$}
Consider the quotient ring 
$$R_{g,n} := \bC[\alpha,\beta,\gamma,\delta_1,\dots,\delta_n]/(\gamma^{g+1},\delta_1^2+\beta-2,\dots,\delta_n^2+\beta-2).$$
Since both $\hat J_{g,n}$ and $J_{g,n}^+$ satisfy Condition \ref{cond_gamma_delta^2+beta-2_in_J} in Theorem  \ref{thm_uniqueness_algebraic_structure_+}, we only need to show that $\hat J_{g,n}$ and $J_{g,n}^+$ have the same quotient image in $R_{g,n}$. 

When there is no risk of confusion, we will use a polynomial $f\in \bC[\alpha,\beta,\gamma,\delta_1,\dots,\delta_n]$ to denote its quotient image in $R_{g,n}$. 

Since the ideal $(\gamma^{g+1},\delta_1^2+\beta-2,\dots,\delta_n^2+\beta-2)$ is invariant under flip symmetries, the flip symmetry group acts on $R_{g,n}$. Since \eqref{eqn_prod_delta_i_notin_I} is satisfied by $R_{g,n}$, it makes sense to refer to the parity of flip symmetries on $R_{g,n}$. 

The map $\pi^{g,n+2}_{g,n}$ induces a map from $R_{g,n}$ to $R_{g,n+2}$, and we will use the same notation to denote the induced map. 

Recall that $\omega = \alpha + (\delta_1+\dots+\delta_n)/2$ and $\bC[\alpha,\beta,\gamma,\delta_1,\dots,\delta_n]$ is identified with $\bC[\omega,\beta,\gamma,\delta_1,\dots,\delta_n]$.
Also recall that for $I\subset\{1,\dots,n\}$, we use $\delta^I$ to denote $\prod_{i\in I} \delta_i$, and use $I^c$ to denote the complement of $I$ in $\{1,\dots,n\}$.
For each $I\subset \{1,\dots,n\}$ with $|I|\le m$, define
$R_{g,n}^I\subset R_{g,n}$ to be the quotient image of 
$$\delta^I\cdot \bC[\omega,\beta,\gamma] \oplus \delta^{I^c}\cdot \bC [\omega,\beta,\gamma]\subset  \bC[\omega,\beta,\gamma,\delta_1,\dots,\delta_n].$$
Note that $R_{g,n}^I$ is a $\bC$--linear subspace of $R_{g,n}$ but not an ideal. We have 
\begin{equation}
	\label{eqn_isotypic_decom_R_gn}
	R_{g,n} = \bigoplus_{|I|\le m} R_{g,n}^I,
\end{equation}
which is the isotypic decomposition of $R_{g,n}$ with respect to the actions by even flip symmetries. Let 
$$
p_I : R_{g,n} \to R_{g,n}^I
$$
be the projection map.

We establish several technical lemmas about $R_{g,n}$. 
\begin{Lemma}
	\label{lem_R_gn_canonical_rep}
	Let $V$ be the $\bC$--linear subspace of $\bC[\omega,\beta,\gamma,\delta_1,\dots,\delta_n]$ spanned by all monomials of the form $\omega^j\beta^k\gamma^s\delta_1^{a_1}\cdots\delta_n^{a_n}$ with $s\le g$ and $a_i\in \{0,1\}$ for all $i$. Then
	\begin{enumerate}
		\item 	 The composition map
		\begin{equation}
			\label{eqn_R_gn_canonical_rep}
			V\hookrightarrow \bC[\omega,\beta,\gamma,\delta_1,\dots,\delta_n] = \bC[\alpha,\beta,\gamma,\delta_1,\dots,\delta_n]\twoheadrightarrow R_{g,n}
		\end{equation}
		is a linear isomorphism. 
		\item The map \eqref{eqn_R_gn_canonical_rep} is equivariant with respect to flip symmetries.
	\end{enumerate}
\end{Lemma}
\begin{proof}
	Part (1) follows from the observation that 
	$$R_{g,n} = \bC[\omega,\beta,\gamma,\delta_1,\dots,\delta_n]/(\gamma^{g+1},\delta_1^2+\beta-2,\dots,\delta_n^2+\beta-2) $$
	 is a free $\bC[\omega,\beta,\gamma]/(\gamma^{g+1})$--module generated by monomials of the form $\delta_1^{a_1}\cdots\delta_n^{a_n}$ with $a_i\in \{0,1\}$ for all $i$. Part (2) is an immediate consequence of the definitions. 		
\end{proof}

\begin{Definition}
	Let $V$ be the linear space in Lemma \ref{lem_R_gn_canonical_rep}.
For each $x\in R_{g,n}$, there exists a unique $f_x\in V$, such that the quotient image of $f_x$ in $R_{g,n}$ is $x$. We call $f_x$ the \emph{canonical representative} of $x$.
\end{Definition}

\begin{Lemma}
	\label{lem_canonical_repre_lower_deg}
	Suppose $f\in \bC[\omega,\beta,\gamma,\delta_1,\dots,\delta_n]$ and $f'$ is the canonical representative of the image of $f$ in $R_{g,n}$. Then $\deg f' \le \deg f$.
\end{Lemma}

\begin{proof}
We only need to verify that case where $f$ is a monomial with coefficient $1$. Write
$f = \omega^j\beta^k\gamma^s \delta_{1}^{2b_1+a_1}\cdots\delta_{n}^{2b_n+a_n},$
where $a_i\in \{0,1\}$ and $b_i$ are non-negative integers. If $s > g$, then $f'=0$. If $s\le g$, then $f' =  \omega^j\beta^k\gamma^s \delta_{1}^{a_1}\cdots\delta_{n}^{a_n}(2-\beta)^{b_1+\dots+b_n}$. In either case, we have $\deg f' \le \deg f$. 
\end{proof}

\begin{Lemma}
	\label{lem_pi_injective_on_invariant_part}
Assume $n\ge 3$, and $I$ is a subset of $\{1,\dots,n\}$ with $|I|\le (n-1)/2$. Then the restriction of $\pi^{g,n}_{g,n-2}:R_{g,n}\to R_{g,n-2}$ to $R_{g,n}^I$ is injective.
\end{Lemma}

\begin{proof}
Let $V$ be as in Lemma \ref{lem_R_gn_canonical_rep}, let $V^I$ be the subspace of $V$ spanned by all monomials of the form $\omega^j\beta^k\gamma^s\delta^I$ ($s\le g$) and $\omega^j\beta^k\gamma^s\delta^{I^c}$  ($s\le g$). It is straightforward to verify that $\pi^{g,n}_{g,n-2}$ restricts to an injection on (the quotient image of) $V^I$. Since the canonical representatives of elements in $R_{g,n}^I$ are contained in $V^I$, the result is proved. 
\end{proof}

The following group of lemmas study the properties of the polynomial $\hat f_{g,n}$ \eqref{eqn_def_f_gn_in_Condition_for_hatJ} under the isotypic decomposition by even flip symmetry actions.
Recall that $p_I: R_{g,n}\to R_{g,n}^I$ is the projection with respect to the decomposition  \eqref{eqn_isotypic_decom_R_gn}. 
We take the lexicographical order for monomials of the form $\omega^j\beta^k\gamma^s\delta_1^{a_1}\cdots\delta_n^{a_n}$ with respect to the ordered list $(\omega, \delta_1,\dots,\delta_n,\beta,\gamma)$. In other words, for two monomials, their comparison is determined by the exponents on the first variable in $(\omega, \delta_1,\dots,\delta_n,\beta,\gamma)$ such that the respective exponents on the two monomials are different. 

\begin{Definition}
For $0\neq x\in R_{g,n}$, let $f_x$ be the canonical representative of $x$. Then $f_x$ can be uniquely decomposed as a sum of monomials:
$$
f_x = \sum_{s\le g,\, a_i\in\{0,1\}} c_{j,k,s,a_1,\dots,a_n}\cdot \omega^j\beta^k\gamma^s\delta_1^{a_1}\cdots\delta_n^{a_n}
$$
Suppose $c_{j,k,s,a_1,\dots,a_n}\cdot \omega^j\beta^k\gamma^s\delta_1^{a_1}\cdots\delta_n^{a_n}$ is the highest non-zero monomial in this sum with respect to the lexicographical order given by $(\omega, \delta_1,\dots,\delta_n,\beta,\gamma)$. We say that  $c_{j,k,s,a_1,\dots,a_n}\cdot \omega^j\beta^k\gamma^s\delta_1^{a_1}\cdots\delta_n^{a_n}$ is the \emph{highest monomial} in $x$.
\end{Definition}

\begin{Lemma}
	\label{lem_highest_term_of_pI(xi)}
	Let $I$ be a subset of $\{1,\dots,n\}$ with $|I|\le m$. Then the highest monomial in $p_I(\xi_{g+m,n})$ is equal to
	\begin{equation}
		\label{eqn_proj_I_xi_leading_term}
	2^{-m-1}\cdot (-1)^k \frac{2^{(r+1)/2}}{k!}\omega^{k}\delta^I
	\end{equation}
	with $k = g+m-|I|$ and $r = n-2|I|$.
\end{Lemma}
\begin{proof}
	Note that $\xi_{g+m,n}$ has degree $2(g+m)$, and among all monomials in $R_{g,n}^I$ with degree no greater than $2(g+m)$, the highest ones are scalar multiples of $\omega^{g+m-|I|}\delta^I$.
	Therefore the result follows from Lemma \ref{lem_formula_rho_kn} and Lemma \ref{lem_rho_leading_term_nonzero}.
\end{proof}

\begin{Lemma}
	\label{lem_projection_f_nonzero_I}
Let $\hat f_{g,n}$ be a polynomial of the form \eqref{eqn_def_f_gn_in_Condition_for_hatJ}, let $I$ be a subset of $\{1,\dots,n\}$ with $|I|\le m$. Then the highest monomial in $p_I(\hat f_{g,n})$ is given by \eqref{eqn_proj_I_xi_leading_term}. In particular, we have  $p_I(\hat f_{g,n})\neq 0$.
\end{Lemma}
\begin{proof}
	By Lemma \ref{lem_canonical_repre_lower_deg} and the definition of $\hat f_{g,n}$, the canonical representative of $p_I(\hat f_{g,n}) - p_I(\xi_{g+m,n})$ has degree strictly less than $2(g+m)$, therefore its highest monomial is lower than $\omega^{g+m-|I|}\delta^I$. Hence the result follows from Lemma \ref{lem_highest_term_of_pI(xi)} 
\end{proof}

\begin{Lemma}
	\label{lem_projection_tau_J_f_nonzero_I}
	Let $\hat f_{g,n}$ be a polynomial of the form \eqref{eqn_def_f_gn_in_Condition_for_hatJ}, let $I,J$ be subsets of $\{1,\dots,n\}$ such that $|I|\le m$. Then  $p_I(\tau_J(\hat f_{g,n}))\neq 0$.
\end{Lemma}
\begin{proof}
	Note that $p_I(\tau_J(\hat f_{g,n}))=\tau_J(p_I(\hat f_{g,n}))$.
	By Lemma \ref{lem_projection_f_nonzero_I}, $\tau_J^2(p_I(\hat f_{g,n}))=p_I(\hat f_{g,n})\neq 0$. So $\tau_J(p_I(\hat f_{g,n}))\neq 0.$
\end{proof}

%

\subsection{Properties of $\hat J_{g,n}$}
\begin{Lemma}
	\label{lem_lower_order_terms_hatJ}
If $f\in \hat J_{g,n}$ satisfies $\deg f < 2(g+m)$, then the image of $f$ in $R_{g,n}$ is zero.
\end{Lemma}
\begin{proof}
	Assume there exists $0\neq f\in \hat J_{g,n}$ such that $\deg f < 2(g+m)$ and the image of $f$ in $R_{g,n}$ is non-zero. Take such an $f$ with the minimal possible degree. Let $f_0$ be the leading-order term of $f$. By Condition \ref{cond_leading_terms_ideal} in Theorem \ref{thm_uniqueness_algebraic_structure_+}, we have $f_0\in I_{g,n}^0$. Since $\deg f_0 = \deg f < 2(g+m)$, by Theorem \ref{thm_generators_of_Ign}, we know that
	$$
	f_0\in (\gamma^{g+1}, \delta_1^2+\beta,\dots,\delta_n^2+\beta),
	$$
	so there exists $f'\in (\gamma^{g+1},\delta_1^2+\beta-2,\dots,\delta_n^2+\beta-2)$ such that the leading order term of $f'$ is also $f_0$. By Condition \ref{cond_gamma_delta^2+beta-2_in_J}, we have $f'\in \hat{J}_{g,n}$, therefore $f-f'\in \hat J_{g,n}$. Since $f-f'\notin  (\gamma^{g+1},\delta_1^2+\beta-2,\dots,\delta_n^2+\beta-2) $, its image in $R_{g,n}$ is non-zero. Since $\deg(f-f')<\deg(f)$, this contradicts the assumption that $f$ minimizes the degree.
\end{proof}

Let $\hat X_{g,n}$ be the subset of $R_{g,n}$ consisting of the quotient images of all $f\in \hat J_{g,n}$ such that $\deg f\le g+m$. By definition, $\hat X_{g,n}$ is a $\bC$--linear space.

\begin{Lemma}
	\label{lem_basis_of_hatX}
	Let $\hat f_{g,n}$ be a polynomial of the form \eqref{eqn_def_f_gn_in_Condition_for_hatJ}
such that Condition \ref{cond_lowest_deg_polynomials} in Theorem \ref{thm_uniqueness_algebraic_structure_+} holds with respect to $\hat f_{g,n}$.  
\begin{enumerate}
	\item If $m$ is odd, then $\hat X_{g,n}$ is the linear space spanned by $\tau_I(\hat f_{g,n})$ for all $I\subset\{1,\dots,n\}$ with $|I|$ even.
	\item If $m$ is even, then $\hat X_{g,n}$ is the linear space spanned by $\tau_I(\hat f_{g,n})$ for all $I\subset\{1,\dots,n\}$ with $|I|$ odd. 
\end{enumerate}
Moreover, we have $\dim_\bC \hat X_{g,n} = 2^{n-1}$. 
\end{Lemma}

\begin{proof}
	We first prove the result when $m$ is odd.
	
	By the definition of $\hat f_{g,n}$, we have $\tau_I(\hat f_{g,n})\in \hat J_{g,n}$ whenever $|I|$ is even. 
	By Lemma \ref{lem_projection_f_nonzero_I}, the linear space spanned by all $\tau_I(\hat f_{g,n})$ with $|I|$ even has dimension at least $2^{n-1}$ in $R_{g,n}/(\beta-2,\gamma)$, therefore the space has dimension at least $2^{n-1}$ in $R_{g,n}$. Hence these elements are linearly independent in $R_{g,n}$. 
	
	We still need to show that $\{\tau_I(\hat f_{g,n})\}_{|I|\text{ even}}$ is a spanning set of $\hat X_{g,n}$. 
	Suppose $0\neq f\in \hat J_{g,n}$ and $\deg f= g+m$, and let $f_0$ be the leading-order term of $f$. By Condition \ref{cond_leading_terms_ideal} and Theorem \ref{thm_generators_of_Ign}, there exists $f'$ which is a linear combination of $\tau_I(\hat f_{g,n})\in \hat J_{g,n}$ for $|I|$ even, such that the leading-order term of $f'$ is also $f_0$. By Lemma \ref{lem_lower_order_terms_hatJ}, the images of $f'$ and $f$ in $R_{g,n}$ are the same. Hence the desired statement is proved.
	
	The proof for $m$ odd is almost verbatim, except that we need to replace the citation to Lemma \ref{lem_projection_f_nonzero_I} with Lemma \ref{lem_projection_tau_J_f_nonzero_I}, and replace all the words ``even'' with ``odd''.
\end{proof}

Note that if $g'\ge g$, then $R_{g,n}$ is a quotient ring of $R_{g',n}$, therefore it makes sense to refer to the quotient image of $\hat X_{g',n}$ in $R_{g,n}$. 
\begin{Lemma}
	\label{lem_Xgn_generates_hatJ}
	The quotient image of $\hat J_{g,n}$ in $R_{g,n}$ is an ideal generated by the quotient images of $\hat X_{g,n}, \hat X_{g+1,n}, \hat X_{g+2,n}$.
\end{Lemma}
\begin{proof}
	By Condition \ref{cond_inclusion_hatJ_in_g} in Theorem \ref{thm_uniqueness_algebraic_structure_+}, the quotient images of $\hat X_{g,n}, \hat X_{g+1,n}, \hat X_{g+2,n}$ in $R_{g,n}$ are all contained in the quotient image of $\hat J_{g,n}$. We only need to show that the quotient images of $\hat X_{g,n}, \hat X_{g+1,n}, \hat X_{g+2,n}$ in $R_{g,n}$ generate the image of $\hat J_{g,n}$ as an ideal. 
	
	Assume there exists $0\neq f\in \hat J_{g,n}$ such that its image in $R_{g,n}$ is not contained in the ideal generated by the quotient images of $\hat X_{g,n}, \hat X_{g+1,n}, \hat X_{g+2,n}$. Take such an $f$ with the minimal possible degree. Then by Condition \ref{cond_leading_terms_ideal}, Theorem \ref{thm_generators_of_Ign}, and Condition \ref{cond_lowest_deg_polynomials}, there exists $f'\in \hat J_{g,n}$ such that $f$ and $f'$ have the same leading-order term, and that the image of $f'$ in $R_{g,n}$ is contained in the ideal generated by the images of $\hat X_{g,n}, \hat X_{g+1,n}, \hat X_{g+2,n}$. Therefore, $f-f'\in \hat J_{g,n}$ and its image in $R_{g,n}$ is not contained in the ideal generated by the quotient images of $\hat X_{g,n}, \hat X_{g+1,n}, \hat X_{g+2,n}$. Since $\deg(f'-f)<\deg f$, this contradicts the definition of $f$. 
\end{proof}

\begin{remark}
	A priori, the polynomials $\hat f_{g,n}$ in Lemma \ref{lem_basis_of_hatX} may depend on the ideals $\hat J_{g,n}$. So Lemma \ref{lem_basis_of_hatX} and Lemma \ref{lem_Xgn_generates_hatJ} do not immediately imply the uniqueness of the ideal $\hat J_{g,n}$. 
\end{remark}

Note that when $g\ge 1$, if the quotient image of $f\in \bC[\alpha,\beta,\gamma,\delta_1,\dots,\delta_n]$ in $R_{g-1,n}$ is zero, then the quotient image of $\gamma \cdot f$ in $R_{g,n}$ is zero. Therefore, the multiplication by $\gamma$ is a well-defined map from $R_{g-1,n}$ to $R_{g,n}$.
\begin{Lemma}
	\label{lem_hatJ_generating_set_gamma}
		Assume $g\ge 1$. Then the quotient image of $\hat J_{g,n}$ in $R_{g,n}$ is an ideal generated by the images of $\hat X_{g,n}, \hat X_{g+1,n}, \gamma \cdot \hat X_{g-1,n}$ in $R_{g,n}$.
\end{Lemma}
\begin{proof}
	By Condition \ref{cond_inclusion_hatJ_in_g} in Theorem \ref{thm_uniqueness_algebraic_structure_+}, the images of $\hat X_{g,n}, \hat X_{g+1,n}, \gamma \cdot \hat X_{g-1,n}$ in $R_{g,n}$ are contained in the quotient image of $\hat J_{g,n}$. The proof of the lemma is then verbatim as the proof of Lemma \ref{lem_Xgn_generates_hatJ}, except that we need to replace the citation of Theorem \ref{thm_generators_of_Ign} with Corollary \ref{cor_generating_set_Ign}.
\end{proof}

\subsection{Proof of Theorem \ref{thm_uniqueness_algebraic_structure_+}}
Since $J_{g,n}^+$ also satisfies all the conditions listed in Theorem \ref{thm_uniqueness_algebraic_structure_+}, all the previous lemmas proved for $\hat J_{g,n}$ also apply to $J_{g,n}^+$. Let $X_{g,n}$ be the subset of $R_{g,n}$ consisting of the quotient images of all $f\in J_{g,n}^+$ such that $\deg f\le g+m$. By Lemma \ref{lem_Xgn_generates_hatJ}, we only need to show that $\hat X_{g,n} = X_{g,n}$ for all $(g,n)$. The proof is divided into 3 steps.\\

\textbf{Step 1.}
First, note that if $(g,n) = (0,3)$ or $(1,1)$, then $g+m=1$. In this case, we must have $h_{g,n}=0$ in \eqref{eqn_def_f_gn_in_Condition_for_hatJ}. So the polynomial $\hat f_{g,n}$ is uniquely determined by Condition \ref{cond_lowest_deg_polynomials} in Theorem \ref{thm_uniqueness_algebraic_structure_+}. By Lemma \ref{lem_basis_of_hatX}, this implies $X_{g,n} = \hat X_{g,n}$.\\

\textbf{Step 2.}
Now we show that $X_{g,1} = \hat X_{g,1}$ for all $g\ge 1$ using induction on $g$. The case for $g=0$ follows from Condition \ref{cond_J01}, and the case for $g=1$ is proved in Step 1.
We also directly verify the case for $g=2$. In this case, $h_{g,n}$ must be a constant in \eqref{eqn_def_f_gn_in_Condition_for_hatJ}. 
By Condition \ref{cond_eigenvalues} in Theorem \ref{thm_uniqueness_algebraic_structure_+} and Lemma \ref{lem_eigenvalue_and_ideal}, for each $i=1,\dots,g+m$, the homomorphism
$$
\bC[\alpha,\beta,\gamma,\delta_1,\dots,\delta_n] \to \bC
$$
sending $(\alpha,\beta,\gamma,\delta_1,\dots,\delta_n)$ to $((-1)^{i+1}(2i-1),2,0,\dots,0)$ maps $\hat J_{g,n}$ to $(0)$. 
The value of $h_{2,1}$ is then uniquely determined by this condition.

Now assume $g\ge 3$ and $X_{k,1} = \hat X_{k,1}$ for all $k<g$, we show that $X_{g,1} = \hat X_{g,1}$. By Lemma \ref{lem_basis_of_hatX}, $\hat X_{g,1}$ is a one-dimensional linear space generated by (the quotient image) of a polynomial of the form \eqref{eqn_def_f_gn_in_Condition_for_hatJ}, and the same statement holds for $X_{g,1}$. Let $\hat f_{g,1}$ be the polynomial given by Condition \ref{cond_lowest_deg_polynomials} with respect to $\hat J_{g,n}$, and let $f_{g,1}$ be the corresponding polynomial for $J_{g,1}^+$. Although the choices of $\hat f_{g,1}$ and $f_{g,1}$ may not be unique, Lemma \ref{lem_lower_order_terms_hatJ} shows that the images of $\hat f_{g,1}$ and $ f_{g,1}$ in $R_{g,1}$ are uniquely determined by the ideals.

By Condition \ref{cond_inclusion_hatJ_in_g}, we have
$f_{g,1} \in J_{g-2,n}^+,
\hat f_{g,1} \in \hat J_{g-2,n}.$
By the induction hypothesis, $X_{k,1} = \hat X_{k,1}$ for $k=g-1,g-2,g-3$. Therefore, by Lemma \ref{lem_hatJ_generating_set_gamma}, $J_{g-2,n}^+ = \hat J_{g-2,n} $. By the definition of $f_{g,1}$ and $\hat f_{g,1}$, we have $\deg (f_{g,1} -\hat f_{g,1})\le 2(g-2)$. Therefore, the image of $f_{g,1}-\hat f_{g,1}$ in $R_{g-2,1}$ lies in $X_{g-2,1} = \hat X_{g-2,1}$. By Lemma \ref{lem_basis_of_hatX}, there exists $c\in \bC$ such that 
\begin{equation}
	\label{eqn_linear_relation_fg1_f(g-2)1}
f_{g,1}-\hat f_{g,1} - cf_{g-2,1}\in (\delta_1^2+\beta-2, \gamma^{g-1}).
\end{equation}
By Condition \ref{cond_eigenvalues}, for each $i=1,\dots,g$, the homomorphism
$$
\bC[\alpha,\beta,\gamma,\delta_1,\dots,\delta_n] \to \bC
$$
sending $(\alpha,\beta,\gamma,\delta_1,\dots,\delta_n)$ to $((-1)^{i+1}(2i-1),2,0,\dots,0)$ maps both $f_{g,1}$ and $\hat f_{g,1}$ to $0$. Therefore by \eqref{eqn_linear_relation_fg1_f(g-2)1}, it maps $c f_{g-2,1}$ to $0$. By \eqref{eqn_def_f_gn_in_Condition_for_hatJ}, under the homomorphism
$$
\bC[\alpha,\beta,\gamma,\delta_1,\dots,\delta_n] \to \bC[\alpha]
$$
sending $\alpha$ to $\alpha$ and $(\beta,\gamma,\delta_1,\dots,\delta_n)$ to $(2,0,\dots,0)$, 
 the image of $f_{g-2,1}$ is a polynomial in $\alpha$ with leading term $\alpha^{g-2}/(g-2)!$. Therefore, it has at most $g-2$ distinct roots. As a result, we must have $c=0$.
 
 In summary, we have proved that $f_{g,1}-\hat f_{g,1} \in (\delta_1^2+\beta-2, \gamma^{g-1}).$ Since 
 $$\deg \gamma^{g-1} = 6(g-1) > \deg (f_{g,1}-\hat f_{g,1}),$$
Lemma \ref{lem_high_power_gamma_no_effect} below will show that $f_{g,1}-\hat f_{g,1} \in (\delta_1^2+\beta-2)$. By Lemma \ref{lem_basis_of_hatX}, this implies $X_{g,1} = \hat X_{g,1}$. 
 \\
 
 \textbf{Step 3.} Now we show that $X_{g,n} = \hat X_{g,n}$ for all $g,n$. The case when $n=1$ was proved in Step 2. We prove the general case by induction on $n$. Assume $X_{g,n} = \hat X_{g,n}$ \emph{for all} $g$, we show that $X_{g,n+2} = \hat X_{g,n+2}$ \emph{for all} $g$.  By the induction hypothesis and Lemma \ref{lem_Xgn_generates_hatJ}, we have $J_{g,n}^+=\hat J_{g,n}$ for all $g$. 
 
 Let $\hat f_{g,n}$ be a polynomial given by Condition \ref{cond_lowest_deg_polynomials} with respect to $\hat J_{g,n}$, and let $f_{g,n}$ be a corresponding polynomial for $J_{g,n}^+$. 
 
For $I\subset\{1,\dots,n+2\}$,
let $\hat X_{g,n+2}^I$
 be the intersection of $\hat X_{g,n+2}$ with $R_{g,n+2}^I$. Since $\hat X_{g,n+2}$ is invariant under even flip symmetries, we have
 $$
 \hat X_{g,n+2} = \bigoplus_{|I|\le (n+1)/2}\hat X_{g,n+2}^I,
 $$
 and this is the isotypic decomposition of $\hat X_{g,n+2}$.
 
 By Lemma \ref{lem_projection_f_nonzero_I}, $\dim \hat X_{g,n+2}^I\ge 1$ for all $I$ with $|I|\le (n+1)/2$. By Lemma \ref{lem_basis_of_hatX}, we have $\dim_\bC \hat X_{g,n+2} = 2^{n+1}$, so  $\dim_\bC \hat X_{g,n+2}^I =1 $ for each $I$ with $|I|\le (n+1)/2$. 
 
 Similarly, let $X_{g,n+2}^I = X_{g,n+2}\cap R_{g,n+2}^I$, then $\dim_\bC X_{g,n+2}^I =1$.
 
 By the definitions of $\hat f_{g,n+2}$ and $f_{g,n+2}$, we know that   
  $p_I(\hat f_{g,n+2})-p_I(f_{g,n+2})$ is represented by a polynomial with degree no greater than $2(g+m-1)$, where $m=(n-1)/2$.  By Condition \ref{cond_n+2_to_n} and the induction hypothesis, 
 \begin{equation}
 	\label{eqn_pi_image_discrepancy}
 \pi^{g,n+2}_{g,n}(p_I(\hat f_{g,n+2})-p_I(f_{g,n+2}))  \in R_{g,n}
 \end{equation}
 is contained in the quotient image of $J_{g,n}^+ = \hat J_{g,n}$. 
 Since \eqref{eqn_pi_image_discrepancy} is represented by a polynomial with degree at most $2(g+m-1)$, Lemma \ref{lem_lower_order_terms_hatJ} implies that it is zero. By Lemma \ref{lem_pi_injective_on_invariant_part}, we have $p_I(\hat f_{g,n+2})-p_I(f_{g,n+2}) = 0$. 
 Hence $X_{g,n}^I = \hat X_{g,n}^I$.
\qed

\vspace{\baselineskip}
We finish the arguments by establishing the following lemma, which was used in Step 2 of the proof of Theorem \ref{thm_uniqueness_algebraic_structure_+}.
\begin{Lemma}
	\label{lem_high_power_gamma_no_effect}
	Suppose $f\in\bC[\alpha, \beta,\gamma,\delta_1]$ satisfies 
	$f\in (\delta_1^2+\beta-2,\gamma^k)$ and $\deg f < 6k,$
	then $f\in (\delta_1^2+\beta-2).$
\end{Lemma}

\begin{proof}
	Let $V$ be the $\bC$--linear subspace of $\bC[\alpha, \beta,\gamma,\delta_1]$ spanned by the monomials of the form $\alpha^j\beta^k\gamma^s\delta^a$ with $a\in \{0,1\}$, then the same argument as Lemma \ref{lem_R_gn_canonical_rep} shows that the composition map
	$$
	V\hookrightarrow \bC[\alpha, \beta,\gamma,\delta_1] \twoheadrightarrow \bC[\alpha, \beta,\gamma,\delta_1]/(\delta_1^2+\beta-2)
	$$
	is a linear isomorphism. Hence for each $x\in  \bC[\alpha, \beta,\gamma,\delta_1]/(\delta_1^2+\beta-2)$, there exists a unique $f_x\in V$ such that $f_x$ represents $f$, and we call it the \emph{canonical representative} of $x$. 
	
	Let $f'$ be the canonical representative of the quotient image of $f$ in the ring $\bC[\alpha, \beta,\gamma,\delta_1]/(\delta_1^2+\beta-2)$. By the same argument as Lemma \ref{lem_canonical_repre_lower_deg}, we see that $\deg f' \le \deg f$. 
	
	Since $f\in  (\delta_1^2+\beta-2,\gamma^k)$, the image of $f$ in $\bC[\alpha, \beta,\gamma,\delta_1]/(\delta_1^2+\beta-2)$ is contained in the ideal generated by $\gamma^k$. Since the multiplication by $\gamma^k$ maps $V$ into $V$, the canonical representative $f'$ must be a multiple of $\gamma^k$. Since $\deg f'\le \deg f < 6k = \deg \gamma^k$, this implies $f'=0$, therefore $f\in (\delta_1^2+\beta-2).$
\end{proof}

\begin{remark}
	\label{rmk_algorithm}
	It is straightforward to turn the proof of Theorem \ref{thm_uniqueness_algebraic_structure_+} into an algorithm that computes the generating sets of $J_{g,n}^\pm$. In Section \ref{sec_n=1}, we will explicitly carry out the argument of Step 2 and write down a recursive formula for the generating set of $J_{g,n}^+$ when $n=1$. For $n\ge 3$, one can recursively compute $X_{g,n+2}$ from $X_{g,n}$ by the following argument.  For each $I\subset \{1,\dots,n+2\}$ with $|I|\le (n+1)/2$,  Conditions \ref{cond_n+2_to_n} and \ref{cond_lowest_deg_polynomials} impose a system of linear equations on the coefficients of the canonical representative of  $p_I(f_{g,n})$, where $f_{g,n}$ is a polynomial given by Condition \ref{cond_lowest_deg_polynomials}. Step 3 of the proof of Theorem  \ref{thm_uniqueness_algebraic_structure_+} shows that the linear system always has a unique solution.
\end{remark}

\section{The one-punctured case}
\label{sec_n=1}

This section finds a recursive formula for the generating set of $J_{g,1}^+$.
When $n=1$, we use $\tau$ to denote the unique non-trivial flip symmetry, and we abbreviate $\delta_1$ to $\delta$. 

Define $f_0=1$. When $g\ge 1$, let $f_g$ be a polynomial given by Condition \ref{cond_lowest_deg_polynomials} in Theorem \ref{thm_uniqueness_algebraic_structure_+} with respect to the ideal $J_{g,1}^+$, then $f_g$ has the form
\begin{equation}
	\label{eqn_fg1_form}
f_{g} = \xi_{g,1}+(-1)^{g}\tau(\xi_{g-1,-1}) + h_{g,1}
\end{equation}
with $\deg h_{g,1} \le 2(g-2)$. By Lemma \ref{lem_basis_of_hatX}, the linear space $X_{g,1}$ defined in Section \ref{sec_uniqueness} is the one-dimensional space spanned by the quotient image of $\tau f_{g}$ in $R_{g,1}$. 
By Lemma \ref{lem_Xgn_generates_hatJ}, we have 
\begin{equation}
	\label{eqn_generators_Jg1+}
J_{g,1}^+ = (\tau f_{g},\,\tau  f_{g+1},\, \tau f_{g+2},\, \gamma^{g+1},\, \delta^2+\beta-2).
\end{equation}

We start by computing $f_1$ and $f_2$. Since $\xi_{1,1}=\alpha$, $\xi_{0,-1}=1$, by \eqref{eqn_fg1_form}, we have
$$
f_1 = \alpha - 1.
$$
Similarly, since $\xi_{2,1}= (\alpha^2-\beta)/2$, $\xi_{1,-1}= \alpha$, we have 
$$f_2 = (\alpha^2-\beta)/2 + (\alpha+\delta) + \text{constant}.$$
 By Condition \ref{cond_eigenvalues} in Theorem \ref{thm_uniqueness_algebraic_structure_+}, the value of $\tau f_2$ is zero when mapping $(\alpha,\beta,\gamma,\delta)$ to $(\lambda,2,0,0)$ with $\lambda=1,-3$. Therefore, we have
$$
f_2 = (\alpha^2-\beta)/2 + (\alpha+\delta) -1/2.
$$

For $g\ge 3$, the choice of $f_g$ is in general not unique. For example, one can add $f_g$ by a multiple of $\delta^2+\beta-2$ and all the conditions on $f_g$ are still satisfied.  However, by Lemma \ref{lem_lower_order_terms_hatJ}, the quotient image of $f_g$ in $R_{g,1}$ is unique. Since
$\deg f_g = 2g < \deg \gamma^{g+1},$
Lemma \ref{lem_high_power_gamma_no_effect} implies that the quotient image of $f_g$ in 
$$
R_1:= \bC[\alpha,\beta,\gamma,\delta]/(\delta^2+\beta-2)
$$
is unique.

We have the following lemma:
\begin{Lemma}
	\label{lem_fg_inductive_formula}
For each $g\ge 3$, there exist $b_1,\dots,b_4\in \bC$ such that the following equation holds in the quotient ring $R_1$:
\begin{equation}
	\label{eqn_fg_inductive_formula}
f_{g} = \Big(\frac{1}{g}\alpha +b_1\Big)f_{g-1} + \Big(\frac{1-g}{g}\beta + b_2 \alpha + b_3 \delta+ b_4\Big)f_{g-2} - \frac{1}{2g}\gamma f_{g-3}.
\end{equation}
\end{Lemma}

\begin{proof}
We first show that \eqref{eqn_fg_inductive_formula} holds in $R_{g-2,1}$. This is weaker than the original statement because $R_{g-2,1}$ is a quotient ring of $R_1$. We will lift the identity to $R_1$ later by Lemma \ref{lem_high_power_gamma_no_effect}.
	
By \eqref{eqn_fg1_form} and \eqref{eqn_recursive_xi}, the degree of 
\begin{equation}
	\label{eqn_diff_f_g_leading}
\tau\Big(f_g - \frac{1}{g} \alpha f_{g-1} - \frac{1-g}{g}\beta f_{g-2} + \frac{1}{2g}\gamma f_{g-3}\Big)
\end{equation}
is no greater than $2(g-1)$. Since the ideals $J_{g,1}^+$ satisfy Condition \ref{cond_inclusion_hatJ_in_g} in Theorem \ref{thm_uniqueness_algebraic_structure_+}, the polynomial \eqref{eqn_diff_f_g_leading} is contained in $J_{g-2,1}^+$. Since the ideals $J_{g,1}^+$ satisfy Condition \ref{cond_leading_terms_ideal} in Theorem \ref{thm_uniqueness_algebraic_structure_+}, 
the leading-order term of \eqref{eqn_diff_f_g_leading} is a homogeneous polynomial in $I_{g-2,1}^0$ with degree at most 
$2(g-1)$. Note that this degree is strictly less than the degree of $\gamma^{g-1}$. 

Applying Theorem \ref{thm_generators_of_Ign} to $I_{g-2,1}^d$ with $d$ odd and invoking Proposition \ref{prop_singular_coh_flip_symmetry}, we see that there exist $b_1,b_2,b_3\in \bC$ and a homogeneous polynomial $h\in\bC[\alpha,\beta,\gamma,\delta]$ of degree $2(g-3)$, such that 
$$
\tau\Big(f_g - \frac{1}{g}\alpha f_{g-1} - \frac{1-g}{g}\beta f_{g-2} + \frac{1}{2g}\gamma f_{g-3}\Big) - \tau(b_1 f_{g-1} + (b_2\alpha+b_3\delta)f_{g-2}) - (\delta^2+\beta)h
$$
has degree no greater than $2(g-2)$. As a result,
$$
\tau\Big(f_g - \frac{1}{g}\alpha f_{g-1} - \frac{1-g}{g}\beta f_{g-2} + \frac{1}{2g}\gamma f_{g-3}\Big) - \tau(b_1 f_{g-1} + (b_2\alpha+b_3\delta)f_{g-2}) - (\delta^2+\beta-2)h
$$
has degree no greater than $2(g-2)$. Since it is also an element in $J_{g-2,1}^+$, it represents an element in $X_{g-2,1}$. Since $X_{g-2,1}$ is a one-dimensional linear space generated by the quotient image of $\tau f_{g-2}$, there exists $b_4$ such that 
$$
\tau\Big(f_g - \frac{1}{g}\alpha f_{g-1} - \frac{1-g}{g}\beta f_{g-2} + \frac{1}{2g}\gamma f_{g-3}\Big) - \tau(b_1 f_{g-1} + (b_2\alpha+b_3\delta)f_{g-2}) - b_4 \tau f_{g-2} 
$$
is zero in $R_{g-2,1}$. Hence Equation \eqref{eqn_fg_inductive_formula} holds in $R_{g-2,1}$.

The difference of the two sides of \eqref{eqn_fg_inductive_formula}  has degree at most $2g$. Since $g\ge 3$, this is strictly less than the degree of $\gamma^{g-1}$. Therefore, by Lemma \ref{lem_high_power_gamma_no_effect}, Equation  \eqref{eqn_fg_inductive_formula} also holds in $R_{1}$. 
\end{proof}

We now compute the coefficients in  \eqref{eqn_fg_inductive_formula}.
\begin{Lemma}
	\label{lem_b1_to_b4}
The coefficients $b_1,\dots,b_4$ in Lemma \ref{lem_fg_inductive_formula} are unique, and they are given by
$$
b_1 = (-1)^{g}\frac{2g-1}{g}, \quad b_2 = 0, \quad b_3 = (-1)^{g}\frac{2g-2}{g},  \quad b_4 = \frac{2(g-1)}{g}.
$$
\end{Lemma}

\begin{proof}
Let $V$ be the $\bC$--linear subspace of $\bC[\alpha, \beta,\gamma,\delta]$ spanned by the monomials of the form $\alpha^j\beta^k\gamma^s\delta^a$ with $a\in \{0,1\}$. Then the composition map
$$
V\hookrightarrow \bC[\alpha, \beta,\gamma,\delta] \twoheadrightarrow \bC[\alpha, \beta,\gamma,\delta]/(\delta^2+\beta-2)
$$
is a linear isomorphism (see the proof of Lemma \ref{lem_high_power_gamma_no_effect} and Lemma \ref{lem_R_gn_canonical_rep}). Hence for each $x\in  \bC[\alpha, \beta,\gamma,\delta]/(\delta^2+\beta-2)$, there exists a unique $f_x\in V$ such that $f_x$ represents $f$, and we call it the \emph{canonical representative} of $x$. 

To simplify notation, we define
$$
\xi_{k} = \xi_{k,1}, \quad \tilde{\xi}_{k} = (-1)^{k+1}\tau \xi_{k,-1},
$$
and define $\tilde{\xi}_k=0$ if $k<0$.

Recall that we assume $g\ge 3$ in Lemma \ref{lem_fg_inductive_formula}.
Comparing the $2(g-1)$ degree terms of the canonical representatives of both sides of \eqref{eqn_fg_inductive_formula}, we obtain 
\begin{equation}
	\label{eqn_subleading_comparison_Jg1}
\tilde\xi_{g-1} =  \frac{\alpha}{g}\tilde\xi_{g-2}+b_1\xi_{g-1} + \frac{1-g}{g}\beta \tilde{\xi}_{g-3} +  (b_2\alpha+b_3\delta)\xi_{g-2} - \frac{\gamma}{2g}\tilde\xi_{g-4}.
\end{equation}
in the ring $\bC[\alpha,\beta,\gamma,\delta]/(\delta^2+\beta)$. 

Similar to the case of $\bC[\alpha,\beta,\gamma,\delta]/(\delta^2+\beta-2)$, every element in  $\bC[\alpha,\beta,\gamma,\delta]/(\delta^2+\beta)$ can be uniquely represented by a linear combination of monomials of the form $\alpha^j\beta^k\gamma^s\delta^a$ with $a\in \{0,1\}$, and we call it the \emph{canonical representative} of an element in $\bC[\alpha,\beta,\gamma,\delta]/(\delta^2+\beta)$.

In the lexicographical ordering of monomials with respect to the ordering $(\alpha,\delta,\beta,\gamma)$, the canonical representative of $\xi_k$ has the form 
$$
\xi_k = \frac{1}{k!}\alpha^{k}  -\frac{(k-1)k(2k-1)}{6k!}\alpha^{k-2}\beta + \text{lower order terms},
$$
where the coefficient of $\alpha^{k-1}\beta$ is computed by the recursive formula \eqref{eqn_recursive_xi} and induction on $k$. Similarly, the canonical representative of $\xi_{k,-1}$ has the form
$$
\xi_{k,-1} = \frac{1}{k!} \alpha^k -\frac{(k-1)k(k+1)}{3k!}\alpha^{k-2}\beta + \text{lower order terms}.
$$
Therefore, when $k\ge 0$, the canonical representative of $\tilde\xi_{k} = (-1)^{k+1} \tau \xi_{k,-1}$ satisfies
\begin{align*}
(-1)^{k+1}\tilde\xi_{k} & = \frac{1}{k!}\Big(\alpha^k + k \alpha^{k-1}\delta - {k \choose 2}\alpha^{k-2}\beta  \Big)
\\
&\qquad 
- \frac{(k-1)k(k+1)}{3k!}\alpha^{k-2}\beta + \text{lower order terms}
\\
& = \frac{1}{k!}\Big(\alpha^k + k\alpha^{k-1}\delta - \frac{k(k-1)(2k+5)}{6}\alpha^{k-2}\beta\Big) + \text{lower order terms}.
\end{align*}
Comparing the coefficients of $\alpha^{g-1}$ on both sides of \eqref{eqn_subleading_comparison_Jg1}, we have
$$
(-1)^{g}\frac{1}{(g-1)!} = (-1)^{g-1}\frac{1}{g\cdot (g-2)!}+ b_1\frac{1}{(g-1)!}+ b_2\frac{1}{(g-2)!},
$$
which simplifies to
\begin{equation}
	\label{eqn_b1_b2_1}
b_1 + (g-1)b_2 = (-1)^{g}\frac{2g-1}{g}.
\end{equation}
Comparing the coefficients of $\alpha^{g-2}\delta$ yields
$$
(-1)^{g}\frac{1}{(g-2)!} = (-1)^{g-1}\frac{1}{g\cdot(g-3)!}+b_3\frac{1}{(g-2)!},
$$ 
which simplifies to
$$
b_3 = (-1)^{g}\frac{2g-2}{g}.
$$
Comparing the coefficients of $\alpha^{g-3}\beta$ yields
\begin{align*}
& -(-1)^{g}\cdot \frac{(g-1)(g-2)(2(g-1)+5)}{6(g-1)!} 
\\
= & -\frac{1}{g} (-1)^{g-1}\frac{(g-2)(g-3)(2(g-2)+5)}{6(g-2)!}- b_1 \frac{(g-1)(g-2)(2(g-1)-1)}{6(g-1)!}\\
 &\quad +(-1)^{g-2}\frac{1-g}{g}\frac{1}{(g-3)!} - b_2\frac{(g-2)(g-3)(2(g-2)-1)}{6(g-2)!}
\end{align*}
which simplifies to 
\begin{equation}
	\label{eqn_b1_b2_2}
 (2g-3)b_1 + (g-3)(2g-5)b_2 = (-1)^{g} \frac{(2g-1)(2g-3)}{g}.
\end{equation}
When $g\ge 3$, Equations \eqref{eqn_b1_b2_1}, \eqref{eqn_b1_b2_2} have a unique solution
$$
b_1 = (-1)^{g}\frac{2g-1}{g}, \quad b_2 = 0.
$$

It remains to determine the value of $b_4$. Consider the homomorphism from $R_1$ to $\bC$ that maps
$(\alpha,\beta,\gamma,\delta)$ to $((-1)^{g}(2g-3), 2, 0,0)$.
Then by Condition \ref{cond_eigenvalues}, the images of $f_g$, $f_{g-1}$ are zero, and the image of $f_{g-2}$ is non-zero. Therefore, \eqref{eqn_fg_inductive_formula} implies that the image of $\frac{1-g}{g}\beta + b_2 \alpha + b_3 \delta+ b_4$ is zero, which implies 
$$
b_4 = \frac{2g-2}{g}.
\phantom\qedhere\makeatletter\displaymath@qed
$$
\end{proof}

We summarize the above results in terms of $\omega$, $\delta$, $\beta$, $\gamma$. Recall that $\omega = \alpha+\delta/2$.
\begin{Theorem}
	\label{thm_generators_Jg1}
	Let $r_g\in \bC[\omega,\beta,\gamma,\delta]$ ($g\ge 0$) be given by
	$$
	r_0 = 1, \quad r_1 = \omega+\delta/2-1,\quad r_2 = ((\omega+\delta/2)^2-\beta)/2+\omega-\delta/2-1/2,
	$$
	and when $g\ge 3$,
$$
	gr_{g} = \Big(\omega+\delta/2 +(-1)^{g}(2g-1)\Big)r_{g-1} + (1-g)\Big(\beta  +(-1)^{g} (2 \delta)-2\Big)r_{g-2} - \frac{\gamma}{2} r_{g-3}.
$$
Then we have
$$
J_{g,1}^+ = (r_{g},\,r_{g+1},\, r_{g+2},\, \delta^2+\beta-2).
$$
\end{Theorem}

\begin{proof}
	By the previous results, we have $r_g = \tau f_g$ for all $g$, and
	$$J_{g,1}^+ = (r_{g},\,r_{g+1},\,r_{g+2},\, \gamma^{g+1},\, \delta^2+\beta-2).$$
	The recursive relation on $r_g$ implies $\gamma (r_{g},\,r_{g+1},\,r_{g+2})\subset (r_{g+1},\,r_{g+2},\, r_{g+3})$. By induction, this implies $\gamma^g\in (r_{g},\,r_{g+1},\,r_{g+2})$. Hence the result is proved.
\end{proof}

The recursive relation of $f_g$ implies that
\begin{equation}
	\label{eqn_Jg1_multiply_beta_inclusion}
(\beta + (-1)^{g} (2\delta) - 2) J_{g,1}^+ \subset (J_{g+1,1}^+,\gamma)
\end{equation}
for all $g\ge 1$. It is straightforward to verify that \eqref{eqn_Jg1_multiply_beta_inclusion} also holds for $g=0$.

\begin{Theorem}
	\label{thm_top_generalized_eigenspace_Vg1+}
	Assume $g\ge 1$. The simultaneous generalized eigenspace for the $\mu$ maps $(\muu(\Sigma),\mu(\pt))$ on
	 $\bV_{g,1}^+$ with eigenvalues $((-1)^{g+1}(2g-1),2)$ has dimension one. 
\end{Theorem}
\begin{proof}
As before, we will use $(\alpha,\beta,\gamma,\delta)$ to denote the linear operators on $\bV_{g,1}^+$ defined by multiplications of $(\alpha,\beta,\gamma,\delta)$. 

Consider the homomorphism $\bC[\alpha,\beta,\gamma,\delta]\to \bC[\alpha]$ mapping $(\alpha,\beta,\gamma,\delta)$ to $(\alpha,2,0,0)$. The image of $\tau f_g$ is a polynomial in $\alpha$ with leading term $\alpha^g/g!$. By Condition \ref{cond_eigenvalues} of Theorem \ref{thm_uniqueness_algebraic_structure_+}, we know that the image of $\tau f_g$ has roots $\alpha = (-1)^{i+1}(2i-1)$ for $i=1,\dots,g$. Therefore it must be a scalar multiple of $\prod_{i=1}^g(\alpha-(-1)^{i+1}(2i-1))$. Hence 
\begin{equation}
	\label{eqn_defn_Pg}
P_g(\alpha):= \prod_{i=1}^g (\alpha-(-1)^{i+1}(2i-1))\in (J_{g,1}^+,\beta-2,\gamma,\delta).
\end{equation}

Since $\bV_{g,1}^+\cong \bC[\alpha,\beta,\gamma,\delta]/J_{g,1}^+$ is a finite-dimensional linear space, for each $g$, there exists a $k_g\in \mathbb{N}$ such that 
$(\beta-2)^{k_g}$ and $(\beta-2)^{k_g+1}$ have the same image as linear operators on $\bV_{g,1}^+$. Then we have
\begin{equation}
	\label{eqn_beta-2_high_power_stabilize}
(J_{g,1}^+,(\beta-2)^{k_g})  = (J_{g,1}^+,(\beta-2)^{k_g+1}) = (J_{g,1}^+,(\beta-2)^{k_g+2}) = \cdots,
\end{equation}
and $\bC[\alpha,\beta,\gamma,\delta]/(J_{g,1}^+,(\beta-2)^{k_g}) $ is isomorphic to the generalized eigenspace of $\beta$ with eigenvalue $2$ in $\bV_{g,1}^+$ as a $\bC[\alpha,\beta,\gamma,\delta]$--module.

By Condition \ref{cond_eigenvalues} in Theorem \ref{thm_uniqueness_algebraic_structure_+}, the eigenvalues of $\alpha$ on 
$$\bC[\alpha,\beta,\gamma,\delta]/(J_{g,1}^+,(\beta-2)^{k_g})$$ are given by 
$(-1)^{i+1}(2i-1)$ for $i=1,\dots, g$. So there exists a positive integer $N_{g}$ depending on $g$, such that
\begin{equation}
	\label{eqn_defn_Ng_on_exp_Pg}
P_g(\alpha)^{N_g} \in (J_{g,1}^+,(\beta-2)^{k_g}).
\end{equation}
By \eqref{eqn_beta-2_high_power_stabilize} and \eqref{eqn_Jg1_multiply_beta_inclusion}, we have 
\begin{equation}
	\label{eqn_Pg_power_times_beta_in_Jg+1}
P_{g}(\alpha)^{N_{g}}\cdot (\beta+(-1)^{g}(2\delta)-2)\in (J_{g+1,1}^+,(\beta-2)^{k_{g+1}}).
\end{equation}
Since $\beta+\delta^2-2\in J_{g+1,1}^+$, we have 
$$
\beta+(-1)^{g}(2\delta)-2= (-1)^{g}(2\delta)-\delta^2 = \delta((-1)^{g} 2-\delta)
$$
as operators on $\bV_{g+1,1}^+$. Moreover, on the generalized eigenspace of $\beta$ with eigenvalue $2$, the only eigenvalue of $\delta$ is $0$, so $(-1)^{g} 2-\delta$ is an invertible linear operator. By \eqref{eqn_Pg_power_times_beta_in_Jg+1}, the linear operator defined by the left-hand side of \eqref{eqn_Pg_power_times_beta_in_Jg+1} is the zero operator on the generalized eigenspace for $\beta$ with eigenvalue $2$ in $\bV_{g+1,1}^+$. Therefore the same statement holds for $P_{g}(\alpha)^{N_{g}}\cdot \delta$. Hence we have
\begin{equation}
	\label{eqn_Pg_Ng_delta}
P_g(\alpha)^{N_g} \cdot \delta \in  (J_{g+1,1}^+,(\beta-2)^{k_{g+1}}).
\end{equation}
Since $\beta+\delta^2-2\in J_{g+1,1}^+$, this implies
\begin{equation}
	\label{eqn_Pg_Ng_beta}
P_g(\alpha)^{N_g} \cdot (\beta-2) \in  (J_{g+1,1}^+,(\beta-2)^{k_{g+1}}).
\end{equation}
Since $\gamma J_{g,1}^+\subset J_{g+1,1}^+$, by \eqref{eqn_beta-2_high_power_stabilize} and \eqref{eqn_defn_Ng_on_exp_Pg}, we also have
\begin{equation}
	\label{eqn_Pg_Ng_gamma}
P_g(\alpha)^{N_g} \cdot\gamma \in (J_{g+1,1}^+,(\beta-2)^{k_{g+1}}).
\end{equation}
Combining \eqref{eqn_Pg_Ng_delta}, \eqref{eqn_Pg_Ng_beta}, \eqref{eqn_Pg_Ng_gamma}, with \eqref{eqn_defn_Pg}, we have 
\begin{align*}
P_g(\alpha)^{N_g} \cdot P_{g+1}(\alpha)& \in (P_g(\alpha)^{N_g}\cdot J_{g+1,1}^+,P_g(\alpha)^{N_g}\cdot(\beta-2),P_g(\alpha)^{N_g}\cdot\gamma,P_g(\alpha)^{N_g}\cdot\delta)
\\
& \subset (J_{g+1,1}^+,(\beta-2)^{k_{g+1}}).
\end{align*}
Therefore, in the linear space $\bC[\alpha,\beta,\gamma,\delta]/(J_{g+1,1}^+,(\beta-2)^{k_{g+1}})$, the \emph{generalized eigenspace} of $\alpha$ with eigenvalue $(-1)^{g}(2g+1)$ is equal to the \emph{eigenspace} of $\alpha$ with eigenvalue $(-1)^{g}(2g+1)$.

The eigenspace of $\alpha$ with eigenvalue $(-1)^{g}(2g+1)$ in $$\bC[\alpha,\beta,\gamma,\delta]/(J_{g+1,1}^+,(\beta-2)^{k_{g+1}})$$ is isomorphic to 
$$
\bC[\alpha,\beta,\gamma,\delta]/\Big(J_{g+1,1}^+,(\beta-2)^{k_{g+1}},(\alpha-(-1)^{g}(2g+1))\Big)
$$
as a $\bC[\alpha,\beta,\gamma,\delta]$--module. Since $P_g(\alpha)^{N_g}$ is invertible in $$\bC[\alpha,\beta,\gamma,\delta]/(\alpha -(-1)^{g}(2g+1) ),$$
by \eqref{eqn_Pg_Ng_delta}, \eqref{eqn_Pg_Ng_beta}, \eqref{eqn_Pg_Ng_gamma}, we have
$$
\beta-2,\delta,\gamma \in \Big(J_{g+1,1}^+,(\beta-2)^{k_{g+1}},(\alpha-(-1)^{g}(2g+1))\Big),
$$
therefore 
\begin{align*}
& \dim_\bC \bC[\alpha,\beta,\gamma,\delta]/\Big(J_{g+1,1}^+,(\beta-2)^{k_{g+1}},(\alpha-(-1)^{g}(2g+1))\Big)
\\
\le & \dim_\bC \bC[\alpha,\beta,\gamma,\delta]/(\alpha -(-1)^{g}(2g+1), \beta-2,\gamma,\delta)=1.
\end{align*}
Hence the generalized eigenspace for $\alpha$ with eigenvalue $(-1)^{g}(2g+1)$ in $$\bC[\alpha,\beta,\gamma,\delta]/(J_{g+1,1}^+,(\beta-2)^{k_{g+1}})$$ is at most $1$. By Condition \ref{cond_eigenvalues} of Theorem \ref{thm_uniqueness_algebraic_structure_+}, we know that the generalized eigenspace is non-zero. Hence its dimension is $1$. This proved Theorem \ref{thm_top_generalized_eigenspace_Vg1+} when $g\ge 2$. The case when $g=1$ can be verified directly using Theorem \ref{thm_generators_Jg1}.
\end{proof}

The analogous result also holds for $\bV_{g,1}^-$.
\begin{Theorem}
	\label{thm_top_generalized_eigenspace_Vg1-}
	Assume $g\ge 1$. The simultaneous generalized eigenspace of the $\mu$ maps $(\muu(\Sigma),\mu(\pt))$ on
	$\bV_{g,1}^-$ with eigenvalues $((-1)^{g}(2g-1),2)$ has dimension one. 
\end{Theorem}

\begin{proof}
The theorem can be proved almost verbatim as the proof of Theorem \ref{thm_top_generalized_eigenspace_Vg1+}. Alternatively, it also follows from the symmetry of generalized eigenspaces for degree 2 maps on instanton homology (see, for example, \cite{xie2020instantons}*{Lemma 3.14}).
\end{proof}

Theorems \ref{thm_top_generalized_eigenspace_Vg1+} and \ref{thm_top_generalized_eigenspace_Vg1-} imply the following excision formula for singular instanton Floer homology.

\begin{Corollary}
	The singular excision formula stated in \cite{XZ:excision}*{Theorem 6.4} also holds for $n=1$.\qed 
\end{Corollary}

The proof is verbatim as the original proof of  \cite{XZ:excision}*{Theorem 6.4}, except that we now use Theorems \ref{thm_top_generalized_eigenspace_Vg1+} and \ref{thm_top_generalized_eigenspace_Vg1-} in the place of \cite{XZ:excision}*{Corollary 5.16}.

\section{Some discussions about local coefficients}
\label{sec_local_coef}

In this section, we use  $\bV_{g,n}(u)$ to denote the instanton homology of $(S^1\times \Sigma_g,S^1\times\{p_1,\dots,p_n\},\emptyset)$ with local coefficients, where we use $u$ to denote the local coefficient variable (our $u$ corresponds to $\tau$ in \cite{kronheimer2022relations}). 
Most of the structural results and constructions
for $\bV_{g,n}$ in Section \ref{sec_preliminaries} still hold for $\bV_{g,n}(u)$ after replacing the coefficient ring
$\mathbb{C}$ by  $\mathcal{R}:=\bC[u,u^{-1}]$ (also cf. \cite{kronheimer2022relations}*{Section 3.4}).
The main difference is that we no longer have flip symmetries on $\bV_{g,n}(u)$ as $\mathcal{R}$--homomorphisms.

 Let
 $\mathbb{A}_{g,n}(u) := \mathbb{A}_{g,n}\otimes_{\bC} \mathcal{R}$.
We summarize the basic properties of $\bV_{g,n}(u)$ in the following proposition.

\begin{Proposition}
The instanton homology ring $\mathbb{V}_{g,n}(u)$ satisfies the following:	
\begin{enumerate}
\item There exits a surjective ring homomorphism
$$
\Phi: \mathbb{A}_{g,n}(u)[\epsilon]/(\epsilon^2-1)\twoheadrightarrow \bV_{g,n}(u).
$$

\item Let $\bV_{g,n}^+(u)$ (or $\bV_{g,n}^-(u)$) be the eigenspace for $\Phi(\epsilon)$ with eigenvalue $1$ (or $-1$). Then 
$\bV_{g,n}^\pm(u)$ is isomorphic to $H^\ast(N_{g,n}^0;\bC)\otimes_{\bC} \mathcal{R}$ as (free) $\mathcal{R}$-modules.

\item We still use $\Phi$ to denote the composition of the projection
\begin{align*}
& \bigoplus_{k=0}^g \Lambda_0^k H_1(\Sigma;\bC)\otimes_{\bC} \mathcal{R}[\alpha, \beta, \gamma, \delta_1,\dots, \delta_n,\epsilon]/(\epsilon^2-1)
\nonumber
\\
&\quad \twoheadrightarrow  \bigoplus_{k=0}^g \Lambda_0^k H_1(\Sigma;\bC)\otimes \mathcal{R}[\alpha, \beta, \gamma, \delta_1,\dots, \delta_n,\epsilon]/(\epsilon^2-1)/(\gamma^{g-k+1}) 
=\bA_{g,n}(u)
\end{align*}
with
$\Phi$. Then its kernel has the form
$$
	\ker \Phi = \bigoplus_{k=0}^g \Lambda_0^k  H_1(\Sigma;\bC) \otimes_{\bC} J_{g,n,k}(u),
$$
where $J_{g,n,k}(u)$ are ideals in $ \mathcal{R}[\alpha, \beta, \gamma, \delta_1,\dots, \delta_n,\epsilon]/(\epsilon^2-1)$. Moreover,
$J_{g,n,k}(u) = J_{g-k,n,0}(u)$.

\item $\delta_i^2 + \beta - u^2-u^{-2}\in  J_{g,n,k}(u)$.

\item Denote $J_{g-k,n,0}(u)$ by $J_{g-k,n}(u)$. Then there exists
 $$W^\eta=w^\eta_0+(-1)^{g}\hat{\epsilon} u^{n-2|\eta|} w_1^\eta +w_2^\eta + \cdots \in J_{g,n}(u)$$ where 
 $\hat \epsilon$ is defined in Section \ref{subsec_sign_E},
$w_0^\eta$, $w_1^\eta$ are defined in Section \ref{subsec_construct_w0w1} and $ w_i^\eta\in 
\mathcal{R}[\alpha,\beta, \gamma,\delta_1,\cdots,\delta_n,\epsilon]/(\epsilon^2-1)$ is a homogeneous polynomial of degree $2g+2m-2i$. Here, we assign $u$ with degree zero.

\end{enumerate}

\end{Proposition}
Items (1), (2), (3), (4) are straightforward generalizations of the corresponding results on $\bV_{g,n}$. Item (5) is a generalization of of Corollary \ref{cor_subleading_term_formula}.
There is no essential difference between the proof of this result and Corollary \ref{cor_subleading_term_formula} except for
determining the factor $u^{n-2|\eta|}$ in front of $w_1^\eta$. In the genus $0$ case, this was done in 
\cite{kronheimer2022relations}*{Lemma 5.10}. There is no change in the proof for the general genus case, so we skip the proof.

The rest of this section extends the computations in Section \ref{sec_n=1} and finds an explicit formula for the generators of $J_{g,1}^+(u)$. When $n=1$, we denote $\delta_1$ by $\delta$, and denote the $\mu$--map $\sigma_p$ for $p\in S^1\times \{p_1\}$ by $\sigma$.

For $\theta\in \bC$, let $\bV_{g,n}(\theta)$ be the instanton homology of $(S^1\times \Sigma_g,S^1\times\{p_1,\dots,p_n\},\emptyset)$ with local coefficients, where the local coefficient variable is taken to be $\theta$. 
By definition, $\bV_{g,n}(\theta)$ is a $\bC$--algebra.

\begin{Lemma}
	\label{lem_eigenvalues_Vgn_theta}
For all $\lambda = 2- 2g +2k$ with $k\in\{0,1,\dots, g-1\}$ and all $\epsilon \in \{-1,1\}$, there exists $0\neq v \in \bV_{g,1}(\theta)$ such that $v$ is a simultaneous eigenvector of $(\muu(\Sigma)),\sigma, \mu(\pt))$ with eigenvalues 
$$
\Big(\epsilon \frac{\theta+\theta^{-1}}{2}+\lambda, \epsilon(\theta^{-1}-\theta), 2\Big).
$$
\end{Lemma}

\begin{proof}
Let $(Y_g,L_{g,n}) = S^1\times (\Sigma_g,\{p_1,\dots,p_n\})$, where $\Sigma_g$ has genus $g$.  Fix a point $x_0\in S^1$.
If $g\ge 1$, let $c_1\subset \Sigma_{g}\backslash\{p_1,\cdots,p_n\}$ be a non-separating simple closed curve. Let $\omega$ be the circle 
$\{x_0\} \times c_1$ in $(Y_{g},L_{g,n})$ and let $\bW_{g,n}(\theta)$ be the instanton Floer homology of $(Y_{g},L_{g,n},\omega)$ with local coefficients associated with all the components of $L_{g,n}$, and the local coefficient variable is taken to be $\theta$. 

If $n\ge 2$, let $c_2\subset \Sigma_{g}$ be an arc that connects $p_1$ and $p_2$.
Let $\omega'$ be the arc $ \{x_0\}\times c_2$ in $(Y_{g},L_{g,n})$ and let 
$\bU_{g,n}(\theta)$ be the instanton Floer homology of  $(Y_{g},L_{g,n},\omega')$ with local coefficients associated with $S^1\times\{p_3,\dots,p_n\}$, and the local coefficient variable is taken to be $\theta$. By \cite{XZ:excision}*{Proposition 2.1}, the action of $\mu(\pt)$ on $\bU_{g,n}(\theta)$ is the scalar multiplication by $2$.

It was proved in \cite{XZ:excision}*{Lemma 5.2}\footnote{The original argument was given for $\theta=1$, but the same argument works for general $\theta$.} that all eigenvalues of $\bW_{g,n}(\theta)$ are also eigenvalues of $\bV_{g,n}(\theta)$.
By torus excision, we have
$$
\bW_{1,1}^{(2)}(\theta) \otimes_\bC \bW_{g,0}^{(2)}(\theta) \cong \bW_{g,1}^{(2)}(\theta),
$$
where $\bW_{g,n}^{(2)}$ denotes the generalized eigenspace of $\mu(\pt)$ in $\bW_{g,n}$ with eigenvalue $2$. 
The above isomorphism also intertwines $\muu(\Sigma_1)+\muu(\Sigma_g)$ on the left-hand side with $\muu(\Sigma_g)$ on the right-hand side. 
Note that when $n=0$, the definition of $\bW_{g,n}$ does not depend on the local coefficient system, therefore the space $\bW_{g,0}(\theta)$ is the same as $\bW_{g,0}(1)$, which was studied in \cite{XZ:excision}*{Lemma 5.1}. In particular, we know that every $\lambda = 2- 2g +2k$ with $k\in\{0,1,\dots, g-1\}$ is an eigenvalue for $\muu(\Sigma)$ in $\bW_{g,0}^{(2)}(\theta)$.

So we only need to show that on $\bW_{1,1}^{(2)}(\theta)$, the operators $(\muu(\Sigma_1),\sigma)$ have simultaneous eigenvalues $(\epsilon (\theta+\theta^{-1})/2,\epsilon(\theta^{-1}-\theta),2 )$ for $\epsilon \in \{1,-1\}$. By torus excision again, the eigenvalues for $(\muu(\Sigma_1),\sigma)$ on $\bW_{1,1}^{(2)}(\theta)$ are the same as the eigenvalues for $(\muu(\Sigma_0),\sigma_{p_3})$ on $\bU_{0,3}(\theta)$, where $p_3\in S^2$ is the marked point corresponding to the link component with local coefficients in the definition of $\bU_{0,3}(\theta)$.

The critical points of the perturbed Chern--Simons functional defining $\bU_{0,3}(\theta)$ contains two non-degenerate orbits with degrees differ by $2$ (see \cite{AHI}*{Example 4.2}), therefore $\bU_{0,3}(\theta)$ has dimension $2$. In particular,   $\bU_{0,3}(\theta)\neq 0$. The simultaneous eigenvalues of $(\muu(\Sigma_0),\sigma)$  are symmetric with respect to an overall change of signs (see, for example, \cite{xie2020instantons}*{Lemma 3.14}).   The proof of \cite{XZ:excision}*{Lemma 5.2} showed that every eigenvalue for $(\muu(\Sigma_0),\sigma)$  on $\bU_{0,3}(\theta)$ is an eigenvalue for $(\muu(\Sigma_0),\sigma)$  on $\II(S^1\times S^2, S^1\times \{p_1,p_2,p_3\},\emptyset)$, where the local coefficient system is associated with $S^1\times\{p_3\}$. 
The latter homology was computed in  \cite{kronheimer2022relations}*{Section 7.4}, and it can be seen from the computation that the operators $(\muu(\Sigma_1),\sigma)$ have simultaneous eigenvalues $(\epsilon (\theta+\theta^{-1})/2,\epsilon(\theta^{-1}-\theta),2 )$  for $\epsilon \in \{1,-1\}$. Since $\bU_{0,3}(\theta)\neq 0$ and the eigenvalues for  $(\muu(\Sigma_0),\sigma)$  are symmetric with respect to sign reversals, the desired result follows.
\end{proof}

\begin{Lemma}
	\label{lem_top_eigenvalue_theta}
	Suppose $g\ge 1$.
	There exists an open neighborhood $U$ of $1\in \bC$ such that for all $\theta\in U$, we have 
\begin{equation}
	\label{eqn_eigenvalues_theta}
\Big( (-1)^{g+1}(2g-2+\frac{\theta+\theta^{-1}}{2}), (-1)^{g+1}(\theta^{-1}-\theta), 2, 0\Big)
\end{equation}
are simultaneous eigenvalues for $(\omega,\delta,\beta,\gamma)$ on $\bV_{g,1}^+(\theta)$ but not on $\bV_{g-1,1}^+(\theta)$. 
\end{Lemma}

\begin{proof}
By Lemma \ref{lem_eigenvalues_Vgn_theta} and the fact that $\gamma$ is nilpotent, we know that \eqref{eqn_eigenvalues_theta} are simultaneous eigenvalues for $(\omega,\delta,\beta,\gamma)$ on  $\bV_{g,1}(\theta)=\bV_{g,n}^+(\theta)\oplus \bV_{g,n}^-(\theta)$. 
When $\theta=1$, by Condition \ref{cond_eigenvalues} in Theorem \ref{thm_uniqueness_algebraic_structure_+}, \eqref{eqn_eigenvalues_theta} are not simultaneous eigenvalues for $(\omega,\delta,\beta,\gamma)$ on $\bV_{g,1}(\theta)^-$ or $\bV_{g-1,1}^+(\theta)$. 
Since being an eigenvalue is a closed condition on $\theta$, the result follows. 
\end{proof}

Now we write down a formula for the instanton homology of $(S^1\times \Sigma,S^1\times \{p\})$ with local coefficients. 
\begin{Theorem}
	Let $r_g\in \bC[u,u^{-1}, \omega,\beta,\gamma,\delta]$ ($g\ge 0$) be given by
	$$
	r_0 = 1, \quad r_1 = \omega+\delta/2-u^{-1},\quad r_2 = ((\omega+\delta/2)^2-\beta)/2+u^{-1}(\omega-\delta/2)-u^{-2}/2,
	$$
	and when $g\ge 3$,
	$$
	gr_{g} = \Big(\omega+\delta/2 +(-1)^{g}(2g-1)u^{-1}\Big)r_{g-1} + (1-g)\Big(\beta  +(-1)^{g} (2 \delta)u^{-1}-2u^{-2}\Big)r_{g-2} - \frac{\gamma}{2} r_{g-3}.
	$$
	Then we have
	$$
	J_{g,1}^+(u) = (r_{g},\,r_{g+1},\, r_{g+2},\, \delta^2+\beta-u^2-u^{-2}).
	$$
\end{Theorem}

\begin{proof}The proof is almost the same as the proof of Theorem \ref{thm_generators_Jg1}, and we give a sketch of the argument.
	 By the same argument as Lemma \ref{lem_Xgn_generates_hatJ}, there exists a sequence 
	 $\{r_k\}_{k\ge 0}$ in $\bC[u,u^{-1},\alpha,\beta,\gamma,\delta]$, such that  $r_0=1$, and for $k\ge 1$,
	 $$
	 r_k = \tau \xi_{k,1} + (-1)^g u^{-1} \xi_{k,-1} + (\text{terms with degree} \le 2(k-2)),
	 $$
	 so that $J_{g,1}^+(u) = (r_g,r_{g+1},r_{g+2}, \delta^2+\beta-u^2-u^{-2}, \gamma^{g+1})$. 
	 
	 Moreover, the same argument as  Lemma \ref{lem_fg_inductive_formula} shows that there exist $b_1,\dots,b_4\in \bC[u,u^{-1}]$ such that
	 $$
	 g r_{g} =  \Big(\omega+\delta/2 +b_1\Big)r_{g-1} + \Big((1-g)\beta  +b_2\omega + b_3\delta + b_4\Big)r_{g-2} - \frac{\gamma}{2} r_{g-3}.
	 $$
	 for all $g\ge 3$, where the above equation holds modulo $(\beta+\delta^2-u^2-u^{-2})$.
	The values of $b_1,b_2,b_3$ are then determined by the comparison of the sub-leading terms, as in the proof of Lemma \ref{lem_b1_to_b4}. To determine the value of $b_4$, note that by Lemma \ref{lem_top_eigenvalue_theta}, there exists an open neighborhood $U$ of $1\in \bC$ such that for all $\theta \in U$, the homomorphism from 
	$\bC[u, u^{-1},\omega,\beta,\gamma,\delta]$ to $\bC$ taking $(\omega,\delta,\beta,\gamma,u)$ to 
	$$
	\Big( (-1)^{g}(2g-4+\frac{\theta+\theta^{-1}}{2}), (-1)^{g}(\theta^{-1}-\theta), 2, 0,\theta\Big)
	$$
	maps $r_{g-1}$, $r_g$, and $\gamma r_{g-3}$ to $0$, and maps $r_{g-2}$ to a non-zero number. Therefore, the image of $\beta+b_3\delta + b_4$ must be zero, which yields the value of $b_4$. The recursive relation of $r_k$ then implies 
	$$
	\gamma^g \in (r_g,r_{g+1},r_{g+2}),
	$$
	so $J_{g,1}^+(u) = (r_g,r_{g+1},r_{g+2}, \delta^2+\beta-u^2-u^{-2})$. 
\end{proof}

\bibliographystyle{amsalpha}
\bibliography{references}

\end{document}